\documentclass{tac}

\usepackage{graphicx}
\usepackage{amssymb} 
\usepackage{enumerate}
\usepackage{amsmath}
\usepackage{mathtools}
\usepackage{latexsym, amssymb}
\usepackage{setspace}
\usepackage{titling}
\usepackage[mathscr]{euscript}
\usepackage{etoolbox}
\usepackage{listings}
\usepackage{hyperref}
\usepackage{mathtools}
\usepackage{multicol}
\usepackage{tikz} 
\usepackage{tikz-cd}
\usetikzlibrary{decorations.pathreplacing}
\usepackage{bm}
\usepackage{extpfeil}
\usepackage{xcolor}
\usepackage[utf8]{inputenc}
\usepackage{enumitem}
\usepackage{thmtools}

\makeatletter
\renewcommand\thmt@autorefsetup{%
  \@xa\def\csname\thmt@envname autorefname\@xa\endcsname\@xa{\thmt@thmname}%
  \ifthmt@hassibling
    \@xa\def\@xa\thmt@autoreffix\@xa{%
      \@xa\global\@xa\let\csname the\thmt@envname\@xa\endcsname
        \csname the\thmt@sibling\endcsname
      \def\thmt@autoreffix{}%
    }%
  \fi
}
\makeatother
\RequirePackage{hyperref}
\RequirePackage{cleveref}

\usetikzlibrary{decorations.markings}

\tikzcdset{scale cd/.style={every label/.append style={scale=#1},
    cells={nodes={scale=#1}}}}

\makeatletter
\newcommand{\colim@}[2]{%
  \vtop{\m@th\ialign{##\cr
    \hfil$#1\operator@font colim$\hfil\cr
    \noalign{\nointerlineskip\kern-\ex@}\cr}}%
}
\newcommand{\colim}{%
  \mathop{\mathpalette\colim@{\rightarrowfill@\textstyle}}\nmlimits@
}
\makeatother

\newtheoremrm{rmk}{Remark}
\newtheoremrm{construction}{Construction}

\newcommand{\f}[1]{\mathfrak{#1}}
\newcommand{\s}[1]{\mathscr{#1}}
\renewcommand{\c}[1]{\mathcal{#1}}

\renewcommand{\bf}[1]{\mathbf{#1}}
\newcommand{\maybehat}[1]{#1} 
\newcommand{\ol}[1]{\overline{#1}}

\newcommand{\tx}[1]{\textnormal{#1}}
\newcommand{\ul}[1]{\underline{#1}}

\newcommand{\db}{\mathbf{\Delta}}
\newcommand{\dop}{\mathbf{\Delta}^{\textrm{op}}}

\newcommand{\cat}{\s{C}\tx{at}}
\newcommand{\prl}{\s{P}\tx{r}^\tx{L}}

\DeclareRobustCommand{\atez}{\text{\reflectbox{$\zeta$}}}

\tikzset{%
    symbol/.style={%
        draw=none,
        every to/.append style={%
            edge node={node [allow upside down, auto=false]{$#1$}}}
    }
}

\tikzcdset{scale cd/.style={every label/.append style={scale=#1},
    cells={nodes={scale=#1}}}}
    
\title{Higher zeroes in pointed infinity-categories}
\author{Thomas Rush Brown}
\date{\vspace{-5ex}} 

\begin{document}
\maketitle
 \begin{abstract} We construct an endofunctor on (small) simplicial sets which freely adjoins a zero object to a (small) infinity-category, as the co-cone construction freely adjoins a terminal object.  We then prove that this endofunctor induces a left adjoint to the inclusion of the infinity-category of (small) pointed infinity-categories and morphisms that preserve zero objects into the infinity-category of (small) infinity-categories.  This adjunction, in turn, leads to the definition of a \emph{coherent choice of nullhomotopies} on a pointed infinity-category.  A coherent choice of nullhomotopies is a non-canonical, coherent choice of ``higher zeroes'' on a pointed infinity-category, mirroring the obvious definition of higher zeroes on the differential graded nerve of a differential graded category of chain complexes over an additive category.  We then show that coherent choices of nullhomotopies give an alternate characterization of fiber sequences in pointed infinity-categories, generalizing the 1-categorical definition of kernels in terms of morphisms that compose to zero.
  \end{abstract}
\section{Introduction}
Joyal's formulation of $\infty$-categories allows us to handle homotopy limits as a generalization of $1$-categorical limits by using easy simplicial constructions that generalize the $1$-categorical notions of cones and slice categories.  The generalization is not quite as direct for fiber sequences, where we lack a simplicial notion of ``higher zero morphisms'' in the context of pointed $\infty$-categories.  If $\s{C}$ is a stable $\infty$-category, we know what it means for a $1$-simplex of $\s{C}$ to be a zero morphism (i.e., it gives a zero morphism in the triangulated category $\tx{h}\s{C}$), but only knowing the zero morphisms of $\tx{h}\s{C}$ is too little information to compute fiber sequences in $\s{C}$.  We are thus at a disadvantage in any context\footnote{Such as in the $\infty$-operadic constructions of (\cite{thesis}, \S 2).} where it would be convenient to mimic the $1$-categorical definition of a kernel of $f: x\to y$ as a terminal object of the full subcategory $\{ g\in \s{C}_{/x}\tx{ ; } f\circ g=0 \} \subseteq \s{C}_{/x}$, rather than as a pullback of a morphism $0\to y$ along $f$.  What's missing here is an easy construction on simplicial sets that freely adjoins a zero object to an $\infty$-category, as the co-cone construction freely adjoins a terminal object.

In Definition~\ref{diamonddefn} we give an endofunctor on simplicial sets, denoted $(-)_\diamond : \tx{sSet} \to \tx{sSet}$, which is an analogue of the functor $(-)^\triangleright : \tx{sSet} \to \tx{sSet}$.  Where the functor $(-)^\triangleright$ freely adjoins a terminal object to a small $\infty$-category, the functor $(-)_\diamond$ freely adjoins a zero object to a small $\infty$-category while keeping it small.  Note that although the $(-)_\diamond$ construction freely adjoins a zero object, it should not be confused with the construction $\s{C} \mapsto \s{C}_* \simeq \s{C} \otimes \s{S}_*$ (\cite{ha}, 1.4.2.18 and 4.8.1.21) of the free pointed $\infty$-category on a presentable $\infty$-category $\s{C}$.  The functor $-\otimes \s{S_*} : \prl \to \prl_{\tx{pt}}$ is left adjoint to the inclusion $\prl_{\tx{pt}} \hookrightarrow \prl$ of the full subcategory in $\prl$ spanned by the presentable $\infty$-categories which contain a zero object (\cite{univloopmachine}, \S 4).  However, the inclusion $\tx{Cat}_\infty^\tx{pt} \hookrightarrow \tx{Cat}_\infty$ is not full, and its left adjoint is not given by a version of $(-)_*$. 

In Theorem~\ref{diamondisfreezerothm} we show, via the classification of coCartesian fibrations, that the induced functor of $\infty$-categories $(-)_\diamond : \tx{Cat}_\infty \to \tx{Cat}_\infty^\tx{pt}$ gives a left adjoint to the inclusion $\tx{Cat}_\infty^\tx{pt} \hookrightarrow \tx{Cat}_\infty$, where $\tx{Cat}_\infty^\tx{pt} \subseteq \tx{Cat}_\infty$ is the subcategory of (small) pointed $\infty$-categories with functors that preserve zero objects (Definition~\ref{catinftyptdefn}).  This adjunction, in turn, shows that an object $x$ of an $\infty$-category $\s{C}$ is a zero object if and only if the canonical inclusion $\s{C}\hookrightarrow \s{C}_\diamond$ admits a left inverse $\f{z} : \s{C}_\diamond \to \s{C}$ whose restriction to $\{x \}_\diamond$ is the constant map to $x$ (Proposition~\ref{sectiondiamondzeroobj}).  We call such a map $\f{z} : \s{C}_\diamond \to \s{C}$ a \emph{coherent choice of nullhomotopies} (Definition~\ref{cohnulldefn}) because $\f{z}$ gives a non-canonical, coherent choice of ``higher zeroes'' on $\s{C}$, generalizing the obvious definition of ``higher zeroes'' in the case of the differential graded nerve of a differential graded category of chain complexes; see Remark~\ref{dgnervermk}.  

In the final section, we apply coherent choices of nullhomotopies to give an alternate description of fiber sequences in pointed $\infty$-categories, with this description mimicking the $1$-categorical definition of a kernel of $f: x\to y$ (a morphism in $1$-category $\s{C}$ enriched over pointed sets) as a terminal object of the full subcategory $\{ g\in \s{C}_{/x}\tx{ ; } f\circ g=0 \} \subseteq \s{C}_{/x}$.  This approach likely extends to (possibly non-pointed) $\infty$-categories which are nerves of categories enriched in pointed Kan complexes; see (\cite{thesis}, \S 1.2).  However, to avoid having too many parallel and variant constructions, we will only address the case of pointed $\infty$-categories in this paper.

The content of this paper comes mainly from the first chapter of the author's thesis (\cite{thesis}, \S 1).

\section{Notation and terminology}

In this paper, we will use the notation, terminology, and size conventions (\cite{htt}, \S 1.2.15) of \cite{htt}, \cite{ha}, and \cite{kerodon}, with a couple of exceptions given in this section.  
\begin{itemize}
\item For ease of indexing, we will write $[-1]$ for the empty partially ordered set, and $\Delta^{-1}$ for the empty simplicial set.
\item Where there might be confusion, we call ordinary categories $1$-categories, to distinguish them from $\infty$-categories or catgories enriched over a monoidal category other than $\tx{Set}$.
\item The $1$-category of simplicial sets is denoted by $\tx{sSet}$, rather than by $\s{S}\tx{et}_\Delta$.  We will also use $\tx{sSet}$ to denote the simplicial category of simplicial sets.
\item Let $\tx{Kan}$ be the full simplicial subcategory of the simplicial category $\tx{sSet}$ spanned by the simplicial sets which are Kan complexes.  We denote by $\s{S}$ the homotopy coherent (a.k.a. simplicial) nerve $\tx{N}(\tx{Kan})$.  In particular, $\s{S}$ is an $\infty$-category.
\item For a topological category $\s{C}$, the \emph{homotopy category} $\tx{h}\s{C}$ is a $1$-category defined as in (\cite{htt}, 1.1.3.2), with $\tx{Hom}_{\tx{h}\s{C}} (x,y) = \pi_0 \tx{Map}_\s{C} (x,y)$ for $x,y\in \s{C}$.  

In particular, let $\s{{CW}}$ be the topological category whose objects are CW complexes, where $\tx{Map}_{\s{{CW}}} (x,y)$ is the set of continuous maps from $x$ to $y$, equipped with the compactly generated version of the compact-open topology.  We denote the homotopy category of $\s{{CW}}$ by $\s{H}$, and we call $\s{H}$ the \emph{homotopy category of spaces} (\cite{htt}, 1.1.3.3).
\end{itemize}
\begin{definition}[\tx{}\cite{htt}, 1.1.5.14 and 1.2.2.1] For a simplicial set $S$, we denote by $\tx{h}S$ the homotopy category $\tx{h}\f{C}[S]$, considered as a category enriched over the homotopy category $\s{H}$ of spaces (\cite{htt}, 1.1.3.4), where $\f{C}[-]$ is the rigidification functor of (\cite{htt}, \S 1.1.5).  For vertices $x,y\in S$, we have $\tx{Map}_{\tx{h}S} (x,y) = | \tx{Map}_{\f{C}[S]} (x,y) | \in \s{H}$.  

We will also denote $\tx{Map}_{\tx{h}S} (x,y)$ by $\tx{Map}_{S} (x,y)$.
\end{definition}
\begin{definition} Let $\s{C}$ be an $\infty$-category.  We call an object $x\in \s{C}$ a \emph{zero object} of $\s{C}$ if $x$ is both initial and terminal. If the $\infty$-category $\s{C}$ contains a zero object, we say $\s{C}$ is a \emph{pointed} $\infty$-category.
\end{definition}
\begin{definition}[\tx{}\cite{kerodon}, 01PG and \cite{htt}, \S 1.2.2] Let $Y$ be a simplicial set with vertices $x_0, x_1 , \dots , x_n \in Y$ ($n\geq 0$).  We denote by $\tx{Hom}^\tx{M}_Y (x_0 , x_1, \dots , x_n)$ the fiber product \[\tx{Fun}(\Delta^n , Y ) \times_{\tx{Fun}(\tx{sk}_0 \Delta^n, Y )} \{ (x_0 , x_1, \dots , x_n) \} \tx{\; in }\tx{sSet}.\]
Note that if $Y$ is an $\infty$-category, then $\tx{Hom}^\tx{M}_Y (x_0 , x_1, \dots , x_n)$ is a Kan complex (\cite{kerodon}, 01JC).  

If $\s{C}$ is an $\infty$-category with objects $x,y,z \in \s{C}$, then the composition map 
\[ -\circ - :  \tx{Hom}^\tx{M}_\s{C} (y,z) \times \tx{Hom}^\tx{M}_\s{C} (x,y) \to \tx{Hom}^\tx{M}_\s{C} (x,z) \tx{\; in hKan} \]
is defined as in (\cite{kerodon}, 01PQ).
\end{definition}
\begin{definition}[\tx{}\cite{kerodon}, 028V] Let $U : \s{E} \to \s{C}$ be a coCartesian fibration of simplicial sets, and let $V : \s{{QC}}_{\tx{Obj}} \to \tx{Cat}_\infty = \s{{QC}}$ be the forgetful map of (\cite{kerodon}, 0212 and 028V).\footnote{Although \cite{kerodon} uses $\s{{QC}}$ to denote the $\infty$-category $\tx{Cat}_\infty = \tx{N}(\tx{Cat}_\infty^\Delta )$ of \cite{htt}, the simplicial set $\s{{QC}}_{\tx{Obj}}$ of (\cite{kerodon}, 0212) is not the same as $\s{Z}^\tx{op}$ of (\cite{htt}, \S 3.3.2), though they are closely related.  The forgetful map $V : \s{{QC}}_{\tx{Obj}} \to \tx{Cat}_\infty = \s{{QC}}$ from \cite{kerodon} is not quite the same as $(-)^\tx{op}$ of the universal Cartesian fibration of \cite{htt}, although they give equivalent classifications of coCartesian fibrations.}  We say a commutative diagram of simplicial sets
\[\begin{tikzcd}
\s{E} \arrow[]{r}{\tilde{\s{F}}} \arrow[]{d}{U} & \s{{QC}}_{\tx{Obj}} \arrow[]{d}{V} \\
\s{C} \arrow[]{r}{\s{F}} & \tx{Cat}_\infty
\end{tikzcd}\]
\emph{witnesses $\s{F}$ as a covariant transport representation of }$U$ if the induced map
\[ \s{E} \to \s{C}\times_{\tx{Cat}_\infty}  \s{{QC}}_{\tx{Obj}} \]
is an equivalence of coCartesian fibrations over $\s{C}$ (\cite{kerodon}, 0281).  We say $\s{F} : \s{C} \to \tx{Cat}_\infty$ \emph{is a covariant transport representation of }$U$ if there is a diagram of simplicial sets that witnesses $\s{F}$ as a covariant transport representation of $U$.
\end{definition}



\section{Convexity and the $(-)_\diamond$ construction}

In this section, we construct a functor $(-)_\diamond : \tx{sSet} \to \tx{sSet}$ (Definition~\ref{diamonddefn}), and show that for any simplicial set $Y$, the simplicial set $Y_\diamond$ is a (pointed) $\infty$-category if and only if $Y$ is an $\infty$-category.

\begin{definition} Let $\tx{Poset}$ be the category of partially ordered sets.  We say a morphism $f: P\to R$ in $\tx{Poset}$ is \emph{convex} if for every relation $p\leq p'$ in $P$ and every $\ol{r}\in R$ such that $f(p) \leq \ol{r} \leq f(p')$ in $R$, there exists $r\in P$ such that $p\leq r\leq p'$ in $P$ and $f(r) = \ol{r}$.  Equivalently, $f : P \to R$ is convex if it has the right lifting property with respect to the (coface) inclusion $\delta_1 : [1]\to [2]$.  Thus convex morphisms in $\tx{Poset}$ form a weakly cosaturated class (\cite{rezk}, 17.3).
\end{definition}
Let $\tx{Spine}[n] \subseteq \Delta^n$ be the spine of $\Delta^n$, i.e.
\[ \tx{Spine}[n]_k = \left\{ \sigma : [k]\to [n] \tx{ in }\db \tx{ ; }\sigma(k) \leq \sigma(0)+1  \right\}. \]
For $0\leq a \leq b \leq n$, we will use $\tx{Spine}\{ a\leq \dots \leq b \}$ to denote $\Delta^{\{ a, a+1 , \dots , b \}} \times_{\Delta^n} \tx{Spine}[n]$.
\begin{proposition}\label{convspineinnderan} If $i:P \to [n]$ is a monomorphism in $\tx{Poset}$, then the following conditions are equivalent: 
\begin{enumerate}
\item\label{convspineinnderan1} $i : P \to [n]$ is convex.
\item\label{convspineinnderan2} The pullback $\tx{Spine}[n] \times_{\Delta^n} \tx{N}(P) \hookrightarrow \tx{N}(P)$ of $\tx{Spine}[n] \hookrightarrow \Delta^n$ is inner anodyne.
\item\label{convspineinnderan3} $\tx{Spine}[n] \times_{\Delta^n} \tx{N}(P) \hookrightarrow \tx{N}(P)$ is a categorical equivalence.
\end{enumerate}
\end{proposition}
\begin{proof} \emph{\ref{convspineinnderan1} $\Rightarrow$~\ref{convspineinnderan2}}: Suppose $i$ is convex.  If $0\leq n\leq 1$, then $\tx{Spine}[n] \times_{\Delta^n} \tx{N}(P) \hookrightarrow \tx{N}(P)$ is an isomorphism, so in particular it's inner anodyne.

Now suppose $n>1$. If $i$ is an isomorphism, then $ \tx{Spine}[n] \times_{\Delta^n} \tx{N}(P) \hookrightarrow \tx{N}(P)$ is inner anodyne by (\cite{kerodon}, 00JA), so we may assume $i$ is not an isomorphism.  Consider the diagram of simplicial sets
\[\begin{tikzcd}
 \tx{Spine}\{1\leq \dots \leq n \}\times_{\Delta^n} \tx{N}(P) \arrow[hook]{r} & \Delta^{\{1,\dots , n \}}\times_{\Delta^n} \tx{N}(P)  \\
 \tx{Spine}\{1\leq \dots \leq n-1 \}\times_{\Delta^n} \tx{N}(P) \arrow[hook]{r} \arrow[hook]{u} \arrow[hook]{d} & \Delta^{\{1,\dots , n-1 \}}\times_{\Delta^n} \tx{N}(P)  \arrow[hook]{u}  \arrow[hook]{d} \\
 \tx{Spine}\{0\leq \dots \leq n-1 \}\times_{\Delta^n} \tx{N}(P)  \arrow[hook]{r}& \Delta^{\{0,\dots , n-1 \}}\times_{\Delta^n} \tx{N}(P) . 
\end{tikzcd}\]
By induction on $n$, each horizontal morphism is inner anodyne.  Since $i$ is a convex monomorphism but not an isomorphism, there can be no pair $p\leq p'$ in $P$ such that $i(p)=0$ and $i(p')=n$, so both
\[\begin{tikzcd} \tx{Spine}\{1\leq \dots \leq n-1 \}\times_{\Delta^n} \tx{N}(P) \arrow[hook]{r}\arrow[hook]{d} & \tx{Spine}\{1\leq \dots \leq n \}\times_{\Delta^n} \tx{N}(P) \arrow[hook]{d} \\
 \tx{Spine}\{0\leq \dots \leq n-1 \}\times_{\Delta^n} \tx{N}(P)\arrow[hook]{r} & \tx{Spine}[n]  \times_{\Delta^n} \tx{N}(P)\end{tikzcd} \]
and
\[\begin{tikzcd} \Delta^{\{1,\dots , n-1 \}}\times_{\Delta^n} \tx{N}(P) \arrow[hook]{r}\arrow[hook]{d}& \Delta^{\{1,\dots , n \}}\times_{\Delta^n} \tx{N}(P)\arrow[hook]{d} \\
 \Delta^{\{0,\dots , n-1 \}}\times_{\Delta^n} \tx{N}(P)\arrow[hook]{r} & \tx{N}(P)\end{tikzcd} \]
are pushout squares.  Thus $\tx{Spine}[n] \times_{\Delta^n} \tx{N}(P) \hookrightarrow \tx{N}(P)$ is inner anodyne by the glueing lemma for inner anodyne morphisms (\cite{stevenson}, Lemma 2.5).

\hfill 

\emph{\ref{convspineinnderan2} $\Rightarrow$~\ref{convspineinnderan3}}: Every inner anodyne morphism is a categorical equivalence (\cite{kerodon}, 01EJ).

\hfill

\emph{\ref{convspineinnderan3} $\Rightarrow$~\ref{convspineinnderan1}}: Finally, suppose $i : P\to [n]$ is not convex, so there exist $0\leq a < b <c \leq n$ and a commutative square in $\tx{Poset}$
\[\begin{tikzcd}
{[1]} \arrow[]{r}{g} \arrow[]{d}[swap]{\delta_1} & P \arrow[]{d}{i} \\
{[2]} \arrow[dashed]{ur}{\not\exists} \arrow[hook]{r}[swap]{a<b<c} & {[n]}
\end{tikzcd}\]
with no lift.  Thus the connected component of $g(0)$ in $\tx{Spine}[n] \times_{\Delta^n} \tx{N}(P) $ does not contain $g(1)$, and there is an isomorphism of simplicial sets $\varphi : \tx{Spine}[n] \times_{\Delta^n} \tx{N}(P) \xrightarrow{\sim} A \coprod B$ with $\varphi (g(0))\in A$ and $\varphi (g(1)) \in B$.  The rigidification functor $\f{C}[ - ] : \tx{sSet} \to \cat_{\tx{sSet}}$ preserves coproducts, so $\tx{Hom}_{\f{C}[\tx{Spine}[n] \times_{\Delta^n} \tx{N}(P) ]} (g(0),g(1))$ is empty.

However, $g \in \tx{Hom}_{\f{C}[\tx{N}(P)]} (g(0), g(1))$, so 
\[ \f{C}[\tx{Spine}[n] \times_{\Delta^n} \tx{N}(P) ] \to  \f{C}[\tx{N}(P)] \]
cannot be a weak equivalence of simplicial categories under the Bergner model structure (\cite{bergner}, \S 4.3), so $\tx{Spine}[n] \times_{\Delta^n} \tx{N}(P) \hookrightarrow \tx{N}(P)$ cannot be a categorical equivalence.
\end{proof}

The following construction, $(-)_\diamond : \tx{sSet} \to \tx{sSet}$, is an analogue of the (co)cone constructions $\Delta^0 \star -$ and $-\star \Delta^0$.  When $\s{C}$ is an $\infty$-category, the constructions $\s{C} \mapsto \Delta^0 \star \s{C}$ and $\s{C}\mapsto \s{C} \star \Delta^0 $ freely adjoin initial and terminal objects.  As we will make precise below (in Theorem~\ref{diamondisfreezerothm}), the construction $\s{C} \mapsto \s{C}_\diamond$ freely adjoins a zero object to an $\infty$-category $\s{C}$.

\begin{definition}\label{diamonddefn} We say two spans $\sigma = \big( A\leftarrow C \to B \big)$, $\sigma' = \big( A \leftarrow C' \to B \big)$ of simplicial sets are \emph{isomorphic rel endpoints} if there is an isomorphism of spans from $\sigma$ to $\sigma'$ which restricts to the identity on $A$ and $B$:
\[\begin{tikzcd} A \arrow[equal]{d}{} & C \arrow[]{d}{\rotatebox[origin=c]{90}{$\sim$}} \arrow[]{l} \arrow[]{r} & B \arrow[equal]{d} \\
A & C' \arrow[]{r} \arrow[]{l} & B . \end{tikzcd}\]
In other words, $\sigma$ and $\sigma'$ are isomorphic rel endpoints if they are isomorphic as objects of $\tx{sSet}_{A/ /B}$.

For a given span of simplicial sets $\sigma = \big( A\leftarrow C \to B \big)$, we denote the class of all spans isomorphic to $\sigma$ rel endpoints by $[A\leftarrow C \to B ]$.

\hfill

For a simplicial set $Y$, define the simplicial set $Y_\diamond$ as follows:
\begin{itemize}
\item An $n$-simplex of $Y_\diamond$ is a class of spans of simplicial sets $[\tx{N}([n]) \leftarrow \tx{N}(P) \to Y ]$, where $\tx{N}(P) \to \tx{N}([n])$ is the nerve of a convex monomorphism of partially ordered sets $P\to [n]$ in $\tx{Poset}$.  
\item For $\alpha : [m] \to [n]$ in $\db$, define $\alpha^* : (Y_\diamond)_n \to (Y_\diamond)_m$ by
\[ [ \tx{N}([n]) \leftarrow N(P) \to Y ] \mapsto [ \tx{N}([m]) \leftarrow \tx{N}([m])\times_{\tx{N}([n])} \tx{N}(P) \to Y ]. \]
\end{itemize}
Note that if $\kappa$ is an infinite cardinal and $Y$ is a $\kappa$-small simplicial set, then $Y_\diamond$ is also $\kappa$-small.  This construction gives a functor $(-)_\diamond : \tx{sSet} \to \tx{sSet}$, where for a morphism $f : Y \to Z$ in $\tx{sSet}$, the morphism $f_\diamond : Y_\diamond \to Z_\diamond$ is given by 
\[\begin{tikzcd}
\tx{Hom}_\tx{sSet} (\Delta^n ,  Y_\diamond) \arrow[]{r}{f_\diamond \circ -} & \tx{Hom}_\tx{sSet} (\Delta^n ,  Z_\diamond) \\[-1.5em]
\Big[ \Delta^n \xhookleftarrow{\psi} \tx{N}(P) \xrightarrow{\sigma} Y  \Big] \arrow[mapsto]{r} & \Big[ \Delta^n \xhookleftarrow{\psi} \tx{N}(P) \xrightarrow{f\circ \sigma} Z \Big] .
\end{tikzcd}\]
We have a canonical inclusion $\iota : Y \hookrightarrow Y_\diamond$, natural in $Y$, given by 
\[\begin{tikzcd} Y_n \arrow[]{r} & (Y_\diamond)_n \\[-1.5em]
\left( \tx{N}([n]) \xrightarrow{\sigma} Y \right) \arrow[mapsto]{r} & {[ \tx{N}([n]) \xleftarrow{\tx{id}} \tx{N}([n]) \xrightarrow{\sigma} Y ]}. \end{tikzcd} \]
There is exactly one vertex of $Y_\diamond$ which is not in the image of $\iota$, namely $[\Delta^0 \leftarrow \tx{N}(\emptyset) \to Y]$.  Let $z_Y$ denote this vertex.
\end{definition}
\begin{construction}\label{triangletodiamondinclusion} We have a natural transformation $\zeta : (-)^\triangleright \Rightarrow (-)_\diamond$ of functors of $1$-categories from $(-)^\triangleright : \tx{sSet} \to \tx{sSet}$ to $(-)_\diamond : \tx{sSet} \to \tx{sSet}$, whose component at $Y \in \tx{sSet}$ is given by 
\[\begin{tikzcd}
\tx{Hom}_{\tx{sSet}} (\Delta^n , Y^\triangleright ) \arrow[hook]{r}{\zeta_Y \circ -} & \tx{Hom}_{\tx{sSet}} (\Delta^n , Y_\diamond ) \\[-1.5em]
\left( \Delta^{\{ 0 ,\dots , m \}} \xrightarrow{\sigma} Y \right) \arrow[mapsto]{r} & \left[ \Delta^n \hookleftarrow \Delta^{\{ 0 ,\dots , m \}} \xrightarrow{\sigma} Y \right]
\end{tikzcd}\]
where we have identified an $n$-simplex of $Y^\triangleright$, given by   $\Delta^{\{ 0 ,\dots , m \}} \star \Delta^\ell \xrightarrow{\sigma \star !} Y \star \Delta^0$ for some $-1\leq m,\ell \leq n$ such that $m + \ell +1=n$, with $\sigma : \Delta^{\{ 0 ,\dots , m \}} \to Y$.

\hfill

In the same way, we define a natural transformation $\atez : (-)^\triangleleft \Rightarrow (-)_\diamond$ by
\[\begin{tikzcd}
\tx{Hom}_{\tx{sSet}} (\Delta^n , Y^\triangleleft ) \arrow[hook]{r}{\atez_Y \circ -} & \tx{Hom}_{\tx{sSet}} (\Delta^n , Y_\diamond ) \\[-1.5em]
\left( \Delta^{\{ m ,\dots , n \}} \xrightarrow{\sigma} Y \right) \arrow[mapsto]{r} & \left[ \Delta^n \hookleftarrow \Delta^{\{ m ,\dots , n \}} \xrightarrow{\sigma} Y \right] 
\end{tikzcd}\]
for $Y \in \tx{sSet}$.  These natural transformations give commutative diagrams of natural transformations of ($1$-categorical) functors $\tx{sSet} \to \tx{sSet}$
\begin{equation}\label{triangletodiamondinclusioncommutativetriangles}\begin{tikzcd}
\tx{id}_{\tx{sSet}} \arrow[Rightarrow]{r} \arrow[Rightarrow]{dr}[swap]{\iota} & (-)^\triangleright \arrow[Rightarrow]{d}{\zeta} & & \tx{id}_{\tx{sSet}} \arrow[Rightarrow]{d} \arrow[Rightarrow]{dr}{\iota} & \\
& (-)_\diamond & & (-)^\triangleleft \arrow[Rightarrow]{r}[swap]{\atez} & (-)_\diamond .
\end{tikzcd}\end{equation}
In Construction~\ref{triangletodiamondinclusionsimplicial}, we'll extend these to simplicial natural transformations of simplicial functors.
\end{construction}
\begin{rmk}\label{homrdiamondrmk} For any simplicial set $Y$ and any $x,y\in Y \subseteq Y_\diamond$, we have isomorphisms of simplicial sets
$ \tx{Hom}^\tx{L}_Y (x,y) \coprod \Delta^0  \cong \tx{Hom}^\tx{L}_{Y_\diamond} (x,y)$ and $
\tx{Hom}^\tx{R}_Y (x,y) \coprod \Delta^0 \cong \tx{Hom}^\tx{R}_{Y_\diamond} (x,y)$, and similarly $\tx{Hom}^\tx{R}_{Y_\diamond} (z_Y ,z_Y) \cong \tx{Hom}^\tx{L}_{Y_\diamond} (x, z_Y) \cong \tx{Hom}^\tx{R}_{Y_\diamond} (z_Y ,y) \cong \Delta^0$.\end{rmk}
\begin{proposition}\label{convmakesinfty} For any simplicial set $Y$, $Y_\diamond$ is an $\infty$-category if and only if $Y$ is an $\infty$-category.  \end{proposition}
\begin{proof} For $0<k<n$, a map $\ell : \Lambda^n_k \to Y_\diamond$ is given by a family 
\begin{equation}\label{convhorncompat} \Big( [ \Delta^{\{ 0,\dots , \widehat{j} , \dots , n \}} \leftarrow \tx{N}(P_j) \to Y ] \Big)_{0\leq j \leq n, \tx{ } j\neq k} ,\end{equation}
with the condition that for all $0\leq j \leq j' \leq n$ with $j\neq k \neq j'$, there is an isomorphism rel endpoints 
 \[\begin{tikzcd}
{ \Delta^{\{0, \dots , \widehat{j} , \dots , \widehat{j'} , \dots , n \}} } \arrow[equal]{d} & { \Delta^{\{0, \dots , \widehat{j} , \dots , \widehat{j'} , \dots , n \}} \times_{\Delta^{\{0, \dots ,  \widehat{j} , \dots , n \}}} \tx{N}(P_j) } \arrow[]{r}\arrow[]{l}\arrow[]{d}{\rotatebox[origin=c]{90}{$\sim$}} & Y  \arrow[equal]{d}\\
{ \Delta^{\{0, \dots , \widehat{j} , \dots , \widehat{j'} , \dots , n \}}} & {\Delta^{\{0, \dots , \widehat{j} , \dots , \widehat{j'} , \dots , n \}} \times_{\Delta^{\{0, \dots ,  \widehat{j'} , \dots , n \}}} \tx{N}(P_{j'}) } \arrow[]{r}\arrow[]{l}& Y 
 \end{tikzcd}\]

\hfill

To prove $Y$ is an $\infty$-category if $Y_\diamond$ is an $\infty$-category, it suffices to show that for any $0<k<n$ and any commutative diagram of simplicial sets
\[ \begin{tikzcd} \Lambda^n_k \arrow[hook]{d} \arrow[bend left]{rr}{\ell} \arrow[]{r} & Y\arrow[hook]{r}{\iota} & Y_\diamond  \\
\Delta^n \arrow[dashed]{ur} \arrow[bend right=20]{urr}[swap]{\sigma}& & \end{tikzcd} \]
in which $\iota : Y \to Y_\diamond$ is the canonical inclusion, the map $\sigma$ factors through $\iota$.  If $\sigma$ is represented by a span $\Delta^n \xleftarrow{\nu} \tx{N}(P) \to Y$, then $\sigma$ factoring through $\iota$ is equivalent to $\nu$ being an isomorphism.

Since $\ell$ factors through $\iota$, each map $\tx{N}(P_j) \to \Delta^{\{ 0,\dots , \widehat{j} , \dots , n \}}$ is an isomorphism, hence each map $ \Delta^{\{ 0,\dots , \widehat{j} , \dots , n \}} \times_{\Delta^n} \tx{N}(P) \to \Delta^{\{ 0,\dots , \widehat{j} , \dots , n \}}$ is an isomorphism.  In particular, this holds for $j \in \{ 0,n\}$.  But $\tx{Spine}[n] \subseteq \big( \Delta^{\{ 0,\dots , n-1\}} \cup \Delta^{\{ 1, \dots , n \}} \big)  \subseteq \Delta^n$, so $\tx{Spine}[n] \hookrightarrow \Delta^n$ factors through $\nu$.  Since $\nu$ is also the nerve of a monomorphism in $\tx{Poset}$, $\nu$ must be an isomorphism.

\hfill

Conversely, suppose $Y$ is an $\infty$-category, and suppose a map $\ell : \Lambda^n_k \to Y_\diamond$, with $0<k<n$, given by a family of classes as in (\ref{convhorncompat}).  Let $\s{P}$ be the set of nonempty subsets of $\{ 0, 1, \dots , \widehat{k} , \dots  ,n \}$, considered as a category whose morphisms are subset inclusions.  For each object $S = \{ s_1 , \dots , s_m \}$ in $\s{P}$, pick a partially ordered set $P_S \cong P_{s_1} \times_{[n]} \dots \times_{[n]} P_{s_m}$, where the map $P_j \to [n]$ is the composition $P_j \hookrightarrow [n-1] \xrightarrow{\delta_j} [n]$ for each $j\in \{ 0,  \dots , \widehat{k} , \dots  ,n \}$.  We have a functor 
\[ \begin{tikzcd} \s{P}^{\tx{op}} \arrow[]{r}{F} & \tx{Poset} \\[-1.5em]
S \arrow[mapsto]{r} & P_S  \end{tikzcd} \]
Let $\bigcup_{j\neq k} P_j := \colim_{\s{P}^\tx{op}} F$ in $\tx{Poset}$, and $\bigcup_{j\neq k} \tx{N}(P_j) := \colim_{\s{P}^\tx{op}} \tx{N}\circ F$ in $\tx{sSet}$.  We have monomorphisms $\psi : \bigcup_{j\neq k} P_j \hookrightarrow [n]$ in $\tx{Poset}$ and $\bigcup_{j\neq k} \tx{N}(P_j) \hookrightarrow \tx{N} \big( \bigcup_{j\neq k} P_j \big) \hookrightarrow \Delta^n$ in $\tx{sSet}$.  The family (\ref{convhorncompat}) gives a map $\bigcup_{j\neq k} \tx{N}(P_j) \to Y$.  

To show that $\ell : \Lambda^n_k \to Y_\diamond$ extends along $\Lambda^n_k \hookrightarrow \Delta^n$ to an $n$-simplex $\sigma : \Delta^n \to Y_\diamond$, it suffices to show that $\bigcup_{j\neq k} P_j \hookrightarrow [n]$ is convex and $\bigcup_{j\neq k} \tx{N}(P_j)  \hookrightarrow \tx{N}\left( \bigcup_{j\neq k} P_j \right)$ is inner anodyne.  Indeed, since $Y$ is an $\infty$-category, the second condition would guarantee a dashed extension
\[\begin{tikzcd} & \bigcup_{j\neq k} \tx{N}(P_j) \arrow[hook]{d} \arrow[hook, bend right]{dl} \arrow[bend left]{dr} & \\ 
\Delta^n & \tx{N}\left( \bigcup_{j\neq k} P_j \right) \arrow[hook]{l} \arrow[dashed]{r} & Y,
\end{tikzcd}\]  
and, by the first condition, the class $[\Delta^n \leftarrow \tx{N}\big( \bigcup_{j\neq k} P_j \big) \to Y]$ would give an extension $\sigma$ of $\ell$.
 
Consider the pullbacks of simplicial sets
\[\begin{tikzcd} \tx{Spine}[n] \times_{\Delta^n} \left( \bigcup_{j\neq k} \tx{N}(P_j) \right) \arrow[hook]{r}{\alpha} \arrow[hook]{d}[swap]{\beta} & \bigcup_{j\neq k} \tx{N}(P_j) \arrow[hook]{d}\\
\tx{Spine}[n] \times_{\Delta^n} \left( \tx{N}\left( \bigcup_{j\neq k} P_j \right) \right) \arrow[hook]{d}\arrow[hook]{r}{\gamma}& \tx{N}\left( \bigcup_{j\neq k} P_j \right) \arrow[hook]{d}{\tx{N}(\psi)} \\
\tx{Spine}[n] \arrow[hook]{r} & \Delta^n .
 \end{tikzcd}\]
 
 If a map of simplicial sets $f : A \to B$ is a monomorphism, bijective on vertices, and a categorical equivalence, and if moreover $B$ is an $\infty$-category, then $f$ is inner anodyne (\cite{stevensoncounter}).  Thus to show $\bigcup_{j\neq k} \tx{N}(P_j)  \hookrightarrow \tx{N}\left( \bigcup_{j\neq k} P_j \right)$ is inner anodyne, it suffices to show it's a categorical equivalence.  By Proposition~\ref{convspineinnderan}, if $\psi : \bigcup_{j\neq k} P_j \hookrightarrow [n]$ is convex, then $\gamma$ is a categorical equivalence.  Thus it remains to show that $\alpha$ and $\beta$ are categorical equivalences, and that $\psi$ is convex.

\hfill

\emph{$\alpha$ is a categorical equivalence}: Consider the two functors
\[ \begin{tikzcd} \s{P}^{\tx{op}} \arrow[]{r}{\tx{N}\circ F} & \tx{sSet}  & & \s{P}^{\tx{op}} \arrow[]{r}{G} & \tx{sSet} \\[-1.5em]
S \arrow[mapsto]{r} & \tx{N}(P_S)   & & S \arrow[mapsto]{r} & \tx{Spine}[n] \times_{\Delta^n} \tx{N}(P_S) . \end{tikzcd} \]
We have a natural transformation $\eta : G \to \tx{N}\circ F$, and by Proposition~\ref{convspineinnderan}, each $\eta_S$ is a categorical equivalence.  The map $\alpha$ is the induced map $\colim G \to \colim (\tx{N}\circ F)$, so to show that $\alpha$ is a categorical equivalence, it suffices to show that $G$ and $\tx{N}\circ F$ are cofibrant in the projective model structure on $\tx{sSet}^{\s{P}^\tx{op}}$, where $\tx{sSet}$ is given the Joyal model structure.  

We can consider $\s{P}^\tx{op}$ as a Reedy category, with $(\s{P}^\tx{op} )_+ = \s{P}^\tx{op}$.  For any $S \in \s{P}^\tx{op}$, let $\s{P}^{\tx{op}}_S \subseteq \s{P}^{\tx{op}}$ be the full subcategory spanned by the $S' \in \s{P}$ that properly contain $S$.  For $S \in \s{P}^\tx{op}$, the comparison maps $L_S (\tx{N}\circ F) \to \tx{N}\circ F(S)$ and $L_S G \to G(S)$ from the latching objects are $\colim_{\s{P}^\tx{op}_S} (\tx{N}\circ F)|_{\s{P}^\tx{op}_S} \to \tx{N}(P_S)$ and $\colim_{\s{P}^\tx{op}_S} G |_{\s{P}^\tx{op}_S} \to \tx{Spine}[n] \times_{\Delta^n} \tx{N}(P_S)$ respectively.  Both are monomorphisms for each $S\in \s{P}^\tx{op}$, so $\tx{N}\circ F$ and $G$ are Reedy cofibrant.  Thus $\tx{N}\circ F$ and $G$ are cofibrant in the projective model structure (\cite{htt}, A.2.9.22). 

\hfill

\emph{$\beta$ is a categorical equivalence}:  Suppose $p', p''$ in $ \bigcup_{j\neq k} P_j$ such that $p'\leq p''$ in $ \bigcup_{j\neq k} P_j$ and $\psi (p') +1=\psi (p'')$.  Then there is some family $j_1 , \dots , j_m \in \{ 0, 1,\dots , \widehat{k} , \dots ,n \}$ and elements $p_{0}, \dots , p_m \in \bigcup_{j\neq k} P_j$ such that $p' =p_{0}$, $p'' = p_{m}$ and $p_{i-1} \leq p_{i}$ in $P_{j_i}$ for $1\leq i \leq m$.  Then $\psi (p_{i-1}) \leq \psi (p_i)$ for each $1\leq i\leq m$, and there is some $i$ such that $p_{i-1} = p'$ and $p_{i} = p''$, with $p' \leq p''$ in $P_{j_i}$.  Thus $\tx{Spine}[n] \times_{\Delta^n} \left( \bigcup_{j\neq k} \tx{N}(P_j) \right)$ contains every edge of $\tx{Spine}[n] \times_{\Delta^n} \left( \tx{N}\left( \bigcup_{j\neq k} P_j \right) \right)$, and both are $1$-skeletal, so $\beta$ is an isomorphism.

\hfill

\emph{$\psi$ is convex}:  First suppose $n=2$.  If $\psi$ is an isomorphism, then in particular it's convex.  If $\psi$ is not an isomorphism, then $\gamma$ is an isomorphism, so $\psi$ is convex by Proposition~\ref{convspineinnderan}.

Next suppose $n\geq 4$, and suppose $\psi$ is not convex, so there exist $0\leq a < b < c \leq n$ and a commutative square in $\tx{Poset}$
\begin{equation}\label{convinfcontradictionsq}\begin{tikzcd}
{[1]} \arrow[]{r}{g} \arrow[hook]{d}[swap]{\delta_1} & \bigcup_{j\neq k} P \arrow[hook]{d}{\psi} \\
{[2]} \arrow[dashed]{ur}{\not\exists} \arrow[hook]{r}[swap]{a<b<c} & {[n]}
\end{tikzcd}\end{equation}
with no lift.  Since $n\geq 4$, $\Lambda^n_k$ contains the ($2$- or $3$-)simplex $\Delta^{\{a,b,c,k\}}$.  Say $\Delta^{\{a,b,c,k\}}$ is a simplex of $ \Delta^{\{ 0 , \dots, \widehat{j_0} , \dots , n \}} \subseteq \Delta^n$.   But 
\[\begin{tikzcd}
P_{j_0}\arrow[hook]{d} \arrow[hook]{r} & \bigcup_{j\neq k} P \arrow[hook]{d}{\psi} \\
{[n-1] } \arrow[]{r}{\delta_{j_0}}& {[n]}
\end{tikzcd}\]
is a pullback square in $\tx{Poset}$, so the non-existence of the lift in (\ref{convinfcontradictionsq}) contradicts the assumption that $P_{j_0} \to [n-1]$ is convex.

Finally, suppose $n=3$, and suppose $\psi$ is not convex, so there exist $0\leq a<b<c \leq 3$ and a square (\ref{convinfcontradictionsq}) in $\tx{Poset}$ with no lift.  If $\Delta^{\{a,b,c \}}$ is contained in $\Lambda^3_k$, then we have a contradiction by the same reasoning as in the case of $n\geq 4$.  Thus we may assume $\Delta^{\{a,b,c \}}$ is not contained in $\Lambda^3_k$, i.e. $k\not\in \{a,b,c\}$.  Thus $a=0$ and $c=3$, so $P_j \to [2]$, being convex, must be an isomorphism for $j\in \{1,2 \}\setminus \{k\}$.  If $k=1$, this means $P_0 \to [2]$ must be an isomorphism; if $k=2$, this means $P_3 \to [2]$ must be an isomorphism.  In either case, $\tx{Spine}[n] \hookrightarrow \Delta^n$ factors through $N(\psi)$, so $\psi$ must be an isomorphism, hence convex.
\end{proof}
\begin{corollary}\label{diamondhaszerocor} If $\s{C}$ is an $\infty$-category, then $\s{C}_\diamond$ is a pointed $\infty$-category with a unique zero object $z_\s{C} \in \s{C}_\diamond$ given by $[ \Delta^0 \leftarrow \tx{N}(\emptyset) \to \s{C} ]$.
\end{corollary}
\section{The $(-)_\diamond$ construction freely adjoins a zero object}
In this section we prove Theorem~\ref{diamondisfreezerothm}, which says that the simplicial functor $\s{C} \mapsto \s{C}_\diamond$ induces a left adjoint to the inclusion $\tx{Cat}_\infty^\tx{pt} \hookrightarrow \s{C}_\infty$.
\begin{definition}\label{catinftyptdefn} We define two simplicial subcategories of $\tx{Cat}_\infty^\Delta$, $\tx{Cat}_{\infty}^{\Delta, \tx{term}}$ and $\tx{Cat}_{\infty}^{\Delta, \tx{pt}}$, fibrant with respect to the Bergner model structure on simplicial categories (\cite{bergner}, \S 4.3), as follows:
\begin{itemize}
\item The objects of $\tx{Cat}_{\infty}^{\Delta, \tx{term}}$ (respectively $\tx{Cat}_{\infty}^{\Delta, \tx{pt}}$) are the small $\infty$-categories which contain a terminal (respectively zero) object.
\item For objects $\s{C}, \s{D}$ in $\tx{Cat}_{\infty}^{\Delta, \tx{term}}$, let $\tx{Fun}' (\s{C}, \s{D}) \subseteq \tx{Fun}(\s{C}, \s{D})$ be the subcategory spanned by the functors $\s{C} \to \s{D}$ which preserve terminal objects, and let $\tx{Map}_{\tx{Cat}_{\infty}^{\Delta, \tx{term}}} (\s{C}, \s{D})$ be the core $\tx{Fun}' (\s{C}, \s{D})^\simeq$ of $\tx{Fun}' (\s{C}, \s{D})$.  
  
  Likewise for objects $\s{C}, \s{D}$ in $\tx{Cat}_{\infty}^{\Delta, \tx{pt}}$, let $\tx{Fun}'' (\s{C}, \s{D}) \subseteq \tx{Fun}(\s{C}, \s{D})$ be the subcategory spanned by the functors $\s{C} \to \s{D}$ which preserve zero objects, and let $\tx{Map}_{\tx{Cat}_{\infty}^{\Delta, \tx{pt}}} (\s{C}, \s{D})$ be the core $\tx{Fun}'' (\s{C}, \s{D})^\simeq$ of $\tx{Fun}'' (\s{C}, \s{D})$.
\end{itemize}
We denote by $\maybehat{\tx{Cat}}_\infty^\tx{term}$ the homotopy coherent nerve $\tx{N}(\tx{Cat}_{\infty}^{\Delta, \tx{term}}) \subseteq \tx{N}(\tx{Cat}_{\infty}^{\Delta}) = \tx{Cat}_\infty$ and by  $\maybehat{\tx{Cat}}_\infty^\tx{pt}$ the homotopy coherent nerve $\tx{N}(\tx{Cat}_{\infty}^{\Delta, \tx{pt}}) \subseteq  \tx{Cat}_\infty$, so we have inclusions $\maybehat{\tx{Cat}}_\infty^\tx{pt}\subseteq \maybehat{\tx{Cat}}_\infty^\tx{term} \subseteq \maybehat{\tx{Cat}}_\infty$.  

Note that $\maybehat{\tx{Cat}}_\infty^\tx{pt}\hookrightarrow \maybehat{\tx{Cat}}_\infty^\tx{term}$ is fully faithful, but $\maybehat{\tx{Cat}}_\infty^\tx{term} \hookrightarrow \maybehat{\tx{Cat}}_\infty$ is not.
\end{definition}
\begin{construction}\label{prodtrianglerightnattransf} We define a (1-categorical) natural transformation $\chi$
\[\begin{tikzcd} \tx{sSet} \times \tx{sSet} \arrow[bend left=20]{rr}{(X,Y) \mapsto X \times (Y^\triangleright)}  \arrow[bend right=20]{rr}[swap]{(X,Y) \mapsto (X \times Y)^\triangleright} & \Downarrow \chi & \tx{sSet} \end{tikzcd}\]
as follows: For simplicial sets $X$ and $Y$,  let $\chi_{X,Y} : X \times (Y^\triangleright) \to (X \times Y)^\triangleright$ be the morphism that takes an $m$-simplex 
\begin{align*} (\alpha, \beta) \in \tx{Hom}(\Delta^m , X) \times \tx{Hom}(\Delta^k , Y) &\subseteq \coprod_{-1\leq\ell\leq m} \tx{Hom}(\Delta^m , X) \times \tx{Hom}(\Delta^\ell , Y) \\
&  \cong \tx{Hom}\big( \Delta^m , X \times (Y^\triangleright) \big)
\end{align*}
to
\begin{align*} (\alpha |_{\Delta^{\{ 0,\dots , k \}}}, \beta)  \in \tx{Hom}(\Delta^k , X) \times \tx{Hom}(\Delta^k , Y) &\subseteq \coprod_{-1\leq\ell\leq m} \tx{Hom}(\Delta^\ell , X) \times \tx{Hom}(\Delta^\ell , Y) \\
&  \cong \tx{Hom}\big( \Delta^m , (X \times Y)^\triangleright \big),
\end{align*}
where if $k=-1$, then $\alpha |_{\Delta^{\{ 0,\dots , k \}}}$ means the unique element of $\tx{Hom}(\Delta^{-1}, X)$.  
\end{construction}
\begin{construction}\label{diamondprodhelpernattransf} We define a (1-categorical) natural transformation  $\varpi$ 
\[\begin{tikzcd} \tx{sSet} \times \tx{sSet} \arrow[bend left=20]{rr}{(X,Y) \mapsto X \times (Y_\diamond )}  \arrow[bend right=20]{rr}[swap]{(X,Y) \mapsto (X \times Y)_\diamond} & \Downarrow \varpi & \tx{sSet} \end{tikzcd}\]
as follows:  For simplicial sets $X$ and $Y$, let $\varpi_{X,Y} : X \times (Y_\diamond ) \to (X\times Y)_\diamond$ be the morphism that takes an $m$-simplex
\[ \left(  \Delta^m \xrightarrow{\alpha} X , \Big[ \Delta^m \xhookleftarrow{\psi} \tx{N}(P) \xrightarrow{\beta} Y \Big]  \right) \in \tx{Hom}(\Delta^m, X ) \times \tx{Hom}(\Delta^m , Y_\diamond ) \]
to the $m$-simplex 
\[ \Big[  \Delta^m \xhookleftarrow{\psi} \tx{N}(P) \xrightarrow{(\alpha\circ\psi, \beta)} X\times Y \Big] \in \tx{Hom}(\Delta^m , (X\times Y)_\diamond ). \]
\end{construction}
\begin{construction}\label{triangletodiamondinclusionsimplicial} We extend the leftmost commutative triangle of (\ref{triangletodiamondinclusioncommutativetriangles}) in Construction~\ref{triangletodiamondinclusion} to a commutative triangle of simplicial natural transformations of simplicial functors $\tx{sSet} \to \tx{sSet}$ as follows, using the notation of Constructions~\ref{prodtrianglerightnattransf} and~\ref{diamondprodhelpernattransf}:

For $\s{C}, \s{D} \in \tx{sSet}$ and $[n] \in \db$, the natural transformations $\chi$ and $\varpi$ respectively give maps
\[ \begin{tikzcd} 
\tx{Hom} (\Delta^n \times \s{C}, \s{D} ) \arrow[]{rr}{f \mapsto f^\triangleright \circ \chi_{\Delta^n ,\s{C}}} & & \tx{Hom} ( \Delta^n \times (\s{C}^\triangleright) , \s{D}^\triangleright  ) \\[-0.5em]
\tx{Hom} (\Delta^n \times \s{C}, \s{D} ) \arrow[]{rr}{f \mapsto f_\diamond \circ \varpi_{\Delta^n ,\s{C}}} & & \tx{Hom} ( \Delta^n \times (\s{C}_\diamond) , \s{D}_\diamond  ),
\end{tikzcd}\]
each natural in the triple $\left( [n] , \s{C}, \s{D} \right) \in \dop \times \tx{sSet}^\tx{op} \times \tx{sSet}$, hence maps of simplicial sets \[\theta_{\s{C},\s{D}}  : \tx{Fun}(\s{C} , \s{D}) \to \tx{Fun}(\s{C}^\triangleright , \s{D}^\triangleright), \tx{ } \vartheta_{\s{C},\s{D}}  : \tx{Fun}(\s{C} , \s{D}) \to \tx{Fun}(\s{C}_\diamond , \s{D}_\diamond ) \] natural in the pair $( \s{C}, \s{D})\in \tx{sSet}^\tx{op}\times \tx{sSet}$.  These respectively give extensions of $(-)^\triangleright$ and $(-)_\diamond$ to simplicial functors $\tx{sSet}\to\tx{sSet}$, and the simplicial natural transformations $\tx{id}_{\tx{sSet}} \Rightarrow (-)^\triangleright \xRightarrow{\zeta} (-)_\diamond$ are evident.
\end{construction}
\begin{lemma}\label{coningsimplfunctor} The functor $(-)^\triangleright = (-)\star \Delta^0 : \tx{sSet}\to \tx{sSet}$ induces a functor of $\infty$-categories $\maybehat{\tx{Cat}}_\infty \to\maybehat{\tx{Cat}}_\infty^\tx{term}$.
\end{lemma}
\begin{proof} The simplicial functor $(-)^\triangleright : \tx{sSet}\to \tx{sSet}$ of Construction~\ref{triangletodiamondinclusionsimplicial} restricts to a simplicial functor $\tx{Cat}_\infty^\Delta \to \tx{Cat}_\infty^\Delta$.  Taking the homotopy coherent nerve of the latter gives a functor of $\infty$-categories $(-)^\triangleright : \tx{Cat}_\infty \to \tx{Cat}_\infty$ which factors through the inclusion $\tx{Cat}_\infty^\tx{term} \hookrightarrow \tx{Cat}_\infty$.
\end{proof}
\begin{lemma}\label{coneptsimplnatr} The natural transformation $\iota : \tx{id}_\tx{sSet} \Rightarrow (-)^\triangleright$ of functors $\tx{sSet}\rightrightarrows \tx{sSet}$ given by the natural inclusion $X \hookrightarrow X^\triangleright$ induces a natural transformation $\Delta^1 \to \tx{Fun}(\tx{Cat}_\infty , \tx{Cat}_\infty)$ from $\tx{id}_{\tx{Cat}_\infty}$ to the functor $(-)^\triangleright : \tx{Cat}_\infty \to \tx{Cat}_\infty$ of Lemma~\ref{coningsimplfunctor}.
\end{lemma}
\begin{proof} 
Restricting the simplicial natural transformation $\tx{id}_{\tx{sSet}} \Rightarrow (-)^\triangleright$ to the simplicial subcategory $\tx{Cat}_\infty^\Delta$ gives a simplicial natural transformation from $\tx{id}_{\tx{Cat}_\infty^\Delta}$ to the simplicial functor $(-)^\triangleright : \tx{Cat}_\infty^\Delta \to \tx{Cat}_\infty^\Delta$.  Let $[1]_{\Delta^0}$ be the simplicial category with two objects $X$ and $Y$ and
\[ \tx{Map}_{[1]_{\Delta^0}} (Z,Z') = \begin{cases}
\Delta^{-1} & \tx{\; if }Z=Y, Z'=X \\
\Delta^0 & \tx{\; otherwise}.
\end{cases} \]
We can identify the above simplicial natural transformation of parallel simplicial functors $ \tx{Cat}_\infty^\Delta \rightrightarrows \tx{Cat}_\infty^\Delta$ with a simplicial functor $[1]_{\Delta^0} \times \tx{Cat}_\infty^\Delta \to \tx{Cat}_\infty^\Delta$.  Taking coherent nerves gives a morphism of simplicial sets $\Delta^1 \times \tx{Cat}_\infty \cong \tx{N}([1]_{\Delta^0}) \times \tx{N}(\tx{Cat}_\infty^\Delta) \to \tx{N}(\tx{Cat}_\infty^\Delta) = \tx{Cat}_\infty$, which we can curry to an edge $\Delta^1 \to \tx{Fun}(\tx{Cat}_\infty , \tx{Cat}_\infty)$ from $\tx{id}_{\tx{Cat}_\infty}$ to $(-)^\triangleright$.
\end{proof}
By the same token, we have corresponding lemmas for the natural inclusion $\tx{id} \Rightarrow (-)_\diamond$:
\begin{lemma}\label{diamondinftyfunctor} The simplicial functor $(-)_\diamond$ of Construction~\ref{triangletodiamondinclusionsimplicial} induces a functor of $\infty$-categories $\tx{Cat}_\infty \to \tx{Cat}_\infty^\tx{pt}$.
\end{lemma}
\begin{lemma}\label{diamondinftynaturaltransf} The simplicial natural transformation $\tx{id} \Rightarrow (-)_\diamond$ of Construction~\ref{triangletodiamondinclusionsimplicial} induces a natural transformation $\Delta^1 \to \tx{Fun}(\tx{Cat}_\infty , \tx{Cat}_\infty)$ from $\tx{id}_{\tx{Cat}_\infty}$ to $(-)_\diamond$.  Likewise, the simplicial natural transformation $\zeta : (-)^\triangleright \Rightarrow (-)_\diamond$ induces a natural transformation $\Delta^1 \to \tx{Fun}(\tx{Cat}_\infty , \tx{Cat}_\infty^\tx{term})$ from $(-)^\triangleright$ to $(-)_\diamond$.
\end{lemma}
\begin{lemma}\label{frakztriangle} We have a morphism of simplicial sets $\f{Z}: \Delta^2 \to \tx{Fun}(\tx{Cat}_\infty , \tx{Cat}_\infty )$ where $\f{Z} |_{\Delta^{\{0,1 \}}}$ is the natural transformation from $\tx{id}_{\tx{Cat}_\infty}$ to  $(-)^\triangleright$ of Lemma~\ref{coneptsimplnatr},  $\f{Z} |_{\Delta^{\{0,2 \}}}$ is the natural transformation $\tx{id}_{\tx{Cat}_\infty}$ to $(-)_\diamond$ of Lemma~\ref{diamondinftynaturaltransf}, and $\f{Z} |_{\Delta^{\{1,2 \}}}$ is the natural transformation $\zeta$ of Lemma~\ref{diamondinftynaturaltransf}.
\end{lemma}
\begin{proposition}\label{catptadj} The functor $-\star \Delta^0 = (-)^\triangleright : \maybehat{\tx{Cat}}_\infty \to \maybehat{\tx{Cat}}_\infty^\tx{term}$ of Lemma~\ref{coningsimplfunctor} is left adjoint to the forgetful functor $U:\maybehat{\tx{Cat}}_\infty^\tx{term} \hookrightarrow \maybehat{\tx{Cat}}_\infty$.
\end{proposition}
\begin{proof} By (\cite{htt}, 5.2.2.8), it suffices to show that the natural transformation of Lemma~\ref{coneptsimplnatr}, $\tx{id}_{ \maybehat{\tx{Cat}}_\infty} \to U\circ (-\star \Delta^0)$, is a unit.  By (\cite{kerodon}, 01EZ), it suffices to show that for any $\infty$-category $\s{C}$ and any $\s{D} \in \maybehat{\tx{Cat}}_\infty^\tx{term}$, restriction along the inclusion $\s{C} \hookrightarrow \s{C}\star \Delta^0$ induces a trivial Kan fibration
\[ \tx{Fun}' (\s{C}\star \Delta^0 , \s{D}) \to \tx{Fun}(\s{C},\s{D}), \]
where $\tx{Fun}' (\s{C}\star \Delta^0 , \s{D}) \subseteq \tx{Fun} (\s{C}\star \Delta^0 , \s{D})$ is the full subcategory spanned by the functors $\s{C}\star \Delta^0 \to \s{D}$ that take the vertex of $\Delta^0$ to a terminal object of $\s{D}$. 

Let $A$ be the collection of non-degenerate simplices of $\s{C}\star \Delta^0$ which are not contained in $\s{C}$, and for $(\sigma : \Delta^n \to \s{C}\star\Delta^0 ) \in A$, let $d(\sigma)$ be $n$, the dimension of the simplex.  Pick a well-ordering on $A$ such that $\sigma <\tau$ if $d(\sigma) < d(\tau)$, and let $\alpha$ be the order type of $A$, so we have an order-preserving bijection $\beta \mapsto \sigma_\beta$, where $\beta$ ranges over the ordinals $\{\beta ; \beta <\alpha \}$.  For each $\beta\leq \alpha$, let $\s{C}_\beta \subseteq \s{C}\star \Delta^0$ be the largest simplicial subset whose nondegenerate simplicies are either contained in $\s{C}$ or are of the form $\sigma_{\gamma}$ for some $\gamma < \beta$, so $\s{C}_0 = \s{C}$, $\s{C}_1 = \s{C} \coprod \Delta^0$, and $\s{C}_\alpha = \s{C}\star \Delta^0$.  Let $v$ be the vertex of $\Delta^0$ in $\s{C}\star \Delta^0$, and for $1\leq \beta\leq \alpha$, let $\tx{Fun}' (\s{C}_\beta , \s{D}) \subseteq \tx{Fun}(\s{C}_\beta , \s{D})$ be the full subcategory spanned by the functors that take $v$ to a terminal object of $\s{D}$.

The subcategory $\tx{Fun}' (\s{C}_1 , \s{D}) \subseteq \tx{Fun} (\s{C}_1 , \s{D})$ is replete, so this inclusion is an isofibration by (\cite{kerodon}, 01EX), and the restriction $\tx{Fun} (\s{C}_1 , \s{D}) \to \tx{Fun}(\s{C},\s{D})$ is an isofibration by (\cite{kerodon}, 01F3), so the restriction $\tx{Fun}' (\s{C}_1 , \s{D}) \to \tx{Fun} (\s{C} , \s{D}) $ is an isofibration.  For any $f \in \tx{Fun} (\s{C} , \s{D})$, the fiber $\tx{Fun}' (\s{C}_1 , \s{D}) \times_{\tx{Fun}(\s{C},\s{D})} \{ f \}$ is isomorphic to the full subcategory of $\s{D}$ spanned by the terminal objects, so it's a contractible Kan complex.  Thus $\tx{Fun}' (\s{C}_1 , \s{D}) \to \tx{Fun} (\s{C} , \s{D}) $ is a trivial Kan fibration.

It remains to show that for $1\leq \beta\leq \beta' \leq\alpha$, the restriction $\tx{Fun}' (\s{C}_{\beta'} , \s{D}) \to \tx{Fun}' (\s{C}_\beta , \s{D})$ is a trivial Kan fibration, which we'll prove by induction on $\beta'$.  If $\beta' = \beta$, this is clear.  If $\beta' >\beta$ is a limit ordinal then by the induction hypothesis, $\tx{Fun}' (\s{C}_{\beta'} , \s{D}) \to \tx{Fun}' (\s{C}_\beta , \s{D})$ is a transfinite limit of a tower of trivial Kan fibrations, hence a trivial Kan fibration.  It remains to show that if $\tx{Fun}' (\s{C}_{\beta'} , \s{D}) \to \tx{Fun}' (\s{C}_\beta , \s{D})$ is a trivial Kan fibration, then $\tx{Fun}' (\s{C}_{\beta'+1} , \s{D}) \to \tx{Fun}' (\s{C}_\beta , \s{D})$ is a trivial Kan fibration, for which it suffices to show that $\tx{Fun}' (\s{C}_{\beta'+1} , \s{D}) \to \tx{Fun}' (\s{C}_{\beta'} , \s{D})$ is a trivial Kan fibration.  This map is a pullback of the map $\tx{Fun} (\s{C}_{\beta'+1} , \s{D}) \to \tx{Fun} (\s{C}_{\beta'} , \s{D})$, hence an isofibration, so it remains to show that its fibers are contractible.

We have a pushout
\[\begin{tikzcd}
\partial \Delta^n \arrow[]{r} \arrow[hook]{d} & \s{C}_{\beta'}\arrow[]{d} \\
\Delta^n \arrow[]{r}{\sigma_{\beta'}} &  \s{C}_{\beta'+1}
\end{tikzcd}\]
with $\sigma_{\beta'} |_{\Delta^{\{n \}}} = v$.  Suppose $f \in \tx{Fun}'(\s{C}_{\beta'} ,\s{D})$, so $f\circ \sigma_{\beta'} |_{\Delta^{\{n\}}}$ is some terminal object of $\s{D}$, say $z\in\s{D}$.  Then we have a pullback square 
\[\begin{tikzcd}
\tx{Fun}' (\s{C}_{\beta'+1} , \s{D}) \times_{\tx{Fun}' (\s{C}_{\beta'} , \s{D})} \{f \} \arrow[]{r} \arrow[]{d} & \tx{Fun}(\Delta^n , \s{D})\times_{\tx{Fun}(\Delta^{\{n \}} , \s{D}) } \{z\} \arrow[]{d} \\ 
\Delta^0 \arrow[hook]{r}{f\circ ( \sigma_{\beta'}|_{\partial \Delta^n})}& \tx{Fun}(\partial\Delta^n , \s{D})\times_{\tx{Fun}(\Delta^{\{n \}} , \s{D}) } \{z\},
\end{tikzcd}\]
so the fiber $\tx{Fun}' (\s{C}_{\beta'+1} , \s{D}) \times_{\tx{Fun}' (\s{C}_{\beta'} , \s{D})} \{f \}$ is a contractible Kan complex by (\cite{thesis}, 1.3.6).
\end{proof}
\begin{rmk}\label{termtoptrmk} Suppose $\s{C}\in \maybehat{\tx{Cat}}_\infty^\tx{term}$, let $\s{R}$ be the set of inclusions of terminal objects into $\s{C}$, $\{ \ol{p}_\alpha :  (\Delta^{-1})^{\triangleleft} = \Delta^0 \xrightarrow{c_\alpha}\s{C} ; c_\alpha \in\s{C} \tx{ is terminal} \}$. Let $\s{P}^{\{ \Delta^{-1} \}}_{\s{R}}(\s{C})$ be as in (\cite{htt}, 5.3.6.2), so $\s{P}^{\{ \Delta^{-1} \}}_{\s{R}}(\s{C})$ is defined as follows:
\begin{itemize}
\item Let $j_0 : \s{C} \to\s{P}(\s{C})$ be the Yoneda embedding, and let $\ul{\Delta}^{-1} \in \s{P}(\s{C})$ be an initial object of $\s{P}(\s{C})$.  For each terminal object $c_\alpha \in\s{C}$, pick a morphism $\ul{\Delta}^{-1} \to j_0 (c_\alpha)$ in $\s{P}(\s{C})$, and let $S$ be the set of these morphisms.  By (\cite{htt}, 5.5.4.15), we have a Bousfield localization $L: \s{P}(\s{C}) \rightleftarrows S^{-1}\s{P}(\s{C}) : i$,
where $i : S^{-1}\s{P}(\s{C}) \hookrightarrow\s{P}(\s{C})$ is the inclusion of the full subcategory spanned by the $S$-local objects of $\s{P}(\s{C})$.
\item The $\infty$-category $\s{P}^{\{ \Delta^{-1} \}}_{\s{R}}(\s{C})$ is defined to be smallest full subcategory of $S^{-1}\s{P}(\s{C})$ which contains the essential image of $L \circ j_0$ and admits an initial object.
\end{itemize}
Note that an object $(Z : \s{C}^\tx{op} \to \s{S}) \in \s{P}(\s{C})$ is $S$-local if and only if for every (or, equivalently, any) terminal object $c_\alpha \in \s{C}$, the Kan complex $Z(c_\alpha)$ is contractible.  Indeed, for any $Z \in\s{P}(\s{C})$ and any map $s : \ul{\Delta}^{-1} \to j_0 (c_\alpha)$ in $S$, we have 
\[ Z(c_\alpha) \simeq \tx{Map}_{\s{P}(\s{C})} (j_0 (c_\alpha) , Z) \xrightarrow{-\circ s} \tx{Map}_{\s{P}(\s{C})} (\ul{\Delta}^{-1} , Z) \simeq *, \]
so $ \tx{Map}_{\s{P}(\s{C})} (j_0 (c_\alpha) , Z) \xrightarrow{-\circ s} \tx{Map}_{\s{P}(\s{C})} (\ul{\Delta}^{-1} , Z)$ is a weak equivalence if and only if $ Z(c_\alpha) \simeq *$.  Thus we have that:\begin{enumerate}
\item\label{termtoptrmk1} The object $j_0 (c_\alpha) \in \s{P}(\s{C})$ is $S$-local for any terminal object $c_\alpha \in\s{C}$, so there is an equivalence $j_0 (c_\alpha)\simeq L\circ j_0 (c_\alpha)$ in $\s{P}(\s{C})$,
\item\label{termtoptrmk2} The object $j_0 (c_\alpha) \in S^{-1} \s{P}(\s{C})$ is initial for any terminal object $c_\alpha \in\s{C}$, and
\item\label{termtoptrmk3} The $\infty$-category $\s{C}$ has a zero object if and only if $j_0 (x)$ is $S$-local for every $x\in\s{C} $.
\end{enumerate}
Suppose $x\in\s{C}$, and let $c_\alpha \in \s{C}$ be terminal.  By item~\ref{termtoptrmk1}, we have weak equivalences of Kan complexes
\begin{align*} \tx{Map}_{\s{P}^{\{ \Delta^{-1} \}}_{\s{R}}(\s{C}) } (L\circ j_0 (x), L\circ j_0 (c_\alpha)) & \simeq  \tx{Map}_{\s{P}(\s{C}) } ( j_0 (x), L\circ j_0 (c_\alpha)) \\ 
& \simeq  \tx{Map}_{\s{P}(\s{C}) } ( j_0 (x),  j_0 (c_\alpha)) \simeq \tx{Map}_{\s{C}}(x,c_\alpha) \simeq *.
\end{align*}
By item~\ref{termtoptrmk2}, $\s{P}^{\{ \Delta^{-1} \}}_{\s{R}}(\s{C}) \subseteq S^{-1} \s{P}(\s{C})$ is the smallest full subcategory containing the essential image of $L\circ j_0$.  In particular, this means $\s{P}^{\{ \Delta^{-1} \}}_{\s{R}}(\s{C})$ is essentially small.

Thus for any terminal object $c_\alpha \in \s{C}$, the object $j_0 (c_\alpha) \in \s{P}^{\{ \Delta^{-1} \}}_{\s{R}}(\s{C}) $ is a zero object, so $\s{P}^{\{ \Delta^{-1} \}}_{\s{R}}(\s{C}) \in \maybehat{\tx{Cat}}_\infty^\tx{pt}$.

Let $j : \s{C} \to \s{P}^{\{ \Delta^{-1} \}}_{\s{R}}(\s{C})$ be $L\circ j_0$.  By the above, $j$ is a morphism in $\maybehat{\tx{Cat}}_\infty^\tx{term}$.
\end{rmk}
\begin{proposition}\label{termptadj} For $\s{C} \in  \maybehat{\tx{Cat}}_\infty^\tx{term}$, let $\s{R}$ be the collection of inclusions of terminal objects into $\s{C}$, $\{ \ol{p}_\alpha :  (\Delta^{-1})^{\triangleleft} = \Delta^0 \xrightarrow{c_\alpha}\s{C} ; c_\alpha \in\s{C} \tx{ is terminal} \}$. The inclusion $\maybehat{\tx{Cat}}_\infty^\tx{pt} \subseteq \maybehat{\tx{Cat}}_\infty^\tx{term}$ admits a left adjoint $L: \maybehat{\tx{Cat}}_\infty^\tx{term} \to \maybehat{\tx{Cat}}_\infty^\tx{pt}$, with $L(\s{C}) \simeq \s{P}^{\{ \Delta^{-1} \}}_{\s{R}}(\s{C}) $, where $\s{P}^{\{ \Delta^{-1} \}}_{\s{R}}(\s{C})$ is as in (\cite{htt}, 5.3.6.2).
\end{proposition}
\begin{proof} By (\cite{htt}, 5.2.7.8) and Remark~\ref{termtoptrmk}, it suffices to show that for each $\s{C}\in \maybehat{\tx{Cat}}_\infty^\tx{term}$, the map $j : \s{C} \to \s{P}^{\{ \Delta^{-1} \}}_{\s{R}}(\s{C})$ from Remark~\ref{termtoptrmk} exhibits $\s{P}^{\{ \Delta^{-1} \}}_{\s{R}}(\s{C}) \in \maybehat{\tx{Cat}}_\infty^\tx{pt}$ as a $\maybehat{\tx{Cat}}_\infty^\tx{pt}$-localization of $\s{C}$ (\cite{htt}, Definition 5.2.7.6).

Suppose $\s{E} \in \maybehat{\tx{Cat}}_\infty^\tx{pt}$.  By (\cite{htt}, 5.3.6.2), precomposition by $j$ induces an equivalence of $\infty$-categories
\[ \tx{Fun}_{\{ \Delta^{-1}\}} (\s{P}^{\{ \Delta^{-1} \}}_{\s{R}}(\s{C}), \s{E}) \xrightarrow{\sim} \tx{Fun}_{\s{R}} (\s{C}, \s{E}), \] 
where $\tx{Fun}_{\{ \Delta^{-1}\}} (\s{P}^{\{ \Delta^{-1} \}}_{\s{R}}(\s{C}), \s{E})\subseteq \tx{Fun} (\s{P}^{\{ \Delta^{-1} \}}_{\s{R}}(\s{C}), \s{E}) $ is the full subcategory spanned by the functors $\s{P}^{\{ \Delta^{-1} \}}_{\s{R}}(\s{C})\to \s{E}$ that take initial (hence, terminal) objects to initial (hence, terminal) objects, and $ \tx{Fun}_{\s{R}} (\s{C}, \s{E}) \subseteq  \tx{Fun} (\s{C}, \s{E})$ is the full subcategory spanned by the functors $\s{C} \to \s{E}$ that take terminal objects of $\s{C}$ to initial (hence, terminal) objects of $\s{E}$.

Thus $j$ exhibits $\s{P}^{\{ \Delta^{-1} \}}_{\s{R}}(\s{C}) \in \maybehat{\tx{Cat}}_\infty^\tx{pt}$ as a $\maybehat{\tx{Cat}}_\infty^\tx{pt}$-localization of $\s{C}$ by (\cite{htt}, 3.1.3.2).
\end{proof}
\begin{rmk}\label{ptunptadjunctionremark} By (\cite{htt}, 5.2.2.6), the two adjunctions from Propositions~\ref{catptadj} and~\ref{termptadj} compose to give an adjunction of $\infty$-categories \[\begin{tikzcd}\maybehat{\tx{Cat}}_\infty \arrow[bend left]{rr}{\s{C} \mapsto \s{P}^{\{ \Delta^{-1} \}}_{\s{R}}(\s{C}^\triangleright)} & \perp & \maybehat{\tx{Cat}}_\infty^\tx{pt} \arrow[hook, bend left]{ll} \end{tikzcd}\]
where $\s{R} = \{ (\Delta^{-1})^\triangleleft = \Delta^0 \xrightarrow{v} \s{C}^\triangleright \}$ is a singleton set containing the inclusion of the cone point $v \in \s{C}^\triangleright$.  

Let's describe $ \s{P}^{\{ \Delta^{-1} \}}_{\s{R}}(\s{C}^\triangleright)$ a bit more explicitly.  Suppose $\s{C} \in \maybehat{\tx{Cat}}_\infty$, let $\s{B}:= \s{C}^\tx{op}$, write $v$ for the cone point of $\s{B}^\triangleleft$, and let $S$ be the singleton set containing a morphism $\ul{\Delta}^{-1} \to V$, where $V : \s{B}^\triangleleft \to \s{S}$ is corepresented by $v$, and $\ul{\Delta}^{-1} \in \tx{Fun}(\s{B}^\triangleleft , \s{S})$ is the constant functor to $\Delta^{-1}$.  Let $S^{-1} \tx{Fun}(\s{B}^\triangleleft , \s{S}) \subseteq \tx{Fun}(\s{B}^\triangleleft , \s{S}) = \s{P}(\s{C}^\triangleright)$ be the full subcategory spanned by the $S$-local objects.  As in Remark~\ref{termtoptrmk}, $Z : \s{B}^\triangleleft \to \s{S}$ is $S$-local if and only if $Z(v)$ is contractible.

For any $b \in \s{B}^\triangleleft $, $V(b) \simeq \tx{Map}_{\s{B}^\triangleleft} (v,b) \simeq *$ in $\s{S}$, so $V \in \tx{Fun}(\s{B}^\triangleleft , \s{S})$ is terminal (\cite{rezk}, 31.9), as is $\ul{\Delta}^0$, so $\ul{\Delta}^{0} \simeq V$.  On the other hand, for any $Z\in S^{-1} \tx{Fun}(\s{B}^\triangleleft , \s{S})$ we have weak equivalences of Kan complexes
\[ \tx{Map}_{S^{-1} \tx{Fun}(\s{B}^\triangleleft , \s{S})} (V,Z) = \tx{Map}_{\tx{Fun}(\s{B}^\triangleleft , \s{S})} (V,Z) \simeq Z(v) \simeq \Delta^0, \]
so $\ul{\Delta}^{0} \simeq V$ in $S^{-1} \tx{Fun}(\s{B}^\triangleleft , \s{S})$ is initial.  

Now suppose $Y \in \tx{Fun}( \s{B}^\triangleleft , \s{S})$ is any functor such that $Y(v) \simeq \Delta^{-1}$, and consider the map $Y \to Y \coprod \ul{\Delta}^0$ in $\tx{Fun}(\s{B}^\triangleleft , \s{S})$.  This map exhibits $Y \coprod \ul{\Delta}^0$ as an $S^{-1} \tx{Fun}(\s{B}^\triangleleft , \s{S})$-localization of $Y$.  Indeed, $Y \coprod \ul{\Delta}^0$ is $S$-local, and for any $Z \in S^{-1} \tx{Fun}(\s{B}^\triangleleft , \s{S})$ we have (by \cite{htt}, \S 4.2.4) a (homotopy) pullback square 
\[ \begin{tikzcd}
\tx{Map}_{\tx{Fun}( \s{B}^\triangleleft , \s{S})} (Y \coprod \ul{\Delta}^0 , Z) \arrow[]{r} \arrow[]{d} & \tx{Map}_{\tx{Fun}( \s{B}^\triangleleft , \s{S})} (Y,Z) \arrow[]{d} \\ 
\tx{Map}_{\tx{Fun}( \s{B}^\triangleleft , \s{S})} (\ul{\Delta}^0 , Z) \arrow[]{r} & \tx{Map}_{\tx{Fun}( \s{B}^\triangleleft , \s{S})} (\ul{\Delta}^{-1}, Z)
\end{tikzcd}\]
in $\s{S}$, in which the bottom horizontal is an equivalence, and the top horizontal is induced by precomposition by $Y \to Y \coprod \ul{\Delta}^0$.  In particular, this means that if $b \in \s{B} \subseteq \s{B}^\triangleleft$ and $B : \s{B}^\triangleleft \to\s{S}$ is corepresented by $b$, then $B \to B \coprod \ul{\Delta}^0$ exhibits $B \coprod \ul{\Delta}^0$ as an $S^{-1} \tx{Fun}(\s{B}^\triangleleft , \s{S})$-localization of $B$.

Letting $j_0 : \s{C}^\triangleright \to\s{P}(\s{C}^\triangleright )$ be the Yoneda embedding, we see by Remark~\ref{termtoptrmk} that $ \s{P}^{\{ \Delta^{-1} \}}_{\s{R}}(\s{C}^\triangleright) \subseteq \s{P}(\s{C}^\triangleright)$ is the full subcategory spanned by the objects of $\s{P}(\s{C}^\triangleright)$ which are isomorphic to either $j_0 (v)$ or $j_0 (c) \coprod \ul{\Delta}^0$ for some object $c \in \s{C} \subseteq \s{C}^\triangleright$.
\end{rmk}
Our next goal is to compare $\s{P}^{\{ \Delta^{-1} \}}_{\s{R}}(\s{C}^\triangleright)$ of Remark~\ref{ptunptadjunctionremark} with $\s{C}_\diamond$.  To do this, we'll show that the natural inclusion $\zeta_\s{C} : \s{C}^\triangleright \hookrightarrow \s{C}_\diamond$ of Construction~\ref{triangletodiamondinclusion} induces a natural transformation $\Delta^1 \to \tx{Fun}(\s{C}^\triangleright , \s{P}(\s{C}^\triangleright))$ 
 from the Yoneda embedding $j_0 : \s{C}^\triangleright \to \s{P}(\s{C}^\triangleright)$ to the composite functor \[ \s{C}^\triangleright \xrightarrow{\zeta_\s{C}} \s{C}_\diamond \xrightarrow{j_0'} \s{P}(\s{C}_\diamond) \xrightarrow{(\zeta_\s{C}^\tx{op})^*} \s{P}(\s{C}^\triangleright),\] where $j_0' : \s{C}_\diamond \rightarrow \s{P}(\s{C}_\diamond)$ is the Yoneda embedding.  Rather than construct this natural transformation via straightening and unstraightening, which is somewhat harder, we will use the relative join construction of (\cite{kerodon}, 0241) to produce a covariant transport representation for it.  We first review the relative join construction:
\begin{definition}[\tx{}\cite{kerodon}, 0241]\label{reljoindef} Let $\s{E}$ be a simplicial set, and let $\rho = \rho_\s{E} : \Delta^1 \times \s{E} \to\s{E} \star \s{E}$ be the unique map such that the diagram
\[\begin{tikzcd} \Delta^{\{ 0 \}} \times \s{E} \arrow[hook]{r} \arrow[]{d}{\tx{id}_\s{E}} & \Delta^1 \times \s{E} \arrow[]{d}{\rho} &\Delta^{\{ 0 \}} \times \s{E} \arrow[hook]{l} \arrow[]{d}{\tx{id}_\s{E}} \\  
\s{E}\star \Delta^{-1} \arrow[hook]{r} & \s{E} \star\s{E} & \Delta^{-1} \star\s{E} \arrow[hook]{l}  \end{tikzcd}\]
is commutative.  For any cospan of simplicial sets $\s{C} \xrightarrow{f} \s{E} \xleftarrow{g} \s{D}$, let $\s{C}\star_{\s{E}} \s{D}$ denote the fiber product $(\s{C} \star \s{D}) \times_{\s{E}\star\s{E}} (\Delta^1 \times \s{E})$, so we have a pullback square
\[\begin{tikzcd}
\s{C}\star_{\s{E}} \s{D} \arrow[]{r} \arrow[]{d} & \s{C}\star\s{D} \arrow[]{d}{f\star g} \\
\Delta^1 \times \s{E} \arrow[]{r}{\rho} & \s{E} \star \s{E}.
\end{tikzcd}\]
We call $\s{C}\star_{\s{E}} \s{D}$ the \emph{join of $\s{C}$ and $\s{D}$ relative to $\s{E}$}.
\end{definition}
\begin{lemma}[\tx{}\cite{thesis}, 1.4.16]\label{fiberrelativejoins} Suppose a commutative diagram of simplicial sets
\begin{equation}\label{lambdagexamplesquare} \begin{tikzcd}
& \tilde{X} \arrow[]{r}{\tilde{g}} \arrow[]{d}[swap]{\lambda} & \tilde{Y} \arrow[]{d}{\lambda'} & \tilde{Y} \arrow[]{d}{\lambda'} \arrow[]{l}[swap]{\tx{id}} \\
\Delta^0 \arrow[]{r}{x} & X \arrow[]{r}{g} & Y & Y \arrow[]{l}[swap]{\tx{id}}
\end{tikzcd}\end{equation}
and let $\tilde{X}_x := \tilde{X} \times_X \{ x\}$, $\tilde{Y}_{g(x)} := \tilde{Y} \times_Y \{ g(x)\}$.  Then we have a pullback square 
\[\begin{tikzcd} 
\tilde{X}_x \star_{\tilde{Y}_{g(x)}} \tilde{Y}_{g(x)} \arrow[]{d} \arrow[]{r} & \tilde{X} \star_{\tilde{Y}} \tilde{Y} \arrow[]{d}{H'} \\
\Delta^0 \star_{\Delta^0} \Delta^0 \arrow[]{r} & X \star_Y Y,
\end{tikzcd}\]
where $H'$ is the morphism of relative joins induced by the morphism of cospans in (\ref{lambdagexamplesquare}), and the lower horizontal $\Delta^1 \cong \Delta^0 \star_{\Delta^0} \Delta^0 \to X \star_Y Y$ is identified with $\tx{id}_{g(x)}$ under the isomorphism of simplicial sets $\tx{Hom}^\tx{R}_{X\star_Y Y} (x, g(x)) \cong \tx{Hom}^\tx{R}_{Y} (g(x), g(x))$. 
\end{lemma}
\begin{construction}\label{upsilonconstr} Suppose a functor of $\infty$-categories $f: \s{A}\to\s{B}$.  We construct a left fibration $\Upsilon_f$ over $\Delta^1 \times \s{A}^\tx{op} \times \s{A}$ as follows: 

\hfill

For a simplicial set  $K$, let $\tx{Tw}(K)$ be the twisted arrow construction of (\cite{kerodon}, 03JG),\footnote{Note that $\tx{Tw}(-)$ of (\cite{kerodon}, 03JG) is slightly different from $\tx{TwArr}(-)$ of (\cite{ha}, 5.2.1.1).  The difference of definitions is because \cite{kerodon} prioritizes coCartesian fibrations and covariant transport, while \cite{ha} and \cite{htt} prioritize Cartesian fibrations and contravariant transport.} i.e. the functor 
\[\begin{tikzcd} 
 \dop \arrow[]{r}{\tx{Tw}(K)} & \tx{Set}\\[-1.5em]
 J \arrow[mapsto]{r} & \tx{Hom}_{\tx{sSet}} \big( \tx{N}(J^\tx{op} \star J), K \big).
\end{tikzcd}\]
Let $\lambda = \lambda_K : \tx{Tw}(K) \to K^\tx{op} \times K$ be the left fibration of (\cite{kerodon}, 03JK), given by 
\[\begin{tikzcd} 
\tx{Hom}(\Delta^n , K) \arrow[]{r}{\lambda_K \circ -} & \tx{Hom}(\Delta^n , K^\tx{op} \times K) \\[-1.5em]
\big( (\Delta^n)^\tx{op} \star \Delta^n \xrightarrow{\sigma} K \big) \arrow[mapsto]{r} & \big( \sigma |_{(\Delta^n)^\tx{op}}, \tx{\,} \sigma |_{\Delta^n} \big).
\end{tikzcd}\]
By naturality of $\lambda$ in $K$, we have a commutative diagram of simplicial sets
\[\begin{tikzcd}
\tx{Tw}(\s{A}) \arrow[dashed]{r}{\gamma} \arrow[bend left]{rr}{\tx{Tw}(f)} \arrow[]{d}[swap]{\lambda_{\s{A}}} & \tx{Tw}(\s{B}) \times_{\s{B}^\tx{op} \times \s{B}} (\s{A}^\tx{op}\times \s{A}) \arrow[]{d}{\lambda'} \arrow[]{r} & \tx{Tw}(\s{B}) \arrow[]{d}{\lambda_{\s{B}}} \\
\s{A}^\tx{op}\times \s{A} \arrow[equal]{r} & \s{A}^\tx{op}\times \s{A} \arrow[]{r}[swap]{f^\tx{op} \times f} & \s{B}^\tx{op}\times \s{B}
\end{tikzcd}\]
in which the right-hand square is a pullback, and in which all vertical morphisms are left fibrations.  Applying the relative join construction to each of the rows in the commutative diagram
\[\begin{tikzcd}
\tx{Tw}(\s{A}) \arrow[]{r}{\gamma} \arrow[]{d}[swap]{\lambda_{\s{A}}} & \tx{Tw}(\s{B}) \times_{\s{B}^\tx{op} \times \s{B}} (\s{A}^\tx{op}\times \s{A}) \arrow[]{d}{\lambda'}  &  \tx{Tw}(\s{B}) \times_{\s{B}^\tx{op} \times \s{B}} (\s{A}^\tx{op}\times \s{A}) \arrow[]{d}{\lambda'} \arrow[]{l}[swap]{\tx{id}} \\
\s{A}^\tx{op}\times \s{A} \arrow[]{r}[swap]{\tx{id}} & \s{A}^\tx{op}\times \s{A} & \s{A}^\tx{op}\times \s{A} \arrow[]{l}{\tx{id}}
\end{tikzcd}\]
gives a morphism of relative joins, and we denote by $\Upsilon_f$ the composition
\[\begin{tikzcd} \tx{Tw}(\s{A}) \star_{\tx{Tw}(\s{B}) \times_{\s{B}^\tx{op} \times \s{B}} (\s{A}^\tx{op}\times \s{A})} \big(  \tx{Tw}(\s{B}) \times_{\s{B}^\tx{op} \times \s{B}} (\s{A}^\tx{op}\times \s{A}) \big) \arrow[]{d} \arrow[bend left=10]{dr}{:=\Upsilon_f} & \\ 
 (\s{A}^\tx{op}\times \s{A}) \star_{\s{A}^\tx{op}\times \s{A}} (\s{A}^\tx{op}\times \s{A}) \arrow[equal]{r}{\sim} & \Delta^1 \times (\s{A}^\tx{op}\times \s{A}) .\end{tikzcd}\]
Both $\lambda_\s{A}$ and $\lambda'$ are left fibrations, so $\Upsilon_f$ is a left fibration by (\cite{kerodon}, 02RH) and (\cite{htt}, 2.4.2.4). 
\end{construction}
\begin{rmk}\label{upsilonunwindrmk} Let $f:\s{A}\to \s{B}$ be a functor between $\infty$-categories.  Roughly speaking, $\Upsilon = \Upsilon_f$ ``gives the natural transformation $\Delta^1 \to \tx{Fun}(\s{A}^\tx{op} \times \s{A}, \s{S}) $ induced by $f$ from `the functor $(a,a') \mapsto \tx{Map}_\s{A} (a,a')$' to `the functor $(a,a') \mapsto \tx{Map}_{\s{B}} (f(a), f(a'))$'.''\footnote{It's possible to get this natural transformation directly via the straightening and unstraightening of \cite{htt}, but it's somewhat harder than using Construction~\ref{upsilonconstr}.}  Let's make this a bit more precise:

\hfill

Letting $\s{D}:=\tx{Tw}(\s{B}) \times_{\s{B}^\tx{op} \times \s{B}} (\s{A}^\tx{op}\times \s{A})$, we have a commutative diagram of simplicial sets
\[\begin{tikzcd}
\tx{Tw}(\s{A}) \arrow[hook]{r} \arrow[]{d}{\lambda_\s{A}} & \tx{Tw}(\s{A}) \star_{\s{D}} \s{D} \arrow[]{d}{\Upsilon} & \tx{Tw}(\s{B}) \times_{\s{B}^\tx{op} \times \s{B}} (\s{A}^\tx{op}\times \s{A}) \arrow[hook]{l} \arrow[]{d}{\lambda'} \\ 
\Delta^{\{0\}} \times \s{A}^\tx{op}\times \s{A} \arrow[hook]{r} & \Delta^{1} \times \s{A}^\tx{op}\times \s{A} & \Delta^{\{1\}} \times \s{A}^\tx{op}\times \s{A} \arrow[hook]{l}
\end{tikzcd}\]
in which both squares are pullbacks, and $\lambda_{\s{A}}$, $\lambda'$ are as in Construction~\ref{upsilonconstr}. 

Any covariant transport representation $\c{H}_\s{A} : \s{A}^\tx{op} \times \s{A}\to\s{S}$ of the left fibration $\lambda_{\s{A}}$  (in the sense of \cite{kerodon}, 028V) is a $\tx{Hom}$-functor for $\s{A}$ (in the sense of \cite{kerodon}, 046B).  Likewise, by (\cite{ha}, 5.2.1.11), for any commutative diagram
\[ \begin{tikzcd}
\tx{Tw}(\s{A}) \arrow[]{r} \arrow[]{d}[swap]{\lambda_\s{A}} & \s{Z}^\tx{op} \arrow[]{d}{q^\tx{op}} \\
 \s{A}^\tx{op} \times \s{A} \arrow[]{r}{\chi} & \tx{Cat}_\infty
\end{tikzcd}\]
wherein $q : \s{Z} \to \tx{Cat}_\infty^\tx{op}$ is the universal Cartesian fibration of (\cite{htt}, \S 3.3.2) and the induced map $\tx{Tw}(\s{A}) \to  (\s{A}^\tx{op} \times \s{A} )_{\tx{Cat}_\infty} \s{Z}^\tx{op}$ is an equivalence of coCartesian fibrations over $\s{A}^\tx{op} \times \s{A}$, the map $\s{A} \to \s{P}(\s{A})$ determined by $\chi$ is homotopic to the Yoneda embedding $\s{A} \to \s{P}(\s{A})$.  By the same token, if $\c{H}_{\s{B}} : \s{B}^\tx{op} \times \s{B} \to \s{S}$ is a $\tx{Hom}$-functor for $\s{B}$, then any covariant transport representation $\c{H}' : \s{A}^\tx{op} \times \s{A} \to \s{S}$ of $\lambda'$ is isomorphic to $(f^\tx{op} \times f)\circ \c{H}_{\s{B}}$ in the $\infty$-category $\tx{Fun}(\s{A}^\tx{op} \times \s{A}, \s{S})$.  In particular, for any $(a,a') \in \s{A}^\tx{op} \times \s{A}$, we have isomorphisms in $\tx{h}\tx{Kan}$
\[ \tx{Tw}(\s{A}) \times_{\s{A}^\tx{op} \times \s{A}} \{ (a,a') \} \simeq \tx{Map}_{\s{A}} (a,a') \tx{\; and \;} \s{D} \times_{\s{A}^\tx{op} \times \s{A}} \{ (a,a') \} \simeq \tx{Map}_{\s{B}} (f(a),f(a')) . \]
Next suppose a vertex $(a,a') \in \s{A}^\tx{op} \times\s{A}$, and consider the pullback 
\[ \Upsilon_{a,a'} : \big( \tx{Tw}(\s{A}) \star_{\s{D}} \s{D} \big) \times_{\Delta^{1} \times \s{A}^\tx{op}\times \s{A} } \big( \Delta^{1} \times \{ (a,a') \} \big) \to \Delta^1 \times \{ (a,a') \} \cong \Delta^1 \]
of $\Upsilon$.  Taking the diagram
\[\begin{tikzcd}
& \tx{Tw}(\s{A}) \arrow[]{r}{\gamma} \arrow[]{d}[swap]{\lambda_{\s{A}}} & \tx{Tw}(\s{B}) \times_{\s{B}^\tx{op} \times \s{B}} (\s{A}^\tx{op}\times \s{A}) \arrow[]{d}{\lambda'} & \tx{Tw}(\s{B}) \times_{\s{B}^\tx{op} \times \s{B}} (\s{A}^\tx{op}\times \s{A}) \arrow[]{d}{\lambda'} \arrow[]{l}[swap]{\tx{id}}  \\
\Delta^0 \arrow[]{r}{(a,a')} & \s{A}^\tx{op}\times \s{A} \arrow[]{r}{\tx{id}} & \s{A}^\tx{op}\times \s{A} &  \s{A}^\tx{op}\times \s{A} \arrow[]{l}[swap]{\tx{id}},
\end{tikzcd}\]
for (\ref{lambdagexamplesquare}) in Lemma~\ref{fiberrelativejoins}, we can identify the pullback $\Upsilon_{a,a'}$ of $\Upsilon$ with the morphism of relative joins induced by the diagram
\[\begin{tikzcd}
\tx{Tw}(\s{A})_{(a,a')} \arrow[]{rr}{\tx{Tw}(f)_{(a,a')}} \arrow[]{d} & & \tx{Tw}(\s{B})_{ (f(a),f(a'))} \arrow[]{d} &  \tx{Tw}(\s{B})_{ (f(a),f(a'))} \arrow[]{l}[swap]{\tx{id}} \arrow[]{d} \\
\Delta^0 \arrow[]{rr} & & \Delta^0 & \Delta^0 .\arrow[]{l}
\end{tikzcd}\]
By (\cite{kerodon}, 024J), $\Upsilon_{a,a'}$ is a coCartesian fibration (moreover, it is a left fibration), and $\tx{Tw}(f)_{(a,a')} : \Delta^1 \to \s{S} \hookrightarrow \tx{Cat}_\infty$ is a covariant transport representation of $\Upsilon_{a,a'}$.

By (\cite{kerodon}, 03JY), there is a morphism $\iota_{a} : \s{A}_{a/} \to \tx{Tw}(\s{A})$, defined as follows: An $n$-simplex of $\s{A}_{a/}$ can be identified with a map $\sigma : \{ v \} \star \Delta^n \to \s{C}$ with $\sigma (v) =a$, and composition with $\iota_{a}$ takes this $n$-simplex of $\s{A}_{a/}$ to the $n$-simplex $(\Delta^n )^{\tx{op}} \star \Delta^n \xrightarrow{! \star \tx{id}} \{ v \} \star \Delta^n \xrightarrow{\sigma} \s{C}$ of $\tx{Tw}(\s{A})$.  The diagram
\[\begin{tikzcd}
\s{A}_{a/} \arrow[]{r}{\iota_{a}} \arrow[]{d} & \tx{Tw}(\s{A})\arrow[]{d}{\tx{pr}_1 \circ \lambda_{\s{A}}} \\
\{a\} \arrow[hook]{r} & \s{A}^\tx{op}
\end{tikzcd}\]
is commutative, so $\iota_a$ factors as $\s{A}_{a/} \to  \tx{Tw}(\s{A}) \times_{\s{A}^{\tx{op}}} \{a \} \hookrightarrow \tx{Tw}(\s{A})$.  Moreover the diagram
\[\begin{tikzcd}
\s{A}_{a/} \arrow[]{rr}{\iota_{a}} \arrow[]{dr} & &  \tx{Tw}(\s{A}) \times_{\s{A}^{\tx{op}}} \{a \} \arrow[]{dl}{\tx{pr}_2 \circ \lambda_{\s{A}} |_{\tx{Tw}(\s{A}) \times_{\s{A}^{\tx{op}}} \{a \}} } \\
& \s{A} &
\end{tikzcd}\]
is an equivalence of left fibrations over $\s{A}$ (\cite{kerodon}, 03K0).  By naturality of the $\iota$ construction, we have a commutative diagram of Kan complexes
\[\begin{tikzcd}
\tx{Hom}_{\s{A}}^\tx{L} (a,a' ) \arrow[]{d}[swap]{f\circ -} \arrow[equal]{r} & \s{A}_{a/} \times_{\s{A}} \{ a' \} \arrow[]{rr}{(\iota_a )_{a'}} & & \tx{Tw}(\s{A})\times_{\s{A}^\tx{op} \times \s{A}} \{ (a,a') \} \arrow[]{d}{\tx{Tw}(f)_{(a,a')}} \\
\tx{Hom}_{\s{B}}^\tx{L} (f(a), f(a') ) \arrow[equal]{r} & \s{B}_{f(a)/} \times_{\s{B}} \{ f(a') \} \arrow[]{rr}{(\iota_{f(a)} )_{f(a')}} & &  \tx{Tw}(\s{B})\times_{\s{B}^\tx{op} \times \s{B}} \{ (f(a),f(a')) \}
\end{tikzcd}\]
in which the horizontal morphisms are weak equivalences of Kan complexes.
\end{rmk}
If $\s{C}$ is an $\infty$-category, then applying Construction~\ref{upsilonconstr} to the inclusion $\zeta : \s{C}^{\triangleright} \hookrightarrow \s{C}_\diamond$ of Construction~\ref{triangletodiamondinclusion} gives a left fibration $\Upsilon_\zeta$ over $\Delta^1 \times (\s{C}^\triangleright )^\tx{op} \times \s{C}^\triangleright$.  Roughly, $\Upsilon_\zeta$ ``gives us a natural transformation $\Delta^1 \to \tx{Fun}( \s{C}^\triangleright, \s{P}(\s{C}^\triangleright ) )$ from the Yoneda embedding `$c \mapsto \tx{Map}_{\s{C}^\triangleright} (-,c)$' to a functor `$c \mapsto  \big( \tx{Map}_{\s{C}_\diamond} (-,c) \big) \big|_{ (\s{C}^\triangleright )^\tx{op}}$'.''  The next proposition says that the components of this natural transformation, which are morphisms in $\s{P}(\s{C}^\triangleright ) $, are localizations relative to the subcategory $S^{-1} \s{P}(\s{C}^\triangleright )  \subseteq \s{P}(\s{C}^\triangleright ) $ from Remark~\ref{ptunptadjunctionremark}.  
\begin{proposition}\label{zetaobjectwiselocalization} Suppose $\s{C}$ is an $\infty$-category, and let $\zeta : \s{C}^{\triangleright} \hookrightarrow \s{C}_\diamond$ be the natural inclusion of Construction~\ref{triangletodiamondinclusion}.   Let $\s{D} := \tx{Tw}(\s{C}_\diamond) \times_{ (\s{C}_\diamond )^\tx{op} \times \s{C}_\diamond} \big( (\s{C}^\triangleright )^\tx{op}\times \s{C}^\triangleright \big)$ and 
\[\Upsilon = \Upsilon_{\zeta} : \tx{Tw}(\s{C}^\triangleright )\star_{\s{D}} \s{D} \to \Delta^1 \times (\s{C}^\triangleright )^\tx{op} \times \s{C}^\triangleright \] 
be as in Construction~\ref{upsilonconstr}, so in particular $\Upsilon$ is a left fibration.  Let $\s{F} : \Delta^1 \times (\s{C}^\triangleright )^\tx{op} \times \s{C}^\triangleright \to \s{S}$ be a covariant transport representation of $\Upsilon$, so there is a commutative diagram
\[\begin{tikzcd}
\tx{Tw}(\s{C}^\triangleright )\star_{\s{D}} \s{D}\arrow[]{d}[swap]{\Upsilon} \arrow[]{rr}  & & \s{{QC}}_{\tx{Obj}} \arrow[]{d}{V} 
\\
 \Delta^1 \times (\s{C}^\triangleright )^\tx{op} \times \s{C}^\triangleright \arrow[]{r}{\s{F}} & \s{S} \arrow[hook]{r} & \tx{Cat}_\infty
\end{tikzcd}\]
where the induced map $\tx{Tw}(\s{C}^\triangleright )\star_{\s{D}} \s{D} \to \big(  \Delta^1 \times (\s{C}^\triangleright )^\tx{op} \times \s{C}^\triangleright \big) \times_{\tx{Cat}_\infty} \s{{QC}}_{\tx{Obj}} $
 is an equivalence of coCartesian (more specifically, left) fibrations over $ \Delta^1 \times (\s{C}^\triangleright )^\tx{op} \times \s{C}^\triangleright$.  

Let $F : \Delta^1 \times \s{C}^\triangleright \to \s{P}(\s{C}^\triangleright )$ be the currying of $\s{F}$.  Let $S \subseteq \tx{Hom}(\Delta^1 , \s{P}(\s{C}^\triangleright ))$ be the singleton set $S = \{ \ul{\Delta}^{-1} \to j_0 (v) \}$, where $j_0 : \s{C}^\triangleright \to \s{P}(\s{C}^\triangleright)$ is the Yoneda embedding, $v \in \s{C}^\triangleright$ is the cone point, and $\ul{\Delta}^{-1}$ is the constant functor at $\Delta^{-1} \in \s{S}$ (cf. Remark~\ref{ptunptadjunctionremark}).  

Then for each $c\in \s{C}^\triangleright$, the edge $F|_{\Delta^1 \times \{ c\}}$ of $\s{P}(\s{C}^\triangleright)$ is a localization relative to the reflective subcategory $S^{-1} \s{P}(\s{C}^\triangleright) \subseteq \s{P}(\s{C}^\triangleright)$.
\end{proposition}
\begin{proof} Let $c$ be a vertex of $\s{C}^\triangleright$. By Remark~\ref{termtoptrmk}, to show that the object $F(1,c) \in \s{P}(\s{C}^\triangleright )$ is in $S^{-1} \s{P}(\s{C}^\triangleright)$, it suffices to show that $F(1,c)(v) \in \s{S}$ is terminal.  As we saw in Remark~\ref{upsilonunwindrmk}, we have an isomorphism 
\[ | \big( \tx{Tw}(\s{C}^\triangleright )\star_{\s{D}} \s{D} \big) \times_{\Delta^1 \times (\s{C}^\triangleright )^\tx{op} \times \s{C}^\triangleright } \{ (1,v,c) \} | \simeq \tx{Map}_{\s{C}_\diamond } (\zeta (v), \zeta (c)) \]
in $\s{H}$, and $\tx{Map}_{\s{C}_\diamond } (\zeta (v), \zeta (c))$ is contractible by Corollary~\ref{diamondhaszerocor}.  Thus, by (\cite{kerodon}, 0290), $F(1,c)(v) = \s{F}(1,v,c) \in \s{S}$ is terminal.

Next suppose $Z \in S^{-1} \s{P}(\s{C}^\triangleright )$, so $Z(v) \simeq *$.  We need to show that composition with $F|_{\Delta^1 \times \{ c\}}$ induces an isomorphism
\[ \tx{Map}_{\s{P}(\s{C}^\triangleright )} \big( F(1,c), Z \big)  \xrightarrow{\sim} \tx{Map}_{\s{P}(\s{C}^\triangleright )} \big(   F(0,c) , Z \big)   \]
in $\s{H}$.  As we saw in Remark~\ref{upsilonunwindrmk}, $F(0,c) \in \s{P}(\s{C}^\triangleright )$ is isomorphic to $j_0 (c)$, where $j_0$ is the Yoneda embedding $\s{C}^\triangleright  \to \s{P}(\s{C}^\triangleright )$, so $c$ represents $F(0,c)$.

\hfill

First suppose $c$ is the cone point of $\s{C}^\triangleright$, which we denote by $v$.  By the Yoneda lemma, the Kan complex $\tx{Map}_{\s{P}(\s{C}^\triangleright )} \big(   F(0,v) , Z \big)$ is contractible, so it remains to show $\tx{Map}_{\s{P}(\s{C}^\triangleright )} \big( F(1,v), Z \big)$ is contractible.  The constant functor $\ul{\Delta}^0 \in S^{-1} \s{P}(\s{C}^\triangleright ) \subseteq \s{P}(\s{C}^\triangleright )$ is initial in $ S^{-1} \s{P}(\s{C}^\triangleright )$, so it suffices to show that $F(1,v)$ is isomorphic to $\ul{\Delta}^0$, i.e. terminal in $\s{P}(\s{C}^\triangleright )$.  

For any $c' \in \s{C}^\triangleright$, we have by Remark~\ref{ptunptadjunctionremark} and (\cite{kerodon}, 0290) that $F(1,v)(c') = \s{F}(1,c' ,v)$ is isomorphic to $\tx{Map}_{\s{C}_\diamond} (c' ,v)$, considered as vertices of $\s{S}$.  Thus for each $c' \in (\s{C}^\triangleright )^\tx{op}$, the vertex $F(1,v)(c') \in \s{S}$ is terminal by Corollary~\ref{diamondhaszerocor}.  Thus $F(1,v) \in \s{P}(\s{C}^\triangleright )$ is terminal by (\cite{rezk}, 31.9).

\hfill

Next suppose $c\in \s{C} \subseteq \s{C}^\triangleright$.  Let $t : \Delta^1 \to \s{C}_\diamond$ be the unique edge from $v$ to $c$, and let $\lambda_{\s{C}_\diamond} : \tx{Tw}(\s{C}_\diamond ) \to \s{C}_\diamond^\tx{op} \times \s{C}_\diamond$ be the left fibration of (\cite{kerodon}, 03JK).  Denote by $\lambda''$ the pullback of $\lambda_{\s{C}_\diamond}$ along $\zeta^\tx{op} \times t$ and by $\lambda'''$ the pullback of $\lambda_{\s{C}_\diamond}$ along $\zeta^\tx{op} \times v$:
\[\begin{tikzcd}
 \tx{Tw}(\s{C}_\diamond )\times_{\s{C}_\diamond^\tx{op} \times \s{C}_\diamond} \big( (\s{C}^\triangleright )^\tx{op} \times \{ v\} \big) \arrow[hook]{r} \arrow[]{d}[swap]{\lambda'''} & \tx{Tw}(\s{C}_\diamond )\times_{\s{C}_\diamond^\tx{op} \times \s{C}_\diamond} \big( (\s{C}^\triangleright )^\tx{op} \times \Delta^1 \big) \arrow[hook]{r} \arrow[]{d}[swap]{\lambda''} & \tx{Tw}(\s{C}_\diamond ) \arrow[]{d}{\lambda_{\s{C}_\diamond}} \\
(\s{C}^\triangleright )^\tx{op} \times \Delta^{\{0\}} \arrow[hook]{r} & (\s{C}^\triangleright )^\tx{op} \times \Delta^1 \arrow[hook]{r}[swap]{\zeta^\tx{op} \times t} & \s{C}_\diamond^\tx{op} \times \s{C}_\diamond .
\end{tikzcd}\]
Let $\s{F}' : (\s{C}^\triangleright )^\tx{op} \times \Delta^1 \to \s{S}$ be a covariant transport representation of the left fibration $\lambda''$, and let $z : \Delta^1 \to \s{P}(\s{C}^\triangleright )$ be the currying of $\s{F}'$.  The pullback of $\lambda_{\s{C}_\diamond }$ along $\zeta^\tx{op} \times c$ is isomorphic, as an object of the ordinary slice category $\tx{sSet}_{/(\s{C}^\triangleright )^\tx{op}}$, to the pullback of $\Upsilon$ along $\Delta^{\{1\}} \times (\s{C}^\triangleright )^\tx{op} \times \{c \} \hookrightarrow \Delta^1 \times (\s{C}^\triangleright )^\tx{op} \times \s{C}^\triangleright$.  On the other hand, for each $c' \in (\s{C}^\triangleright )^\tx{op}$, the fiber of $\lambda''$ over $(c',0)$ is isomorphic (as a simplicial set) to $\Delta^0$.  Thus, by (\cite{kerodon}, 02S8), we may assume $z$ is a morphism $z : \ul{\Delta}^0 \to F(1,c)$ in $\s{P}(\s{C}^\triangleright )$, where $\ul{\Delta}^0 \in \s{P}(\s{C}^\triangleright )$ is the constant functor to $\Delta^0 \in \s{S}$. 

By (\cite{kerodon}, 02RB), there is a lift $H$ of $\lambda''' \times \tx{id}_{\Delta^1}$ that makes the following diagram commute:
\begin{equation}\label{hliftcovarianttransportsq} \begin{tikzcd}
 \tx{Tw}(\s{C}_\diamond )\times_{\s{C}_\diamond^\tx{op} \times \s{C}_\diamond} \big( (\s{C}^\triangleright )^\tx{op} \times \{ v\} \big) \times \Delta^{\{0\} } \arrow[hook]{rr} \arrow[hook]{d} & & \tx{Tw}(\s{C}_\diamond )\times_{\s{C}_\diamond^\tx{op} \times \s{C}_\diamond} \big( (\s{C}^\triangleright )^\tx{op} \times \Delta^1 \big) \arrow[]{d}{\lambda''} \\
 \tx{Tw}(\s{C}_\diamond )\times_{\s{C}_\diamond^\tx{op} \times \s{C}_\diamond} \big( (\s{C}^\triangleright )^\tx{op} \times \{ v\} \big) \times \Delta^1 \arrow[]{rr}[swap]{\lambda''' \times \tx{id}_{\Delta^1}} \arrow[dashed]{urr}{H} & & (\s{C}^\triangleright )^\tx{op} \times \Delta^1 .
\end{tikzcd}\end{equation}
For any $c' \in (\s{C}^\triangleright )^{\tx{op}}$, we can pull back the lower triangle of (\ref{hliftcovarianttransportsq}) along $\{ c'\} \times \Delta^1 \hookrightarrow (\s{C}^\triangleright )^\tx{op} \times \Delta^1$ to get
\[\begin{tikzcd}
\tx{Tw}(\s{C}_\diamond ) \times_{\s{C}_\diamond^\tx{op} \times \s{C}_\diamond } \{ (c' , v) \} \times \Delta^1 \arrow[]{rr}{H_{\{c' \} \times \Delta^1}} \arrow[]{dr}[swap]{! \times \tx{id}_{\Delta^1}} & & \tx{Tw}(\s{C}_\diamond ) \times_{\s{C}_\diamond^{\tx{op}} \times \s{C}_\diamond } (\{ c' \} \times \Delta^1 ) \arrow[]{dl}{\lambda''''} \\
 & \{ c'\} \times\Delta^1 , &
\end{tikzcd}\]
where $\lambda''''$ is the pullback of $\lambda_{\s{C}_\diamond}$ along $c' \times t : \Delta^0 \times \Delta^1 \to \s{C}_\diamond^{\tx{op}} \times \s{C}_\diamond $.  We have an isomorphism of simplicial sets $\tx{Tw}(\s{C}_\diamond ) \times_{\s{C}_\diamond^\tx{op} \times \s{C}_\diamond } \{ (c' , v) \} \cong \Delta^0$, and $H_{\{c' \} \times \Delta^1}$ is an edge of $\tx{Tw}(\s{C}_\diamond )$ given by a map $\alpha : (\Delta^1 )^\tx{op} \star \Delta^1 \to \s{C}_{\diamond }$ such that $\alpha |_{(\Delta^1 )^\tx{op}}$ is $\tx{id}_{c'}$ and $\alpha |_{\Delta^1 }$ is $t$.

\hfill

By Remark~\ref{ptunptadjunctionremark}, it suffices to show that there is a 2-simplex $\beta : \Delta^2 \to \s{P}(\s{C}^\triangleright )$ such that the edge $\beta|_{\Delta^{\{0,1 \}}}$ is $F(0,c) \to F(0,c) \coprod \ul{\Delta}^0$, the edge $\beta|_{\Delta^{\{0,2 \}}}$ is $F|_{\Delta^1 \times \{ c\}}$, and the edge $\beta|_{\Delta^{\{1,2 \}}}$ is an isomorphism.  

Consider the two $1$-simplicies of $\s{P}(\s{C}^\triangleright )_{/ F(1,c)}$ given by the unique map $\ul{\Delta}^{-1} \to z$ and $\tx{id}_{F|_{\Delta^1 \times \{c\}}}$, i.e.
\[ \begin{tikzcd}
\ul{\Delta}^{-1} \arrow[]{rr}{!} \arrow[]{dr}[swap]{!} & & \ul{\Delta}^0 \arrow[]{dl}{z} & F(0,c) \arrow[]{rr}{\tx{id}} \arrow[]{dr}[swap]{F|_{\Delta^1 \times \{c\}}} & &  F(0,c) \arrow[]{dl}{F|_{\Delta^1 \times \{c\}}} \\
 & F(1,c) & , & & F(1,c). &
\end{tikzcd}\]
This pair of $1$-simplicies corresponds to a map $z' : \partial \Delta^1 \to \tx{Fun}\big( \Delta^1 , \s{P}(\s{C}^\triangleright )_{/F(1,c)} \big)$.  By (\cite{kerodon}, 02KC and 05QG), the $\infty$-category $\tx{Fun}\big( \Delta^1 , \s{P}(\s{C}^\triangleright )_{/F(1,c)} \big)$ admits $\partial \Delta^1$-indexed colimits, and moreover we have a colimit of $z'$ in $\tx{Fun}\big( \Delta^1 , \s{P}(\s{C}^\triangleright )_{/F(1,c)} \big)$ which corresponds to a $2$-simplex of $\s{P}(\s{C}^\triangleright )$ of the form 
\[\begin{tikzcd}
F(0,c) \arrow[]{rr}  \arrow[]{dr}[swap]{F|_{\Delta^1 \times \{c\}}}  & & F(0,c) \coprod \ul{\Delta}^0 \arrow[]{dl}{(F|_{\Delta^1 \times \{c\} }, z)} \\
 & F(1,c), & 
\end{tikzcd}\]
where the (homotopy) coproduct is taken in $\s{P}(\s{C}^\triangleright )$.  It remains to show that the morphism $(F|_{\Delta^1 \times \{c\} }, z) : F(0,c) \coprod \ul{\Delta}^0 \to F(1,c)$ in $\s{P}(\s{C}^\triangleright )$ is an isomorphism, for which it suffices to show that $(F|_{\Delta^1 \times \{c\} }, z)(c')$ is an isomorphism in $\s{S}$ for each $c' \in (\s{C}^\triangleright )^\tx{op}$.  

Pick any $c' \in (\s{C}^\triangleright )^\tx{op}$. By Remark~\ref{upsilonunwindrmk}, we have commutative squares in $\s{S}$ (more precisely, a map $\Lambda^2_1 \times \Delta^1 \to \s{S}$) of the form
\[\begin{tikzcd}
\tx{Hom}_{\s{C}^\triangleright }^\tx{L} (c',c) \arrow[]{rr}{\zeta \circ -} \arrow[]{d}[swap]{(\iota_{c'})_{c}} & & \tx{Hom}_{\s{C}_\diamond }^\tx{L} (c',c) \arrow[]{d}{(\iota_{\zeta(c')})_{\zeta (c)}} \\
 \tx{Tw}(\s{C}^\triangleright )\times_{(\s{C}^\triangleright )^\tx{op} \times \s{C}^\triangleright } \{ (c',c) \} \arrow[]{rr}[swap]{\tx{Tw}(\zeta )_{(c',c)}} \arrow[]{d}[swap]{\alpha} & & \tx{Tw}(\s{C}_\diamond )\times_{\s{C}_\diamond^\tx{op} \times \s{C}_\diamond } \{ (c',c) \} \arrow[]{d}{\gamma} \\
F(0,c)(c') \arrow[]{rr}[swap]{F|_{\Delta^1 \times \{c\}} (c') } & & F(1,c)(c')
\end{tikzcd}\]
in which the morphisms $(\iota_{c'})_{c}$, $(\iota_{\zeta(c')})_{\zeta (c)}$, $\alpha$ and $\gamma$ of $\s{S}$ are isomorphisms. 
We define a vertex $\chi_{c',c}$ of $\tx{Hom}_{\s{C}_\diamond }^\tx{L} (c',c)$ as follows: 
\begin{itemize}
\item If $c' \neq v$, then let $\chi_{c',c} : \Delta^0 \to \tx{Hom}_{\s{C}_\diamond }^\tx{L} (c',c)$ be the morphism of Kan complexes
\[ (\Delta^n \to \Delta^0 ) \xmapsto{\chi_{c',c}} 
\left[ \Delta^{n+1} \hookleftarrow \Delta^{\{0 \}} \coprod \Delta^{\{1, \dots , n+1 \}} \xrightarrow{\sigma} \s{C} \right]  \]
where $\sigma|_{\Delta^{\{0\}}} = c'$ and $\sigma|_{\Delta^{\{1, \dots , n+1 \}}}$ is the constant map to $c$.  
\item If $c' = v$, then let $\chi_{c',c} : \Delta^0 \to \tx{Hom}_{\s{C}_\diamond }^\tx{L} (c',c)$ be the morphism of Kan complexes
\[ (\Delta^n \to \Delta^0 ) \xmapsto{\chi_{c',c}} 
\left[ \Delta^{n+1} \hookleftarrow  \Delta^{\{1, \dots , n+1 \}} \xrightarrow{\sigma} \s{C} \right]  \]
where $\sigma$ is the constant map to $c$.  
\end{itemize}
By the same reasoning as in Remark~\ref{homrdiamondrmk}, since $c \in \s{C}^\triangleright $ is not the cone point, the cospan of Kan complexes
\[  \tx{Hom}_{\s{C}^\triangleright }^\tx{L} (c',c) \xrightarrow{\zeta\circ -} \tx{Hom}_{\s{C}_\diamond }^\tx{L} (c',c) \xleftarrow{\chi_{c',c}} \Delta^0 \]
exhibits $\tx{Hom}_{\s{C}_\diamond }^\tx{L} (c',c)$ as a coproduct of simplicial sets $\tx{Hom}_{\s{C}^\triangleright }^\tx{L} (c',c) \coprod \Delta^0$.  Moreover, the diagram of simplicial sets
\[\begin{tikzcd}
\Delta^0\arrow[]{rrr}{\chi_{c' ,c}} \arrow[]{d}[swap]{!} & & & \tx{Hom}^\tx{L}_{\s{C}_\diamond} (c',c) \arrow[]{d}{( \iota_{\zeta (c')} )_{\zeta (c)}} \\
\tx{Tw}(\s{C}_\diamond ) \times_{\s{C}_\diamond^\tx{op} \times \s{C}_\diamond } \{ (c' , v) \} \arrow[]{rrr}[swap]{H_{\{c' \} \times \Delta^1  } |_{\{ (c' , 1) \} }} & & & \tx{Tw}(\s{C}_\diamond ) \times_{\s{C}_\diamond^\tx{op} \times \s{C}_\diamond } \{ (c' , c) \}
\end{tikzcd}\]
commutes.  Thus the diagram $(\partial \Delta^1 )^\triangleright \to \s{S}$ given by
\[ \Delta^0 \xrightarrow{z(c')} F(1,c)(c') \xleftarrow{F|_{\Delta^1 \times \{ c\}} (c')} F(0,c)(c')\]
is a colimit diagram, so $(F|_{\Delta^1 \times \{c\}} , z) : F(0,c) \coprod \ul{\Delta}^0 \to F(1,c)$ is an isomorphism in $\s{P}(\s{C}^\triangleright )$ by (\cite{kerodon}, 05QG).
\end{proof}
\begin{rmk} Using the same notation as Remark~\ref{ptunptadjunctionremark} and Proposition~\ref{zetaobjectwiselocalization}, let $c \in \s{C}^\triangleright $ be any vertex and let 
\[\lambda_c : \tx{Tw}(\s{C}_\diamond ) \times_{\s{C}_\diamond^\tx{op} \times \s{C}_\diamond} \left( (\s{C}^\triangleright )^\tx{op} \times \{ c\} \right) \to (\s{C}^\triangleright )^\tx{op} \]
be the pullback of $\lambda_{\s{C}_\diamond}$ along $(\s{C}^\triangleright )^\tx{op} \cong (\s{C}^\triangleright )^\tx{op} \times \Delta^0 \xhookrightarrow{\zeta^\tx{op} \times c} \s{C}_\diamond^\tx{op} \times \s{C}_\diamond$.  If $g_c :  (\s{C}^\triangleright )^\tx{op} \to \s{S}$ is a covariant transport representation of $\lambda_c$, then by Proposition~\ref{zetaobjectwiselocalization}, $g_c$ is isomorphic in $\s{P}(\s{C}^\triangleright )$ to $L\circ j_0 (c)$, where $j_0 : \s{C}^\triangleright \to \s{P}(\s{C}^\triangleright)$ is the Yoneda embedding and $L : \s{P}(\s{C}^\triangleright ) \to S^{-1} \s{P}(\s{C}^\triangleright )$ is any left adjoint to the inclusion $ S^{-1} \s{P}(\s{C}^\triangleright ) \hookrightarrow   \s{P}(\s{C}^\triangleright )$.  Thus for any $c , c' \in \s{C}_\diamond $, we have isomorphisms in $\s{H}$
\begin{align*} \tx{Map}_{\s{P}(\s{C}^\triangleright ) } (g_c , g_{c'}) & \simeq  \tx{Map}_{S^{-1} \s{P}(\s{C}^\triangleright ) } \left( L\circ j_0 (c), g_{c'} \right) \\
& \simeq \tx{Map}_{\s{P}(\s{C}^\triangleright ) } \left( j_0 (c), g_{c'}\right) \\
& \simeq g_{c'} (c) \simeq \tx{Map}_{\s{C}_\diamond }(c,c'), \end{align*}
where the last isomorphism follows from Remark~\ref{upsilonunwindrmk}.  In Proposition~\ref{diamondisfreezerocor}, we make this a bit more precise, to show there is an equivalence of $\infty$-categories $\s{C}_\diamond \xrightarrow{\sim} \s{P}^{\{ \Delta^{-1} \}}_{\s{R}}(\s{C}^\triangleright)$.
\end{rmk}
\begin{proposition}\label{diamondisfreezerocor}
Let $\s{C}$ be an $\infty$-category, let $\zeta = \zeta_\s{C} : \s{C}^\triangleright \hookrightarrow \s{C}_\diamond$ be the natural inclusion of Construction~\ref{triangletodiamondinclusion}, and let $\lambda_{\s{C}_\diamond} : \tx{Tw}(\s{C}_\diamond) \to \s{C}_\diamond^\tx{op} \times \s{C}_\diamond$ be the left fibration of (\cite{kerodon}, 03JK).  Let $\ul{\lambda}$ be the pullback of $\lambda_{\s{C}_\diamond}$ along $\zeta^\tx{op} \times \tx{id}_{\s{C}_\diamond}$:
\[\begin{tikzcd}
 \tx{Tw}(\s{C}_\diamond) \times_{ \s{C}_\diamond^\tx{op} \times \s{C}_\diamond} (  (\s{C}^\triangleright )^\tx{op} \times \s{C}_\diamond ) \arrow[hook]{r}{\beta} \arrow[]{d}{\ul{\lambda}} &  \tx{Tw}(\s{C}_\diamond) \arrow[]{d}{\lambda_{\s{C}_\diamond}} \\
 (\s{C}^\triangleright )^\tx{op} \times \s{C}_\diamond \arrow[hook]{r}[swap]{\zeta^\tx{op} \times \tx{id}}  & \s{C}_\diamond^\tx{op} \times \s{C}_\diamond .
\end{tikzcd}\]
Let $\s{G} :  (\s{C}^\triangleright )^\tx{op} \times \s{C}_\diamond \to \s{S}$ be a covariant transport representation of $\ul{\lambda}$, and denote by $G : \s{C}_\diamond \to \s{P}(\s{C}^\triangleright )$ the currying of $\s{G}$.  Then $G$ is fully faithful, and the essential image of $G$ is the full subcategory $\s{P}^{\{ \Delta^{-1} \}}_{\s{R}}(\s{C}^\triangleright) \subseteq \s{P}(\s{C}^\triangleright )$ from Remark~\ref{ptunptadjunctionremark}.  Thus $G$ gives an equivalence of $\infty$-categories $\s{C}_\diamond \xrightarrow{\sim} \s{P}^{\{ \Delta^{-1} \}}_{\s{R}}(\s{C}^\triangleright)$.
\end{proposition}
\begin{proof} As in Remark~\ref{ptunptadjunctionremark} and Proposition~\ref{zetaobjectwiselocalization}, let $S \subseteq \tx{Hom}(\Delta^1 , \s{P}(\s{C}^\triangleright ))$ be the singleton set $S = \{ \ul{\Delta}^{-1} \to j_0 (v) \}$, where $j_0 : \s{C}^\triangleright \to \s{P}(\s{C}^\triangleright)$ is the Yoneda embedding, $v \in \s{C}^\triangleright$ is the cone point, and $\ul{\Delta}^{-1}$ is the constant functor at $\Delta^{-1} \in \s{S}$.  Let $L : \s{P}(\s{C}^\triangleright ) \to S^{-1} \s{P}(\s{C}^\triangleright )$ be left adjoint to the inclusion $S^{-1} \s{P}(\s{C}^\triangleright ) \hookrightarrow \s{P}(\s{C}^\triangleright )$.  By Proposition~\ref{zetaobjectwiselocalization} and Remark~\ref{upsilonunwindrmk}, for each $c\in \s{C}_\diamond$ we have that $G (c)$ is isomorphic to $L( j_0 (c))$.  By Remark~\ref{termtoptrmk}, $\s{P}^{\{ \Delta^{-1} \}}_{\s{R}}(\s{C}^\triangleright)$ is the full subcategory of $\s{P}(\s{C}^\triangleright )$ spanned by all the objects of $\s{P}(\s{C}^\triangleright )$ which are isomorphic to $L( j_0 (c))$ for some $c\in \s{C}^\triangleright$, so the essential image of $G$ is  $\s{P}^{\{ \Delta^{-1} \}}_{\s{R}}(\s{C}^\triangleright)$.

\hfill

It remains to show that $G$ is fully faithful. Let $F : \Delta^1 \times \s{C}^\triangleright \to \s{P}(\s{C}^\triangleright )$ be as in the statement of Proposition~\ref{zetaobjectwiselocalization}, and fix a vertex $c\in \s{C}_\diamond$; by abuse of notation, we will consider $c$ as a vertex of $\s{C}^\triangleright$ as well.  Let $\c{H} : \s{C}_\diamond^\tx{op} \times \s{C}_\diamond \to \s{S}$ be a $\tx{Hom}$-functor for $\s{C}_\diamond$, let $H : \s{C}_\diamond \to \s{P} (\s{C}_\diamond )$ be the currying of $\c{H}$, and let $(\zeta^\tx{op})^* : \s{P}(\s{C}_\diamond ) \to \s{P}(\s{C}^\triangleright )$ be the functor of $\infty$-categories
\[ \begin{tikzcd}
\tx{Hom}_{\tx{sSet}} \big( \Delta^n , \tx{Fun} (\s{C}_\diamond^\tx{op} , \s{S}) \big) \arrow[]{rr}{(\zeta^\tx{op})^* \circ - } & & \tx{Hom}_{\tx{sSet}} \big( \Delta^n , \tx{Fun} ((\s{C}^\triangleright )^\tx{op} , \s{S}) \big) \\[-1.5em]
\left( \Delta^n \times \s{C}_\diamond^\tx{op} \xrightarrow{\sigma} \s{S} \right) \arrow[mapsto]{rr} & & \left( \Delta^n \times (\s{C}^\triangleright )^\tx{op} \xhookrightarrow{\tx{id} \times \zeta^\tx{op}} \Delta^n \times \s{C}_\diamond^\tx{op} \xrightarrow{\sigma} \s{S} \right) .
\end{tikzcd}\]
We have a commutative diagram of $\infty$-categories
\begin{equation}\label{diamondisfreezerocorpentagon1} \begin{tikzcd}
\s{C}_\diamond \arrow[equal]{r} \arrow[]{d}[swap]{H} & \s{C}_\diamond \arrow[]{d}{(\zeta^\tx{op})^* \circ H} \\
\s{P}(\s{C}_\diamond ) \arrow[]{r}{(\zeta^\tx{op})^* } \arrow[bend right]{dr}[swap]{\tx{ev}_c } & \s{P}(\s{C}^\triangleright ) \arrow[]{d}{\tx{ev}_{\zeta^\tx{op} (c)} } \\
& \s{S}
\end{tikzcd}\end{equation}
where $\tx{ev}_c$ stands for evaluation at $c$ and $\tx{ev}_{\zeta^\tx{op} (c)}$ stands for evaluation at $c = \zeta^\tx{op} (c) \in (\s{C}^\triangleright )^\tx{op}$.  In particular, we have an equality of vertices
\[ \tx{ev}_c (H(c)) = \c{H}(c,c) = \tx{ev}_{\zeta^\tx{op} (c)} \big( \big( (\zeta^\tx{op})^* \circ H \big) (c) \big) \]
in $\s{S}$.  By the Yoneda lemma, 
 the functor $H$ is fully faithful.  Both $\s{G}$ and $\c{H} \circ (\zeta^\tx{op} \times \tx{id}_{\s{C}_\diamond})$ are covariant transport representations of $\ul{\lambda}$, so $\s{G}$ and $\c{H} \circ (\zeta^\tx{op} \times \tx{id}_{\s{C}_\diamond})$ are isomorphic as objects of the $\infty$-category $\tx{Fun}\left( (\s{C}^\triangleright )^\tx{op} \times \s{C}_\diamond , \s{S} \right)$, so $G$ and $(\zeta^\tx{op})^* \circ H$ are isomorphic as objects of the $\infty$-category $\tx{Fun}(\s{C}_\diamond , \s{P}(\s{C}^\triangleright))$.  By (\cite{rezk}, 46.4), $G$ is fully faithful if and only if $(\zeta^\tx{op})^* \circ H$ is fully faithful, so it suffices to show that $(\zeta^\tx{op})^* \circ H$ is fully faithful.

\hfill

Let $f_c : \Delta^1 \to \s{P}(\s{C}^\triangleright )$ be the composite $\Delta^1 \cong \Delta^1 \times \{c\} \xrightarrow{F|_{\Delta^1 \times \{c\}}} \s{P}(\s{C}^\triangleright ) $, so we have $f_c : F(0,c) \to F(1,c)$.  By Remark~\ref{upsilonunwindrmk}, $F(0,c) : (\s{C}^\triangleright )^\tx{op} \to \s{S}$ is representable by $c \in \s{C}^\triangleright $, in the sense of (\cite{kerodon}, 038B), so we can pick a vertex $\eta'$ of the Kan complex $F(0,c)(c)$ that exhibits $F(0,c)$ as represented by $c$.  Both $\big( (\zeta^\tx{op})^* \circ H \big) (c)$ and $F(1,c)$ are covariant transport representations for the pullback of $\lambda_{\s{C}_\diamond }$ along $(\s{C}^\triangleright )^\tx{op} \cong (\s{C}^\triangleright )^\tx{op} \times \{ c\} \xrightarrow{\zeta^\tx{op} \times c} \s{C}_\diamond^\tx{op} \times \s{C}_\diamond$:
\[\begin{tikzcd}
 \tx{Tw}(\s{C}_\diamond) \times_{ \s{C}_\diamond^\tx{op} \times \s{C}_\diamond} (  (\s{C}^\triangleright )^\tx{op} \times \{ c\} ) \arrow[hook]{r} \arrow[]{d} &  \tx{Tw}(\s{C}_\diamond) \arrow[]{d}{\lambda_{\s{C}_\diamond}} \\
 (\s{C}^\triangleright )^\tx{op} \times \{ c\} \arrow[hook]{r}[swap]{\zeta^\tx{op} \times c }  & \s{C}_\diamond^\tx{op} \times \s{C}_\diamond .
\end{tikzcd}\]
Thus we can pick an isomorphism $w : F(1,c) \xrightarrow{\sim} \big( (\zeta^\tx{op})^* \circ H \big) (c)$ in $\s{P} (\s{C}^\triangleright )$ in the $\infty$-category $\s{P}(\s{C}^\triangleright )$.  Let $\nu_c : F(0,c) \to  \big( (\zeta^\tx{op})^* \circ H \big) (c)$ be a composite of $w$ and $f_c$ in the $\infty$-category $\s{P}(\s{C}^\triangleright )$, and let $\eta : \Delta^0 \to \tx{ev}_{\zeta^\tx{op} (c)} \big( \big( (\zeta^\tx{op})^* \circ H \big) (c) \big)$ be a composite of the two edges 
\begin{align*} \eta' : & \tx{ } \Delta^0  \rightarrow F(0,c)(c)=  \tx{ev}_{\zeta^\tx{op} (c)} \big( F(0,c) \big) , \\
\tx{ev}_{\zeta^\tx{op} (c)} (\nu_c ) : & \tx{ } \tx{ev}_{\zeta^\tx{op} (c)} \big( F(0,c) \big) \to \tx{ev}_{\zeta^\tx{op} (c)} \big( \big( (\zeta^\tx{op})^* \circ H \big) (c) \big) \end{align*}
in the $\infty$-category $\s{S}$.

Now pick any $c' \in \s{C}_\diamond$.  Diagram (\ref{diamondisfreezerocorpentagon1}) induces a commutative diagram of Kan complexes
\[\begin{tikzcd}
\tx{Hom}^\tx{M}_{\s{C}_\diamond } (c,c') \arrow[equal]{r} \arrow[]{d}[swap]{H} & \tx{Hom}^\tx{M}_{\s{C}_\diamond } (c,c') \arrow[]{d}{(\zeta^\tx{op})^* \circ H} \\
\tx{Hom}^\tx{M}_{\s{P}(\s{C}_\diamond )} (H(c),H(c')) \arrow[]{r}{(\zeta^\tx{op} )^*}  \arrow[bend right=10]{dr}[swap]{\tx{ev}_{c} } & \tx{Hom}^\tx{M}_{\s{P}(\s{C}^\triangleright )} (H(c)\circ \zeta^\tx{op} ,H(c')\circ \zeta^\tx{op} )  \arrow[]{d}{\tx{ev}_{\zeta^\tx{op} (c)} } \\
& \tx{Hom}^\tx{M}_{\s{S}} \big( \c{H}(c,c) , \c{H}(c,c') \big),
\end{tikzcd}\]
in which $H : \tx{Hom}^\tx{M}_{\s{C}_\diamond } (c,c') \to \tx{Hom}^\tx{M}_{\s{P}(\s{C}_\diamond )} (H(c),H(c'))$ is a weak equivalence, and we want to show that $(\zeta^\tx{op})^* \circ H :  \tx{Hom}^\tx{M}_{\s{C}_\diamond } (c,c') \to \tx{Hom}^\tx{M}_{\s{P}(\s{C}^\triangleright )} (H(c)\circ \zeta^\tx{op} ,H(c')\circ \zeta^\tx{op} ) $ is a weak equivalence.  In $\tx{hKan}$, we have a commutative diagram
\[\begin{tikzcd}
 \tx{Hom}^\tx{M}_{\s{P}(\s{C}^\triangleright )} (H(c)\circ \zeta^\tx{op} ,H(c')\circ \zeta^\tx{op} )  \arrow[]{d}[swap]{\tx{ev}_{\zeta^\tx{op} (c)} } \arrow[]{rr}{-\circ [\nu_c ]} & & \tx{Hom}^\tx{M}_{\s{P}(\s{C}^\triangleright )} (F(0,c) ,H(c')\circ \zeta^\tx{op} )  \arrow[]{d}{\tx{ev}_{\zeta^\tx{op} (c)} } \\
\tx{Hom}^\tx{M}_{\s{S}} \big( \c{H}(c,c) , \c{H}(c,c') \big) \arrow[bend right=10]{drr}[swap]{-\circ [\eta]} \arrow[]{rr}{-\circ [\tx{ev}_{\zeta^\tx{op} (c)} (\nu_c )]} & & \tx{Hom}^\tx{M}_{\s{S}} \big( F(0,c)(c) , \c{H}(c,c') \big) \arrow[]{d}{-\circ [\eta' ]} \\
& & \tx{Hom}^\tx{M}_{\s{S}} \big( \Delta^0 , \c{H}(c,c') \big).
\end{tikzcd}\]
By Proposition~\ref{zetaobjectwiselocalization}, $H(c')\circ \zeta^\tx{op} \in S^{-1} \s{P}(\s{C}^\triangleright )$ and $\nu_c$ is a localization relative to $S^{-1} \s{P}(\s{C}^\triangleright )$.  Thus
\[ -\circ [\nu_c ] : \tx{Hom}^\tx{M}_{\s{P}(\s{C}^\triangleright )} (H(c)\circ \zeta^\tx{op} ,H(c')\circ \zeta^\tx{op} )  \to \tx{Hom}^\tx{M}_{\s{P}(\s{C}^\triangleright )} (F(0,c) ,H(c')\circ \zeta^\tx{op} )  \]
is an isomorphism in $\tx{hKan}$.  On the other hand, the composite 
\[ \tx{Hom}^\tx{M}_{\s{P}(\s{C}^\triangleright )} (F(0,c) ,H(c')\circ \zeta^\tx{op} ) \xrightarrow{\tx{ev}_{\zeta^\tx{op} (c)}} 
\tx{Hom}^\tx{M}_{\s{S}} \big( F(0,c)(c) , \c{H}(c,c') \big) \xrightarrow{-\circ [\eta' ]}  \tx{Hom}^\tx{M}_{\s{S}} \big( \Delta^0 , \c{H}(c,c') \big) \]
is an isomorphism in $\tx{hKan}$ by (\cite{kerodon}, 03M5).  Thus, to show that 
\[(\zeta^\tx{op})^* \circ H :  \tx{Hom}^\tx{M}_{\s{C}_\diamond } (c,c') \to \tx{Hom}^\tx{M}_{\s{P}(\s{C}^\triangleright )} (H(c)\circ \zeta^\tx{op} ,H(c')\circ \zeta^\tx{op} )\]
is an isomorphism in $\tx{hKan}$, it suffices to show that the composite 
\[\begin{tikzcd}
\tx{Hom}^\tx{M}_{\s{P}(\s{C}_\diamond )} (H(c),H(c'))  \arrow[]{r}{\tx{ev}_{c} } &   \tx{Hom}^\tx{M}_{\s{S}} \big( \c{H}(c,c) , \c{H}(c,c') \big) \arrow[]{r}{-\circ [\eta ]} & \tx{Hom}^\tx{M}_{\s{S}} \big( \Delta^0 , \c{H}(c,c') \big)
\end{tikzcd}\]
is an isomorphism in $\tx{hKan}$.  For this, it suffices to show that $\eta : \Delta^0 \to \c{H}(c,c) = \tx{ev}_c H(c)$ exhibits the functor $H(c) : \s{C}_\diamond^\tx{op} \to \s{S}$ as represented by $c\in \s{C}_\diamond$.  

As in Proposition~\ref{zetaobjectwiselocalization}, we denote by $\s{D}$ the fiber product $\tx{Tw}(\s{C}_\diamond) \times_{ (\s{C}_\diamond )^\tx{op} \times \s{C}_\diamond} \big( (\s{C}^\triangleright )^\tx{op}\times \s{C}^\triangleright \big)$.  Let $\vartheta : \Delta^1 \to  \tx{Tw}( \s{C}^\triangleright ) \star_\s{D} \s{D}$ be the morphism of relative joins $\Delta^1 \cong \Delta^0 \star_{\Delta^0} \Delta^0 \to \tx{Tw}( \s{C}^\triangleright ) \star_\s{D} \s{D}$ induced by the morphism of cospans
\[\begin{tikzcd}
\Delta^0 \arrow[]{d}[swap]{\tx{id}_c} \arrow[]{r} & \Delta^0 \arrow[]{d}{\gamma (\tx{id}_c )} & \Delta^0 \arrow[]{l} \arrow[]{d}{\gamma (\tx{id}_c )} \\
\tx{Tw}(\s{C}^\triangleright ) \arrow[]{r}{\gamma} & \tx{Tw}(\s{C}_\diamond) \times_{ (\s{C}_\diamond )^\tx{op} \times \s{C}_\diamond} \big( (\s{C}^\triangleright )^\tx{op}\times \s{C}^\triangleright \big) & \tx{Tw}(\s{C}_\diamond) \times_{ (\s{C}_\diamond )^\tx{op} \times \s{C}_\diamond} \big( (\s{C}^\triangleright )^\tx{op}\times \s{C}^\triangleright \big) \arrow[]{l}[swap]{\tx{id}},
\end{tikzcd}\] 
where $\gamma := (  \tx{Tw}(\zeta ), \lambda_{\s{C}^\triangleright } )$, as in Construction~\ref{upsilonconstr}.  The edge $\vartheta$ is an ($\Upsilon_\zeta$-coCartesian) lift of $\Delta^1 \times \{ (c,c) \}$ starting at $\tx{id}_c \in  \tx{Tw}( \s{C}^\triangleright ) \cong \left( \tx{Tw}( \s{C}^\triangleright ) \star_\s{D} \s{D} \right) \times_{\Delta^1} \Delta^{\{0 \}} \subseteq  \tx{Tw}( \s{C}^\triangleright ) \star_\s{D} \s{D}$:  
\[\begin{tikzcd}
\Delta^{\{0\}} \arrow[hook]{d} \arrow[]{r}{\tx{id}_c} & \tx{Tw}( \s{C}^\triangleright )  \arrow[hook]{r} & \tx{Tw}( \s{C}^\triangleright ) \star_\s{D} \s{D} \arrow[]{d}{\Upsilon_\zeta} \\
\Delta^1 \arrow[]{urr}{\vartheta} \arrow[]{r}[swap]{\cong} & \Delta^1 \times \{ (c,c) \} \arrow[hook]{r} & \Delta^{1} \times (\s{C}^\triangleright )^\tx{op} \times \s{C}^\triangleright .
\end{tikzcd}\]
We have two pullback squares
\[\begin{tikzcd}
\tx{Tw}(\s{C}^\triangleright ) \arrow[hook]{r} \arrow[]{d}[swap]{\lambda_{\s{C}^\triangleright }} & \tx{Tw}( \s{C}^\triangleright ) \star_\s{D} \s{D} \arrow[]{d}{\Upsilon_\zeta} &   \s{D} \arrow[hook]{l} \arrow[]{d}{\lambda'} \\
\Delta^{\{0\}} \times (\s{C}^\triangleright )^\tx{op} \times \s{C}^\triangleright \arrow[hook]{r} & \Delta^{1} \times (\s{C}^\triangleright )^\tx{op} \times \s{C}^\triangleright & \Delta^{\{1\}} \times (\s{C}^\triangleright )^\tx{op} \times \s{C}^\triangleright \arrow[hook]{l}
\end{tikzcd}\]
where $\lambda'$ is the pullback of $\lambda_{\s{C}_\diamond }: \tx{Tw}(\s{C}_\diamond ) \to \s{C}_\diamond^\tx{op} \times \s{C}_\diamond$ along $\zeta^\tx{op} \times \zeta$, and we have a commutative diagram of simplicial sets
\[\begin{tikzcd}
\s{D}\arrow[]{d}[swap]{\lambda'} \arrow[hook]{r}{\alpha} &  \tx{Tw}(\s{C}_\diamond) \times_{ \s{C}_\diamond^\tx{op} \times \s{C}_\diamond} (  (\s{C}^\triangleright )^\tx{op} \times \s{C}_\diamond ) \arrow[hook]{r}{\beta} \arrow[]{d}{\ul{\lambda}} &  \tx{Tw}(\s{C}_\diamond) \arrow[]{d}{\lambda_{\s{C}_\diamond}} \arrow[]{r} & \s{S}_{\Delta^0 /} \arrow[]{d} \\
 (\s{C}^\triangleright )^\tx{op} \times \s{C}^\triangleright \arrow[hook]{r}{\tx{id} \times \zeta} &  (\s{C}^\triangleright )^\tx{op} \times \s{C}_\diamond \arrow[hook]{r}{ \zeta^\tx{op} \times \tx{id}} &  \s{C}_\diamond^\tx{op} \times \s{C}_\diamond \arrow[]{r}{\c{H}}& \s{S}
\end{tikzcd}\]
where the first two squares are pullback squares and the third is a categorical pullback square.  Under the weak equivalence $\c{H}(c,c) \simeq \tx{Tw}(\s{C}_\diamond ) \times_{\s{C}_\diamond^\tx{op} \times \s{C}_\diamond} \{ (c,c) \}$ of (\cite{kerodon}, 027K), the isomorphism class of $\eta$ in $\c{H}(c,c)$ corresponds to the connected component of the vertex $\beta\circ\alpha ( \vartheta (1) )$ in $\tx{Tw}(\s{C}_\diamond ) \times_{\s{C}_\diamond^\tx{op} \times \s{C}_\diamond} \{ (c,c) \}$.  

But $\beta\circ\alpha ( \vartheta (1) ) = \tx{id}_c \in \tx{Tw}(\s{C}_\diamond ) \times_{\s{C}_\diamond^\tx{op} \times \s{C}_\diamond} \{ (c,c) \} \subseteq \tx{Tw}(\s{C}_\diamond )$, and $\tx{id}_c$ is an initial object of the $\infty$-category $\tx{Tw}(\s{C}_\diamond ) \times_{\s{C}_\diamond^\tx{op} \times \s{C}_\diamond} (\s{C}_\diamond^\tx{op} \times \{ c \})$.  By (\cite{kerodon}, 02J8),  $\eta$ exhibits the functor $H(c) = \c{H}(-,c): \s{C}_\diamond^\tx{op} \to \s{S}$ as represented by $c\in \s{C}_\diamond$.  This completes the proof that $(\zeta^\tx{op})^* \circ H : \s{C}_\diamond \to \s{P}(\s{C}^\triangleright )$ is fully faithful.
\end{proof}
In this last proposition, we've shown that $\s{C}_\diamond$ is categorically equivalent to $\s{P}^{\{ \Delta^{-1} \}}_{\s{R}}(\s{C}^\triangleright)$.  The next thing to show is that for any $\infty$-category $\s{C} \in \tx{Cat}_\infty$, the functor $\zeta_\s{C} : \s{C}^\triangleright \to \s{C}_\diamond$ is a localization in $\tx{Cat}_\infty^\tx{term}$ relative to the reflective subcategory $\tx{Cat}_\infty^\tx{pt} \subseteq \tx{Cat}_\infty^\tx{term}$.  We begin with some lemmas about the essential uniqueness of localizations.
\begin{lemma}\label{righthorninneran} Let $\sigma_0 : \Delta^2 \to \Delta^1$ be the map with $\sigma_0 (0)=\sigma_0 (1) = 0$ and $\sigma_0 (2)=1$.  Then the morphisms
\[\begin{tikzcd}[column sep=small]
( \Lambda^2_2 )^\sharp \arrow[hook]{rr}{i} \arrow[]{dr}& & ( \Lambda^2_2 )^\sharp \coprod_{( \Lambda^2_2 )^\flat} (\Delta^2)^\flat \arrow[]{dl}{\sigma_0} &    ( \Lambda^2_2 )^\sharp \coprod_{( \Lambda^2_2 )^\flat} (\Delta^2)^\flat \arrow[hook]{rr}{j}  \arrow[]{dr} & & (\Delta^2)^\sharp  \arrow[]{dl}{\sigma_0}  \\
  & (\Delta^1)^\sharp & &  &  (\Delta^1)^\sharp &
\end{tikzcd}\]
in $\tx{sSet}^+_{/\Delta^1}$ are Cartesian equivalences in the sense of (\cite{htt}, 3.1.3.3), hence trivial cofibrations in the Cartesian model structure on $\tx{sSet}^+_{/\Delta^1}$ (\cite{htt}, \S 3.1.3).
\end{lemma}
\begin{proof}  By (\cite{htt}, 3.1.3.4), it suffices to show that $i$ and $j$ are marked anodyne as morphisms of marked simplicial sets.  The morphism $i$ is marked anodyne by (\cite{htt}, 2.4.1.4) and (\cite{htt}, 3.1.1.6), and $j$ is marked anodyne by (\cite{htt}, 3.1.1.7).
\end{proof}
\begin{lemma}\label{cocartnzippinglem} Let $\s{B}$ be an $\infty$-category with a coCartesian fibration $p: \s{B}\to \Delta^1$, let $X$ be a simplicial set, and let $\iota : \Lambda^2_0 \hookrightarrow \Delta^2$ be the inclusion.  Suppose a commutative solid-arrow diagram of simplicial sets
\begin{equation}\label{cocartnzippinglemsq}\begin{tikzcd}
\Lambda^2_0 \times X \arrow[hook]{d}[swap]{\iota \times \tx{id}_X} \arrow[]{rr}{s} & & \s{B} \arrow[]{d}{p} \\
\Delta^2 \times X \arrow[]{r}[swap]{\tx{pr}_1} \arrow[dashed]{urr}{\ol{s}} & \Delta^2 \arrow[]{r}[swap]{\sigma_1} & \Delta^1
\end{tikzcd}\end{equation}
such that for each vertex $x \in X$, both edges $s |_{\Delta^{\{ 0,1\}} \times \{x\}}$ and $s |_{\Delta^{\{0,2 \}} \times \{x\}}$ of $\s{B}$ are $p$-coCartesian.  Let $\tx{Fun}'(\Delta^2 \times X , \s{B}) \subseteq \tx{Fun}(\Delta^2 \times X , \s{B})$ be the full subcategory spanned by the maps $t : \Delta^2 \times X \to \s{B}$ such that $t |_{\Delta^{\{1,2 \}} \times X} : \Delta^1 \times X \to \s{B}$ corresponds to an isomorphism in $\tx{Fun}(X,\s{B})$.  Then the fiber $\tx{Fun}' (\Delta^2 \times X,\s{B})_{(\sigma_1 \circ \tx{pr}_1 , s)} $ in the pullback square
\[\begin{tikzcd}
\tx{Fun}' (\Delta^2 \times X,\s{B})_{(\sigma_1 \circ \tx{pr}_1 , s)} \arrow[]{d} \arrow[]{r}& \tx{Fun}' (\Delta^2 \times X,\s{B}) \arrow[]{d} \\
\{ (\sigma_1 \circ \tx{pr}_1 , s) \} \arrow[hook]{r}& \tx{Fun}(\Delta^2 \times X , \Delta^1) \times_{\tx{Fun}(\Lambda^2_0 \times X , \Delta^1)} \tx{Fun}(\Lambda^2_0 \times X, \s{B})
\end{tikzcd}\]
is a contractible Kan complex.

In particular there exists a dashed extension $\ol{s}$ of $s$, making (\ref{cocartnzippinglemsq}) commute, such that the edge $\ol{s}|_{\Delta^{\{1,2 \}} \times \{ x\} }$  is an isomorphism in $\s{B}$ for each vertex $x\in X$.
\end{lemma}
\begin{proof} Restriction along the inclusion $\tx{sk}_0 X \hookrightarrow X$ induces a conservative inner fibration of $\infty$-categories $\tx{Fun}(X,\s{B}) \to \tx{Fun}(\tx{sk}_0 X,\s{B})$.  Thus for $t : \Delta^2 \times X \to \s{B}$, the restriction $t |_{\Delta^{\{1,2 \}} \times X} : \Delta^1 \times X \to \s{B}$ corresponds to an isomorphism in $\tx{Fun}(X,\s{B})$ if and only if for each vertex $x\in X$, the edge $t|_{\Delta^{\{1,2 \}} \times \{x\}}$ of $\s{B}$ is an isomorphism.  If moreover $p\circ t = \sigma_1 \circ \tx{pr}_1$, then this last condition holds if and only if for each vertex $x\in X$, the edge $t|_{\Delta^{\{1,2 \}} \times \{x\}}$ of $\s{B}$ is $p$-coCartesian.  By the same token, for any simplicial set $K$, we can identify a map of simplicial sets $K \to \tx{Fun}' (\Delta^2 \times X,\s{B})_{(\sigma_1 \circ \tx{pr}_1 , s)}$ with a map of simplicial sets $f : \Delta^2 \times X \times K \to \s{B}$ such that the following hold:
\begin{enumerate}[label=(\roman*)]
\item\label{cocartnzippinglemi} The diagram
\[\begin{tikzcd}
\Lambda^2_0 \times X \times K \arrow[]{r} \arrow[hook]{d}[swap]{\iota\times \tx{id}_{X\times K}} & \Lambda^2_0 \times X \arrow[]{r}{s} & \s{B}\arrow[]{d}{p}\\
\Delta^2 \times X \times K \arrow[]{r}[swap]{\tx{pr}_1} \arrow[]{urr}[swap]{f} & \Delta^2  \arrow[]{r}[swap]{\sigma_1}& \Delta^1
\end{tikzcd}\]
commutes, and
\item\label{cocartnzippinglemii} For each vertex $(x,k) \in X\times K$ and each edge $e$ of $\Delta^2$, the edge $f|_{e\times \{(x,k) \}}$ of $\s{B}$ is $p$-coCartesian.
\end{enumerate}

By (\cite{htt}, 3.1.4.3) and Lemma~\ref{righthorninneran}, the morphism 
\[\begin{tikzcd} (\Lambda^2_0)^\sharp \times X^\flat \arrow[]{rr}{\iota\times \tx{id}_X} \arrow[]{dr} & & (\Delta^2)^\sharp \times X^\flat \arrow[]{dl}{\sigma_1 \circ \tx{pr}_1} \\
& (\Delta^1)^\sharp &
\end{tikzcd}\]
in $\tx{sSet}^+_{/ \Delta^1}$ is a trivial cofibration in the coCartesian model structure on $\tx{sSet}^+_{/ \Delta^1}$.  Let $\s{B}^\natural$ be the marked simplicial set $(\s{B}, \s{E})$, where $\s{E}$ comprises the $p$-coCartesian edges of $\s{B}$. Restriction along $\iota \times \tx{id}_X$ induces a trivial fibration of $\infty$-categories
\[ q : \tx{Map}^\flat_{\Delta^1} \big( (\Delta^2)^\sharp \times X^\flat , \s{B}^\natural \big) \to \tx{Map}^\flat_{\Delta^1} \big( (\Lambda^2_0)^\sharp \times X^\flat , \s{B}^\natural \big). \]
The morphism $s$ is the underlying map of simplicial sets of a morphism $s' : (\Lambda^2_0)^\sharp \times X^\flat \to \s{B}^\natural$ in $\tx{sSet}^+_{/ \Delta^1}$.  Consider the fiber of $q$ over $s'$,
\[ \tx{Map}^\flat_{\Delta^1} \big( (\Delta^2)^\sharp \times X^\flat , \s{B}^\natural \big)_{s'} = \tx{Map}^\flat_{\Delta^1} \big( (\Delta^2)^\sharp \times X^\flat , \s{B}^\natural \big) \times_{ \tx{Map}^\flat_{\Delta^1} \left( (\Lambda^2_0)^\sharp \times X^\flat , \s{B}^\natural \right)} \{s' \}, \]
which is a contractible Kan complex.  

For any simplicial set $K$, a morphism $K \to  \tx{Map}^\flat_{\Delta^1} \big( (\Delta^2)^\sharp \times X^\flat , \s{B}^\natural \big)_{s'} $ in $\tx{sSet}$ can be identified with a map $f:  (\Delta^2)^\sharp \times X^\flat \times K^\flat \to \s{B}^\natural$ in $\tx{sSet}^+$ such that the following hold:
\begin{enumerate}[label=(\alph*)]
\item\label{cocartnzippinglema} The triangle 
\[ \begin{tikzcd}
(\Delta^2)^\sharp \times X^\flat \times (\tx{sk}_0 K)^\flat \arrow[]{dr}[swap]{\sigma_1 \circ \tx{pr}_1} \arrow[]{rr}{f|_{ (\Delta^2)^\sharp \times X^\flat \times (\tx{sk}_0 K)^\flat }} & & \s{B}^\natural \arrow[]{dl}{p} \\ 
& (\Delta^1)^\sharp &
\end{tikzcd}\]
in $\tx{sSet}^+$ commutes, and
\item\label{cocartnzippinglemb} The square
\[\begin{tikzcd}
(\Lambda^2_0)^\sharp \times X^\flat \times K^\flat \arrow[]{d} \arrow[]{rr}{ \iota \times \tx{id}}& & (\Delta^2)^\sharp \times X^\flat \times K^\flat \arrow[]{d}{f}\\
 (\Lambda^2_0)^\sharp \times X^\flat  \arrow[]{rr}{s'}  & & \s{B}^\natural
\end{tikzcd}\]
in $\tx{sSet}^+$ commutes.
\end{enumerate}
If $W$ is a simplicial set with a sub-simplicial set $\tx{sk}_0 W \subseteq W' \subseteq W$ and $P$ is a poset, then any map of simplicial sets $W' \to \tx{N}(P)$ has at most one extension $W \to \tx{N}(P)$.  Thus condition~\ref{cocartnzippinglema} is equivalent to the condition
\begin{enumerate}[label=(\alph*')]
\item\label{cocartnzippinglemaprime} The triangle 
\[ \begin{tikzcd}
(\Delta^2)^\sharp \times X^\flat \times K^\flat \arrow[]{dr}[swap]{\sigma_1 \circ \tx{pr}_1} \arrow[]{rr}{f} & & \s{B}^\natural \arrow[]{dl}{p} \\ 
& (\Delta^1)^\sharp &
\end{tikzcd}\]
in $\tx{sSet}^+$ commutes.
\end{enumerate}
Comparing conditions~\ref{cocartnzippinglemaprime} and~\ref{cocartnzippinglemb} with conditions~\ref{cocartnzippinglemi} and~\ref{cocartnzippinglemii}, we have a natural isomorphism
\[\tx{Fun}' (\Delta^2 \times X,\s{B})_{(\sigma_1 \circ \tx{pr}_1 , s)} \cong \tx{Map}^\flat_{\Delta^1} \big( (\Delta^2)^\sharp \times X^\flat , \s{B}^\natural \big)_{s'}, \]
so $\tx{Fun}' (\Delta^2 \times X,\s{B})_{(\sigma_1 \circ \tx{pr}_1 , s)} $ is a contractible Kan complex.
\end{proof}
\begin{proposition}\label{localizationzippingprop} Let $\s{C}$ be an $\infty$-category, $\nu : \s{C}^0 \hookrightarrow \s{C}$ a full, reflective subcategory, $X$ a simplicial set, and suppose a map $\theta : \Lambda^2_0 \to \tx{Fun}(X,\s{C})$ such that for each vertex $x\in X$, both edges $\theta (0)(x) \to \theta (1)(x)$ and $\theta (0)(x) \to \theta (2)(x)$ in $\s{C}$ are localizations relative to $\s{C}^0$.  Then there exists an extension $\ol{\theta} : \Delta^2 \to \tx{Fun}(X,\s{C})$ of $\theta$ such that $\ol{\theta}|_{\Delta^{\{1,2 \}}}$ is an isomorphism.
\end{proposition}
\begin{proof} Let $\s{B}\subseteq \s{C} \times \Delta^1$ be the full subcategory spanned by the vertices $(c,i)$ such that $c\in \s{C}^0$ if $i=1$, and let $p= \tx{pr}_2 |_{\s{B}} : \s{B} \to \Delta^1$.  Then $p$ is a coCartesian fibration, and an edge $\ell : (c,0) \to (c',1)$ of $\s{B}$ is $p$-coCartesian if and only if the edge $\tx{pr}_1 (\ell ) : c\to c'$ of $\s{C}$ is a localization relative to $\s{C}^0$.

Let $t:\Lambda^2_0 \times X \to \s{C}$ be the uncurrying of $\theta$, and let $\iota : \Lambda^2_0 \hookrightarrow \Delta^2$ be the inclusion.  There is a unique dashed lift $s$ of $t$ making the following diagram commute:
\[\begin{tikzcd}
\Lambda^2_0  \arrow[hook]{r}{\iota} & \Delta^2 \arrow[]{r}{\sigma_1} & \Delta^1  \\
\Lambda^2_0 \times X \arrow[]{drr}[swap]{t} \arrow[dashed]{rr}{s} \arrow[]{u}{\tx{pr}_1}& & \s{B} \arrow[]{d}{\tx{pr}_1} \arrow[]{u}[swap]{p} \\
 & & \s{C} 
\end{tikzcd}\]

For each $x\in X$, the edges $\tx{pr}_1 \circ s|_{\Delta^{\{0,1 \}} \times \{x \}}$ and $\tx{pr}_1 \circ s|_{\Delta^{\{0,2 \}} \times \{x \}}$ of $\s{C}$ are localizations relative to $\s{C}^0$, so $  s|_{\Delta^{\{0,1 \}} \times \{x \}} $ and $ s|_{\Delta^{\{0,2 \}} \times \{x \}}$ are $p$-coCartesian edges of $\s{B}$. By Lemma~\ref{cocartnzippinglem}, there is an extension $\ol{s} : \Delta^2 \times X \to \s{B}$ of $s$ such that for each $x\in X$, the edge $\ol{s} |_{\Delta^{\{1,2 \}}  \times \{x \}}$ of $\s{B} \times_{\Delta^1} \{ 1 \}\cong \s{C}^0 \subseteq \s{B}$ is an isomorphism.  Currying the composition $\Delta^2 \times X \xrightarrow{\ol{s}} \s{B} \hookrightarrow \s{C} \times \Delta^1$ gives a map $\ol{\theta} : \Delta^2 \to \tx{Fun}(X,\s{C})$, and $\ol{\theta} |_{\Delta^{\{1,2 \}}}$ is an isomorphism in $\tx{Fun}(X, \s{C})$ by (\cite{rezk}, 37.2).
\end{proof}

\begin{proposition}\label{zetaisptlocalization} For any $\infty$-category $\s{C}\in \tx{Cat}_\infty$, the edge $\zeta = \zeta_{\s{C}} : \s{C}^\triangleright \to \s{C}_\diamond$ of $\tx{Cat}_\infty^{\tx{term}}$ is a localization relative to the full subcategory $\tx{Cat}_\infty^{\tx{pt}} \subseteq \tx{Cat}_\infty^{\tx{term}}$ (which is a reflective subcategory by Proposition~\ref{termptadj}).
\end{proposition}
\begin{proof} As in Remark~\ref{termtoptrmk}, let $j_0 : \s{C}^\triangleright \to \s{P}(\s{C}^\triangleright)$ be the Yoneda embedding, and let $S \subseteq$ $\tx{Hom}_{\tx{sSet}} (\Delta^1 , \s{P}(\s{C}^\triangleright))$ be the singleton set containing the unique edge $\ul{\Delta}^{-1} \to j_0 (v)$ from the constant presheaf at $\Delta^{-1}$ to the Yoneda embedding of the cone point $v \in \s{C}^\triangleright$.  Let $L : \s{P}(\s{C}^\triangleright ) \to S^{-1} \s{P}(\s{C}^\triangleright )$ be a left adjoint to the inclusion $S^{-1} \s{P}(\s{C}^\triangleright ) \subseteq \s{P}(\s{C}^\triangleright )$ of the full subcategory of $S$-local objects, let $\s{P}^{\{\Delta^{-1} \}}_\s{R} (\s{C}^\triangleright ) \subseteq \s{P}(\s{C}^\triangleright )$ be the essential image of the composite $L\circ j_0 : \s{C}^\triangleright \to S^{-1} \s{P}(\s{C}^\triangleright )$, and denote by $j$ the composite $L\circ j_0 : \s{C}^\triangleright \to \s{P}^{\{\Delta^{-1} \}}_\s{R} (\s{C}^\triangleright )$.  By the proof of Proposition~\ref{termptadj}, it suffices to find a $2$-simplex $\Delta^2 \to \tx{Cat}_\infty^\tx{term}$ of the form
\[\begin{tikzcd} \s{C}^\triangleright \arrow[]{r}{\zeta} \arrow[]{dr}[swap]{j}& \s{C}^\diamond \arrow[]{d}{g} \\
& \s{P}^{\{\Delta^{-1} \}}_\s{R} (\s{C}^\triangleright ), \end{tikzcd}\]
where $g : \s{C}^\diamond \to \s{P}^{\{\Delta^{-1} \}}_\s{R} (\s{C}^\triangleright )$ is a categorical equivalence.  Note that any categorical equivalence $\s{C}^\diamond \to \s{P}^{\{\Delta^{-1} \}}_\s{R} (\s{C}^\triangleright )$ would automatically be an edge of $\tx{Cat}_\infty^\tx{pt} \subseteq \tx{Cat}_\infty^\tx{term}$.  

Since $\s{P}^{\{\Delta^{-1} \}}_\s{R} (\s{C}^\triangleright ) \subseteq \s{P}(\s{C}^\triangleright)$ is a full subcategory, these data are equivalent to 
\begin{enumerate}[label=(\alph*)]
\item\label{zetaisptlocalizationa} a categorical equivalence $g : \s{C}^\diamond \to \s{P}^{\{\Delta^{-1} \}}_\s{R} (\s{C}^\triangleright )$, together with 
\item\label{zetaisptlocalizationb} an edge $h : \Delta^1 \to \tx{Fun}(\s{C}^\triangleright , \s{P}(\s{C}^\triangleright))^\simeq$ such that the two endpoints of $h$ are given by the composite functors
\begin{align*}
h(0) & = \left( \s{C}^\triangleright \xrightarrow{j} \s{P}^{\{\Delta^{-1} \}}_\s{R} (\s{C}^\triangleright ) \hookrightarrow \s{P}(\s{C}^\triangleright) \right) \in \tx{Fun}(\s{C}^\triangleright , \s{P}(\s{C}^\triangleright))^\simeq , \tx{ \; and}\\
h(1) & = \left( \s{C}^\triangleright \xrightarrow{\zeta} \s{C}_\diamond \xrightarrow{g} \s{P}^{\{\Delta^{-1} \}}_\s{R} (\s{C}^\triangleright ) \hookrightarrow \s{P}(\s{C}^\triangleright) \right) \in \tx{Fun}(\s{C}^\triangleright , \s{P}(\s{C}^\triangleright))^\simeq .
\end{align*}
\end{enumerate}
We construct $g$ and $h$ as follows: 

\hfill

As in Proposition~\ref{zetaobjectwiselocalization}, let $\s{D} := \tx{Tw}(\s{C}_\diamond) \times_{ (\s{C}_\diamond )^\tx{op} \times \s{C}_\diamond} \big( (\s{C}^\triangleright )^\tx{op}\times \s{C}^\triangleright \big)$, let \[\s{F}: \Delta^1 \times (\s{C}^\triangleright )^\tx{op} \times \s{C}^\triangleright \to \s{S}\] be a covariant transport representation of the left fibration 
\[\Upsilon = \Upsilon_{\zeta} : \tx{Tw}(\s{C}^\triangleright )\star_{\s{D}} \s{D} \to \Delta^1 \times (\s{C}^\triangleright )^\tx{op} \times \s{C}^\triangleright , \]
and let $F : \Delta^1 \times \s{C}^\triangleright \to \s{P}(\s{C}^\triangleright )$ and $f : \Delta^1 \to \tx{Fun}(  \s{C}^\triangleright , \s{P}(\s{C}^\triangleright ) ) $ be curryings of $\s{F}$. We have three pullback squares of simplicial sets:
\[\begin{tikzcd} \tx{Tw}(\s{C}^\triangleright )\star_{\s{D}} \s{D} \arrow[]{d}[swap]{\Upsilon} & \s{D} \arrow[]{d}{\lambda'} \arrow[hook]{l} \\
\Delta^1 \times (\s{C}^\triangleright )^\tx{op} \times \s{C}^\triangleright & \Delta^{\{1\}} \times (\s{C}^\triangleright )^\tx{op} \times \s{C}^\triangleright  \arrow[hook]{l} \end{tikzcd}\]
and
\[\begin{tikzcd} 
\s{D}  \arrow[]{d}[swap]{\lambda'}  \arrow[hook]{r}  & \tx{Tw}(\s{C}_\diamond )\times_{ (\s{C}_\diamond )^\tx{op} \times \s{C}_\diamond} \big( (\s{C}^\triangleright )^\tx{op}\times \s{C}_\diamond \big) \arrow[]{d}{\ul{\lambda}}  \arrow[hook]{r} & \tx{Tw}(\s{C}_\diamond ) \arrow[]{d}{\lambda_{\s{C}_\diamond}} \\
 (\s{C}^\triangleright )^\tx{op} \times \s{C}^\triangleright  \arrow[hook]{r}[swap]{\tx{id} \times \zeta}  &  (\s{C}^\triangleright )^\tx{op} \times \s{C}_\diamond   \arrow[hook]{r}[swap]{\zeta^\tx{op} \times \tx{id}} &  \s{C}_\diamond^\tx{op} \times \s{C}_\diamond .  \end{tikzcd}\]
As in Proposition~\ref{diamondisfreezerocor}, let $\s{G} : (\s{C}^\triangleright )^\tx{op} \times \s{C}_\diamond \to \s{S}$ be a covariant transport representation of $\ul{\lambda}$, and let $G : \s{C}_\diamond \to \s{P}(\s{C}^\triangleright )$ be the currying of $\s{G}$.  By Proposition~\ref{diamondisfreezerocor}, $G$ factors as $\s{C}_\diamond \xrightarrow{g} \s{P}^{\{\Delta^{-1} \}}_\s{R} (\s{C}^\triangleright ) \hookrightarrow \s{P}(\s{C}^\triangleright)$, where $g$ is a categorical equivalence (satisfying~\ref{zetaisptlocalizationa}).  Both $\s{F}|_{\Delta^{\{1\}} \times (\s{C}^\triangleright )^\tx{op} \times \s{C}^\triangleright}$ and $\s{G} \circ (\tx{id}_{(\s{C}^\triangleright )^\tx{op} } \times \zeta_{\s{C}})$ are covariant transport representations of $\lambda'$, so there is an edge $\mu : \Delta^1 \to \tx{Fun}( (\s{C}^\triangleright )^\tx{op} \times \s{C}^\triangleright , \s{S}  )^\simeq$ from $\s{F}|_{\Delta^{\{1\}} \times (\s{C}^\triangleright )^\tx{op} \times \s{C}^\triangleright}$ to $\s{G} \circ (\tx{id}_{(\s{C}^\triangleright )^\tx{op} } \times \zeta_{\s{C}})$.  Composing $\mu$ with the isomorphism of simplicial sets $\tx{Fun}( (\s{C}^\triangleright )^\tx{op} \times \s{C}^\triangleright , \s{S}  )^\simeq \xrightarrow{\cong} \tx{Fun}(  \s{C}^\triangleright , \s{P}(\s{C}^\triangleright )  )^\simeq$ gives an edge $u : \Delta^1 \to \tx{Fun}(  \s{C}^\triangleright , \s{P}(\s{C}^\triangleright )  )^\simeq$ from $F(1,-) = f(1)$ to $G \circ \zeta$.  Thus $u(1)$ is a candidate for $h(1)$ in~\ref{zetaisptlocalizationb}.

\hfill

Consider $j=L\circ j_0$ as a functor $\s{C}^\triangleright \to \s{P}(\s{C}^\triangleright )$. By Proposition~\ref{zetaobjectwiselocalization} and Proposition~\ref{localizationzippingprop}, there is a $2$-simplex $\ol{\theta} : \Delta^2 \to \tx{Fun}(\s{C}^\triangleright , \s{P}(\s{C}^\triangleright ) )$ of the form
\[ \begin{tikzcd} f(0) \arrow[]{dr}[swap]{\theta''} \arrow[]{r}{f} & f(1) \\
 & j \arrow[]{u}[swap]{w} \end{tikzcd}\]
 in which $w$ is an isomorphism in the $\infty$-category $\tx{Fun}(\s{C}^\triangleright , \s{P}(\s{C}^\triangleright ) )$.  Indeed, let $\eta$ be a unit for the adjunction $(L \dashv U) : \s{P}(\s{C}^\triangleright ) \leftrightarrows S^{-1} \s{P}(\s{C}^\triangleright )$, and let $\theta'$ be the composite
 \[\begin{tikzcd} \Delta^1 \arrow[]{r}{\eta} \arrow[bend right=15]{rr}[swap]{\theta'} &  \tx{Fun}(\s{P}(\s{C}^\triangleright ),  \s{P}(\s{C}^\triangleright )) \arrow[]{r}{-\circ j_0} & \tx{Fun}( \s{C}^\triangleright ,  \s{P}(\s{C}^\triangleright )) ,
 \end{tikzcd}\] 
 so $\theta'$ is an edge from $j_0$ to $j$ such that for each $c\in \s{C}^\triangleright$, the edge $\theta' (c)$ of $\s{P}(\s{C}^\triangleright )$ is a localization relative to $S^{-1} \s{P}(\s{C}^\triangleright )$.  By Remark~\ref{upsilonunwindrmk}, there is an edge $v : \Delta^1 \to \tx{Fun}(\s{C}^\triangleright , \s{P}(\s{C}^\triangleright ) )^\simeq$ from $f(0)$ to $j_0$.  If $\theta''$ is any composite of $v$ and $\theta'$ in the $\infty$-category $ \tx{Fun}(\s{C}^\triangleright , \s{P}(\s{C}^\triangleright ) )$, then for each $c\in \s{C}^\triangleright$, the edge $\theta'' (c)$ of $\s{P}(\s{C}^\triangleright )$ is likewise a localization relative to $S^{-1} \s{P}(\s{C}^\triangleright )$, and the edge $f$ has the same property by Proposition~\ref{zetaobjectwiselocalization}.  
 
 Let $\theta : \Lambda^2_0 \to \tx{Fun}(\s{C}^\triangleright , \s{P}(\s{C}^\triangleright ) )$ be the left horn given by $\theta |_{\Delta^{\{0,1 \}}} = \theta''$ and $\theta |_{\Delta^{\{0,2 \}}} = f$.  By Proposition~\ref{localizationzippingprop}, there is an extension $\ol{\theta} : \Delta^2 \to \tx{Fun}(\s{C}^\triangleright , \s{P}(\s{C}^\triangleright ) )$ of $\theta$ such that $w:= \ol{\theta}|_{\Delta^{\{1,2\}}}$ is an isomorphism in the $\infty$-category $\tx{Fun}(\s{C}^\triangleright , \s{P}(\s{C}^\triangleright ) )$.

\hfill

If we let $h : \Delta^1 \to \tx{Fun}(\s{C}^\triangleright , \s{P}(\s{C}^\triangleright ) )^\simeq$ be any composite of $w$ and $u$ in the $\infty$-category $\tx{Fun}(\s{C}^\triangleright , \s{P}(\s{C}^\triangleright ) )^\simeq$, then $h$ satisfies~\ref{zetaisptlocalizationb}.
\end{proof}
\begin{theorem}\label{diamondisfreezerothm} The functor of $\infty$-categories $(-)_\diamond : \tx{Cat}_\infty \to \tx{Cat}_\infty^\tx{pt}$ of Lemma~\ref{diamondinftyfunctor} is left adjoint to the inclusion $U : \tx{Cat}_\infty^\tx{pt} \hookrightarrow  \tx{Cat}_\infty$, and the natural transformation \[\iota : \Delta^1 \to \tx{Fun}(\tx{Cat}_\infty , \tx{Cat}_\infty)\] from $\tx{id}_{\tx{Cat}_\infty}$ to $U\circ (-)_\diamond$ of Lemma~\ref{diamondinftynaturaltransf} is a unit of the adjunction $((-)_\diamond \dashv U)$.
\end{theorem}
\begin{proof} Suppose $\s{C}\in \tx{Cat}_\infty$ and $\s{D} \in  \tx{Cat}_\infty^{\tx{pt}}$.  By (\cite{htt}, 5.2.2.8), it suffices to show that the composition
\begin{equation}\label{diamondisfreezerothmcomp} \tx{Map}_{ \tx{Cat}_\infty^{\tx{pt}}} (\s{C}_\diamond , \s{D}) \to \tx{Map}_{ \tx{Cat}_\infty} (\s{C}_\diamond , \s{D}) \xrightarrow{\iota (\s{C})}  \tx{Map}_{ \tx{Cat}_\infty} (\s{C}, \s{D})\end{equation}
is an isomorphism in $\s{H}$.  Let $U' : \tx{Cat}_\infty^\tx{pt} \hookrightarrow  \tx{Cat}_\infty^\tx{term}$ and $U'' : \tx{Cat}_\infty^\tx{term} \hookrightarrow  \tx{Cat}_\infty$ be the inclusions, and let $\iota' : \Delta^1 \to \tx{Fun}(\tx{Cat}_\infty , \tx{Cat}_\infty )$ from $\tx{id}_{\tx{Cat}_\infty}$ to $U'' \circ (-)^\triangleright$ be the natural transformation of Lemma~\ref{coneptsimplnatr}.  Recall that $U'$ is a full inclusion, so we can identify $\tx{Map}_{\tx{Cat}_\infty^\tx{pt}} (\s{C}_\diamond , \s{D})$ with $\tx{Map}_{\tx{Cat}_\infty^\tx{term}} (\s{C}_\diamond , \s{D})$.

Let $\f{Z} : \Delta^2 \to \tx{Fun}(\tx{Cat}_\infty , \tx{Cat}_\infty)$ be as in Lemma~\ref{frakztriangle}, so $\f{Z} (\s{C})$ is a $2$-simplex of $\tx{Cat}_\infty$ of the form 
\[ \begin{tikzcd} &\s{C}^\triangleright \arrow[]{dr}{\zeta (\s{C})}& \\
\s{C}\arrow[]{ur}{\iota' (\s{C})} \arrow[]{rr}{\iota (\s{C})}&&\s{C}_\diamond .\end{tikzcd} \]
By Proposition~\ref{catptadj}, the composition 
\[ \tx{Map}_{\tx{Cat}_\infty^{\tx{term}}} (\s{C}^\triangleright , \s{D}) \to \tx{Map}_{\tx{Cat}_\infty} (\s{C}^\triangleright , \s{D}) \xrightarrow{\iota' (\s{C})} \tx{Map}_{\tx{Cat}_\infty} (\s{C} , \s{D}) \]
is an isomorphism in $\s{H}$.  By Proposition~\ref{zetaisptlocalization}, the composition
\[ \tx{Map}_{\tx{Cat}_\infty^\tx{pt}} (\s{C}_\diamond , \s{D}) \xrightarrow{=} \tx{Map}_{\tx{Cat}_\infty^\tx{term}} (\s{C}_\diamond , \s{D}) \xrightarrow{\zeta (\s{C})}  \tx{Map}_{\tx{Cat}_\infty^\tx{term}} (\s{C}^\triangleright , \s{D}) \]
is an isomorphism in $\s{H}$.  The $2$-simplex $\f{Z}(\s{C})$ gives a commutative diagram in $\s{H}$
\[\begin{tikzcd}
\tx{Map}_{ \tx{Cat}_\infty^{\tx{pt}}} (\s{C}_\diamond , \s{D}) \arrow[equal]{r} & \tx{Map}_{ \tx{Cat}_\infty^{\tx{term}}} (\s{C}_\diamond , \s{D}) \arrow[]{r}{\zeta (\s{C})} \arrow[]{d} & \tx{Map}_{ \tx{Cat}_\infty^{\tx{term}}} (\s{C}^\triangleright , \s{D})  \arrow[]{d} \\
& \tx{Map}_{ \tx{Cat}_\infty} (\s{C}_\diamond , \s{D})\arrow[]{r}{\zeta (\s{C})} \arrow[bend right=10]{dr}[swap]{\iota (\s{C})} & \tx{Map}_{ \tx{Cat}_\infty} (\s{C}^\triangleright , \s{D})  \arrow[]{d}{\iota' (\s{C})} \\
& &  \tx{Map}_{ \tx{Cat}_\infty} (\s{C} , \s{D})
\end{tikzcd}\]
in which the top horizontal and right vertical compositions are isomorphisms, so the composition (\ref{diamondisfreezerothmcomp}) is an isomorphism.
\end{proof}



\section{Coherent choices of nullhomotopies}
In this section we give an analogue (Proposition~\ref{sectiondiamondzeroobj}) of the following lemma in the case of pointed $\infty$-categories:
\begin{lemma}[\tx{}\cite{kerodon}, 02J5]\label{ptinclsectionterminal} Let $\s{C}$ be an $\infty$-category with an object $x\in \s{C}$.  The following conditions are equivalent:
\begin{enumerate}
\item\label{ptinclsectionterminal1} The object $x\in \s{C}$ is terminal.
\item\label{ptinclsectionterminal2} The forgetful map $\s{C}_{/x} \to \s{C}$ is a trivial Kan fibration.
\item\label{ptinclsectionterminal3} There exists a functor $f : \s{C}^\triangleright \to \s{C}$ such that $f|_{\s{C}} = \tx{id}_{\s{C}}$ and such that the composition
\[ \Delta^1 \cong \{ x\}^\triangleright \hookrightarrow  \s{C}^\triangleright \xrightarrow{f} \s{C} \]
is the identity morphism $\tx{id}_x$ (i.e., the constant map $\Delta^1 \to \Delta^0 \xrightarrow{x} \s{C}$).
\item\label{ptinclsectionterminal4} The inclusion $\Delta^0 \xhookrightarrow{x} \s{C}$ is right anodyne.
\end{enumerate}
\end{lemma}
In order to prove an analogue of~\ref{ptinclsectionterminal1}$\Rightarrow$\ref{ptinclsectionterminal3} in Lemma~\ref{ptinclsectionterminal} for the $(-)_\diamond$ construction, we can't rely on an analogue of the slice construction; see (\cite{thesis}, 1.5.4).  Instead, if $\Delta^0 \xrightarrow{x} \s{C}$ is the inclusion of a zero object into a pointed $\infty$-category, we will apply a triangle identity of the adjunction of Theorem~\ref{diamondisfreezerothm} to the edge $\Delta^0 \xrightarrow{x} \s{C}$ of $\tx{Cat}_\infty^\tx{pt}$ to get a diagram $\Delta^2 \times \Delta^1 \to \tx{Cat}_\infty$ of the form 
\[\begin{tikzcd}
 \Delta^0 \arrow[]{rr}{x} \arrow[]{dr} \arrow[hook]{dd} & & \s{C} \arrow[]{dr}{\tx{id}_{\s{C}}} \arrow[hook]{dd} & \\
 & \Delta^0 \arrow[crossing over]{rr}[near start]{x} & & \s{C} \\
\Delta^0_\diamond  \arrow[]{ur} \arrow[]{rr}[swap]{x_\diamond} & & \s{C}_\diamond \arrow[]{ur} &
\end{tikzcd}\]
We want a strictly commutative diagram of simplicial sets, for which we'll need some rectification.
\begin{lemma}[\tx{}\cite{thesis}, 1.5.2]\label{cofinprojmodelstr} Let $\bf{A}$ be a combinatorial model category, and consider the poset $[1]$ as a category.  If $\bf{A}^{[1]}$ is given the projective model structure, then a morphism $j : a\to b$ in $\bf{A}^{[1]}$ is a projective cofibration if $j(0) : a(0)\to b(0)$ and $(b, j(1)) : b(0)\coprod_{a(0)} a(1) \to b(1)$ are both cofibrations in $\bf{A}$.
\end{lemma}

\begin{proposition}\label{sectiondiamondzeroobj} Let $\s{C}$ be an $\infty$-category with an object $x\in \s{C}$.  Let $\{ x \}_\diamond \subseteq  \s{C}_\diamond$ be the subcategory given by simplicies of the form $\big[ \Delta^n \hookleftarrow \tx{N}(P) \xrightarrow{\sigma} \s{C} \big]$ wherein $\sigma$ factors as $\tx{N}(P) \xrightarrow{!} \Delta^0 \xrightarrow{x} \s{C}$.  The following conditions are equivalent:
\begin{enumerate}
\item\label{sectiondiamondzeroobj1} The object $x\in \s{C}$ is a zero object.
\item\label{sectiondiamondzeroobj2} There exists a functor $f : \s{C}_\diamond \to \s{C}$ such that $f|_{\s{C}} = \tx{id}_{\s{C}}$ and $f|_{\{ x \}_\diamond} : \{ x \}_\diamond \to \s{C}$ factors as $\{ x \}_\diamond \xrightarrow{!} \Delta^0 \xrightarrow{x} \s{C}$.
\end{enumerate}
\end{proposition}
\begin{proof} First suppose there exists a functor $f : \s{C}_\diamond \to \s{C}$ as in~\ref{sectiondiamondzeroobj2}. Letting $f'$ be the composite $f':= f \circ \zeta_{\s{C}} : \s{C}^\triangleright \to \s{C}$, we see that $f'$ satisfies condition~\ref{ptinclsectionterminal3} of Lemma~\ref{ptinclsectionterminal}, so $x \in \s{C}$ is terminal.  Likewise, if $f'' := f \circ \atez_{\s{C}} : \s{C}^\triangleleft \to \s{C}$, then $(f'')^\tx{op}$ satisfies condition~\ref{ptinclsectionterminal3} of Lemma~\ref{ptinclsectionterminal}, so $x \in \s{C}$ is initial.

\hfill

Conversely, suppose $x\in \s{C}$ is a zero object, so $\Delta^0 \xrightarrow{x} \s{C}$ is an edge of $\tx{Cat}_\infty^\tx{pt}$.  Recall from Theorem~\ref{diamondisfreezerothm} the adjunction of $\infty$-categories $\left( (-)_\diamond \dashv U \right) : \tx{Cat}_\infty \rightleftarrows \tx{Cat}_\infty^\tx{pt}$ with unit $\iota : \Delta^1 \to \tx{Fun}(\tx{Cat}_\infty, \tx{Cat}_\infty)$, which (by \cite{kerodon} 02FX and \cite{htt} 5.2.2.8) is equivalently an adjunction with a given unit in the sense of (\cite{kerodon}, 02EL).  

By (\cite{kerodon}, 02F4), there exist a counit $\epsilon : \Delta^1 \to \tx{Fun}(\tx{Cat}_\infty^\tx{pt} , \tx{Cat}_\infty^\tx{pt})$ from $(-)_\diamond \circ U$ to $\tx{id}_{\tx{Cat}_\infty^\tx{pt}}$ and a $2$-simplex $\tau : \Delta^2 \to \tx{Fun}(\tx{Cat}_\infty^{\tx{pt}}, \tx{Cat}_\infty)$ such that the following hold:
\begin{itemize}
\item The restriction $\tau |_{\Delta^{\{0,2 \}}}$ is the edge $\tx{id}_U$.
\item The restriction $\tau |_{\Delta^{\{0,1 \}}}$ is the composite
\[ \Delta^1 \xrightarrow{\iota} \tx{Fun}(\tx{Cat}_\infty , \tx{Cat}_\infty) \xrightarrow{-\circ U} \tx{Fun}(\tx{Cat}_\infty^\tx{pt} , \tx{Cat}_\infty). \]
\item The restriction $\tau |_{\Delta^{\{1,2 \}}}$ is the composite
\[ \Delta^1 \xrightarrow{\epsilon} \tx{Fun}(\tx{Cat}_\infty^\tx{pt} , \tx{Cat}_\infty^\tx{pt}) \xrightarrow{U\circ -} \tx{Fun}(\tx{Cat}_\infty^\tx{pt} , \tx{Cat}_\infty).  \]
\end{itemize}
Let $X : \Delta^1 \to \tx{Cat}_\infty^\tx{pt}$ be the inclusion of the zero object, $\Delta^0 \xrightarrow{x}\s{C}$. Uncurrying $\tau$ to a morphism $\Delta^2 \times \tx{Cat}_\infty^\tx{pt} \to \tx{Cat}_\infty$ and restricting along $\Delta^2 \times \Delta^1 \xrightarrow{\tx{id}_{\Delta^2} \times X} \Delta^2 \times \tx{Cat}_\infty^\tx{pt}$, we get a morphism of simplicial sets $g : \Delta^2 \times \Delta^1 \to \tx{Cat}_\infty$, in which the square $g|_{\Delta^{\{0,1 \}} \times \Delta^1}$ strictly commutes, meaning that if $\tx{Cat}_{\infty, 0}^{\Delta}\subseteq \tx{Cat}_\infty^{\Delta}$ is the underlying $1$-category of the simplicial category $\tx{Cat}_\infty^{\Delta}$, then $g|_{\Delta^{\{0,1 \}} \times \Delta^1}$ factors through the inclusion $\tx{N}(\tx{Cat}_{\infty, 0}^{\Delta} )\hookrightarrow \tx{N}( \tx{Cat}_{\infty}^{\Delta} ) = \tx{Cat}_\infty$.  By (\cite{htt}, 3.1.4.4 and 3.1.5.1), the diagram $g$ witnesses the monomorphism 
\begin{equation}\label{sectiondiamondzeroobjmono}\begin{tikzcd}
\Delta^0 \arrow[]{d}{x} \arrow[hook]{r}{\iota_{\Delta^0}} & \Delta^0_\diamond \arrow[]{d}{x_\diamond} \\
\s{C} \arrow[hook]{r}{\iota_{\s{C}}} & \s{C}_\diamond 
\end{tikzcd}\end{equation}
in $\tx{sSet}^{[1]}$ as a split monomorphism in $\tx{Ho}(\tx{sSet})^{[1]}$, where $\tx{sSet}$ is given the Joyal model structure.\footnote{In this proof, we use $\tx{Ho}(\bf{A})$ to denote the homotopy category of an ordinary model category $\bf{A}$.  Remember that $\tx{sSet}$ with the Joyal model category structure \emph{doesn't} extend to a simplicial model category under the standard simplicial enrichment on $\tx{sSet}$.  Thus we need to be a bit careful when writing something like ``$\tx{h}(\tx{sSet})$,'' if we're using the Joyal model structure on $\tx{sSet}$.}  By (\cite{cordierporterhomotopycoherent}, 4.7) and (\cite{htt}, \S 3.1.4), the induced map $\tx{Ho}(\tx{sSet}^{[1]}) \to \tx{Ho}(\tx{sSet})^{[1]}$ (where $\tx{sSet}^{[1]}$ is given the projective Joyal model structure) is an equivalence of categories.  Equivalences of categories reflect split monomorphisms, so (\ref{sectiondiamondzeroobjmono}) is also a split monomorphism in $\tx{Ho}(\tx{sSet}^{[1]})$.

Let $X$, $X'$, and $\phi$ respectively be the objects $\Delta^0 \xrightarrow{x} \s{C}$, $\Delta^0_\diamond \xrightarrow{x_\diamond} \s{C}_\diamond$, and $\Delta^{-1} \xrightarrow{\phi} \Delta^{-1}$ of $\tx{sSet}^{[1]}$, and let $\iota : X\to X'$ in $\tx{sSet}^{[1]}$ denote the natural transformation $\tx{id}_{\tx{sSet}} \Rightarrow (-)_\diamond$ applied to $X$.  By the above, we have a commutative diagram
\[\begin{tikzcd}
X \arrow[]{dd}[swap]{[\iota]} \arrow[]{dr}{[\tx{id}_X]} & \\
 & X \\
 X' \arrow[]{ur}[swap]{\ol{F}} &
\end{tikzcd}\]
in the homotopy category $\tx{Ho}\big( \tx{sSet}^{[1]} \big)$. By Lemma~\ref{cofinprojmodelstr}, the morphisms $\phi \to X$, $\phi \to X'$ and $\iota : X\to X'$ are cofibrations under the projective Joyal model structure, and $X$ is fibrant, so by (\cite{htt}, A.2.3.1) there is a morphism $F : X' \to X$ in $\tx{sSet}^{[1]}$ with $[F] = \ol{F}$ such that the triangle
\[\begin{tikzcd}
X \arrow[]{dd}[swap]{\iota} \arrow[]{dr}{\tx{id}_X} & \\
 & X \\
 X' \arrow[]{ur}[swap]{F} &
\end{tikzcd}\]
in $\tx{sSet}^{[1]}$ commutes, so we have a commutative diagram of simplicial sets
\[\begin{tikzcd}
 \Delta^0 \arrow[]{rr}{x} \arrow[]{dr} \arrow[hook]{dd} & & \s{C} \arrow[]{dr}{\tx{id}_{\s{C}}} \arrow[hook]{dd} & \\
 & \Delta^0 \arrow[crossing over]{rr}[near start]{x} & & \s{C} \\
\Delta^0_\diamond  \arrow[]{ur}[swap]{F(0)} \arrow[]{rr}[swap]{x_\diamond} & & \s{C}_\diamond \arrow[]{ur}[swap]{F(1)} &
\end{tikzcd}\]
and $f := F(1) : \s{C}_\diamond \to \s{C}$ satisfies condition~\ref{sectiondiamondzeroobj2}.
\end{proof}
\begin{definition}\label{cohnulldefn} Let $Y$ be a simplicial set.  We call a morphism $\f{z} : Y_\diamond \to Y$ a \emph{coherent choice of nullhomotopies} if it satisfies the following conditions:
\begin{enumerate}
\item The restriction $\f{z} |_{Y}$ is $\tx{id}_Y$, and
\item If $x := \f{z}(z_Y) \in Y$ is the image of the vertex $z_Y \in Y_\diamond$, then $\f{z}|_{\{ x\}_\diamond }$ is the constant map $\{ x\}_\diamond \to \Delta^0 \xrightarrow{x} Y$.
\end{enumerate}\end{definition}
\begin{rmk}\label{dgnervermk} Let $\s{A}$ be an additive category with a given zero object $0\in \s{A}$, and let $\tx{Ch}(\s{A})$ be the differential graded category of chain complexes in $\s{A}$ (see \cite{ha}, \S 1.3.1).  Then the differential graded nerve $\tx{N}_{\tx{dg}} (\tx{Ch}(\s{A}) )$ of $\tx{Ch}(\s{A})$ is a pointed $\infty$-category, with a zero object $0_\bullet \in \tx{N}_{\tx{dg}} (\tx{Ch}(\s{A}) )$ given by the constant chain complex at $0 \in \s{A}$.

We can explicitly construct a coherent choice of nullhomotopies 
\[\f{z} : \tx{N}_{\tx{dg}} (\tx{Ch}(\s{A}) )_\diamond  \to \tx{N}_{\tx{dg}} (\tx{Ch}(\s{A}) ) \tx{ \; with \; }\f{z}(z_{ \tx{N}_{\tx{dg}}(\tx{Ch}(\s{A}) )}) = 0_\bullet \]
as follows, using (\cite{ha}, 1.3.1.6) for the construction of the differential graded nerve:

\hfill

Suppose $n\geq 0$ and an $n$-simplex $\tilde{\sigma}$ of $\tx{N}_{\tx{dg}} (\tx{Ch}(\s{A}) )_\diamond$,
\[ \tilde{\sigma} = \big[ \Delta^n \xhookleftarrow{\psi} \tx{N}(P) \xrightarrow{\tau} \tx{N}_{\tx{dg}} (\tx{Ch}(\s{A}) )  \big].  \]
If $I = \{ i_- < i_m < i_{m-1} < \dots < i_1 < i_+ \} \subseteq [n]$ is such that $\Delta^I \hookrightarrow \Delta^n$ factors through $\psi$, then there is a unique, nondegenerate $(m+1)$-simplex $\sigma_I : \Delta^{m+1}\to \tx{N}(P)$ such that $\psi \circ \sigma_I = (\Delta^{m+1} \cong \Delta^I \hookrightarrow \Delta^n)$.  Write $( \{Y_{I, j} \}_{0\leq j \leq m+1} , \{ g_{I,J} \})$ for the $(m+1)$-simplex of $\tx{N}_{\tx{dg}} (\tx{Ch}(\s{A}))$ given by $\tau \circ \sigma_I : \Delta^{m+1} \to \tx{N}_{\tx{dg}} (\tx{Ch}(\s{A}))$, where $J$ ranges over all subsets of $[m+1]$ containing at least two distinct elements.  

We define $\f{z}(\tilde{\sigma}) : \Delta^n \to \tx{N}_{\tx{dg}} (\tx{Ch}(\s{A}))$ to be the $n$-simplex $( \{ X_i \}_{0\leq i \leq n} , \{ f_I \} )$, where
\[ X_i := \begin{cases} \tau (p) & \tx{if there exists }p\in P\tx{ with }\psi(p)=i \\ 
0_\bullet & \tx{otherwise,}\end{cases} \]
and for $I = \{ i_- < i_m < i_{m -1} < \dots < i_1 < i_+ \} \subseteq [n]$,
\[ f_I := \begin{cases} g_{I, [m+1]} \in \tx{Map}_{\tx{Ch}(\s{A})} (X_{i_-}, X_{i_+})_m & \tx{if } \Delta^I \hookrightarrow \Delta^n \tx{ factors through }\psi \\
0 \in \tx{Map}_{\tx{Ch}(\s{A})} (X_{i_-}, X_{i_+})_m  & \tx{otherwise}. \end{cases} \]
\end{rmk}

\section{Fiber sequences via coherent choices of nullhomotopies}

In this section, we apply coherent choices of nullhomotopies to give an alternate description of fiber sequences in pointed $\infty$-categories, with this description mimicking the $1$-categorical definition of a kernel of $f: x\to y$ (a morphism in $1$-category $\s{C}$ enriched over pointed sets) as a terminal object of the full subcategory $\{ g\in \s{C}_{/x}\tx{ ; } f\circ g=0 \} \subseteq \s{C}_{/x}$.

\begin{construction} Let $Y$ be a simplicial set with a coherent choice of nullhomotopies $\f{z} : Y_\diamond \to Y$ and let $b\xrightarrow{g} c$ be an edge of $Y$.  We define a simplicial subset $Y_{\perp g} \subseteq Y_{/g}$ as follows:
\begin{itemize}
\item For $n\geq 0$, an $n$-simplex $\ol{\sigma} : \Delta^n \star \Delta^1 \to Y$ of $Y_{/g}$ is in $Y_{\perp g}$ if there exists an $(n+1)$-simplex $\theta : \Delta^{n+1} \to Y_\diamond$ of the form $\big[ \Delta^n \star \Delta^{\{1\}} \hookleftarrow \Delta^n \coprod \Delta^{\{1\}} \to Y \big]$ such that $\f{z}(\theta) = \ol{\sigma}|_{\Delta^n \star \Delta^{\{1\}}}$.
\end{itemize}
\end{construction}
\begin{lemma}[\tx{}\cite{thesis}, 1.5.8]  Let $\s{C}$ be an $\infty$-category with a coherent choice of nullhomotopies $\f{z} : \s{C}_\diamond \to \s{C}$ and let $b\xrightarrow{g} c$ be a morphism in $\s{C}$.  Then the composite map $\s{C}_{\perp g} \hookrightarrow \s{C}_{/g} \to \s{C}$ is a right fibration.  In particular, $ \s{C}_{\perp g}$ is an $\infty$-category.
\end{lemma}
\begin{rmk}[\tx{}\cite{thesis}, 1.5.9]\label{cohperpslicecomparison} Let $\s{C}$ be a (pointed) $\infty$-category with a coherent choice of nullhomotopies $\f{z} : \s{C}_\diamond \to \s{C}$, and let $b\xrightarrow{g} c$ be a morphism in $\s{C}$.  Let $\tilde{g} : \Lambda^2_2 \to \s{C}$ be the right horn given by $\tilde{g}|_{\Delta^{\{0,2\}}} = g$ and $\tilde{g}|_{\Delta^{\{1,2\}}} = \f{z} ( \big[ \Delta^1 \hookleftarrow \Delta^{\{1 \}} \xrightarrow{c} \s{C} \big] )$.  We have a morphism of right fibrations over $\s{C}$
\[\begin{tikzcd}
\s{C}_{\perp g} \arrow[]{dr} \arrow[hook]{rr}{\ell_g} & & \s{C}_{/\tilde{g}} \arrow[]{dl} \\ 
& \s{C} &
\end{tikzcd}\]
where $\ell_g$ is defined as follows:

For an $n$-simplex $\sigma$ of $\s{C}_{\perp g}$ given by $\ol{\sigma} : \Delta^n \star \Delta^2 \to \s{C}$, the $n$-simplex $\ell_g (\sigma )$ of $\s{C}_{/\tilde{g}}$ is given by the composite
\[ \Delta^n \star \Lambda^2_2 \cong \Delta^n \star \Delta^{\{ 0,2 \}} \coprod_{\Delta^n \star \Delta^{\{ 2 \}}} \Delta^n \star \Delta^{\{ 1,2 \}} \xrightarrow{( \ol{\sigma} , \sigma' )} \s{C}, \]
where $\sigma' : \Delta^n \star \Delta^{\{ 1,2 \}} \to \s{C}$ is the $(n+2)$-simplex of $\s{C}$ given by 
\[ \f{z}\left( \big[ \Delta^n \star \Delta^{\{ 1,2 \}} \hookleftarrow \Delta^n \coprod \Delta^{\{ 2 \}} \xrightarrow{( \ol{\sigma}|_{\Delta^n} , c )} \s{C} \big]  \right). \]
Note that an $n$-simplex $\ol{\sigma} : \Delta^n \star \Lambda^2_2$ of $\s{C}_{/\tilde{g}}$ is in the image of the monomorphism $\ell_g$ if and only if $\ol{\sigma}|_{\Delta^n \star \Delta^{1,2}}$ is of the form $\f{z}\left( \big[ \Delta^n \star \Delta^{\{ 1,2 \}} \hookleftarrow \Delta^n \coprod \Delta^{\{ 2 \}} \to \s{C} \big]  \right)$.

The map $\ell_g$ factors through a simplicial subset of $\s{C}_{/\tilde{g}}$, which will be convenient to define:  Let $\s{C}_{/\tilde{g}}' \subseteq \s{C}_{/\tilde{g}}$ be the simplicial subset given as follows:
\begin{itemize}
\item For $n\geq 0$, an $n$-simplex $\ol{\sigma} : \Delta^n \star \Lambda^2_2 \to \s{C}$ of $\s{C}_{/\tilde{g}}$ is in $\s{C}_{/\tilde{g}}' $ if there exists an $n+1$-simplex of $\s{C}_\diamond$ of the form 
$\theta = \big[  \Delta^n \star \Delta^{\{1 \}} \hookleftarrow  \Delta^n \star \Delta^{-1}  \to \s{C}  \big]$
such that $\ol{\sigma} |_{\Delta^n \star \Delta^{\{1 \}} } = \f{z}(\theta )$.
\end{itemize}
Note that the composite $ \s{C}_{/\tilde{g}}' \hookrightarrow \s{C}_{/\tilde{g}} \to \s{C}$ is a right fibration.
\end{rmk}
\begin{proposition}\label{cohnultriangleprop} Let $\s{C}$ be a (pointed) $\infty$-category with a coherent choice of nullhomotopies $\f{z} : \s{C}_\diamond \to \s{C}$, and let $b\xrightarrow{g} c$ be a morphism in $\s{C}$.  Then the monomorphism $\ell_g : \s{C}_{\perp g} \to \s{C}_{/ \tilde{g}}$ of Remark~\ref{cohperpslicecomparison} is a categorical equivalence.  

Thus, by (\cite{kerodon}, 02HS, 02MY and 02HN), an object $x \in \s{C}_{\perp g}$ is terminal if and only if $\ell_g (x) \in \s{C}_{/ \tilde{g}}$ is terminal (i.e., $\ell_g (x)$ is a limit diagram in $\s{C}$), and an object $y \in  \s{C}_{/ \tilde{g}}$ is terminal if and only if there exists a terminal object $x \in \s{C}_{\perp g}$ such that $y$ is isomorphic to $\ell_g (x)$ in the $\infty$-category $\s{C}_{/ \tilde{g}}$.
\end{proposition}
\begin{proof} Using the notation of Remark ~\ref{cohperpslicecomparison}, we have morphisms of right fibrations over $\s{C}$
\[\begin{tikzcd}
\s{C}_{\perp g} \arrow[bend right]{dr} \arrow[hook, bend left]{rr}{\ell_g} \arrow[hook]{r}[swap]{u_g} &  \s{C}_{/\tilde{g}} \arrow[hook]{r}[swap]{v_g} \arrow[]{d} & \s{C}_{/\tilde{g}} \arrow[bend left]{dl} \\ 
& \s{C} & .
\end{tikzcd}\]
Thus to show that $\ell_g$ is a categorical equivalence, it suffices to show that for each $x \in \s{C}$, the induced morphisms of fibers
\[ ( u_g)_x : \{ x \} \times_\s{C} \s{C}_{\perp g} \hookrightarrow  \{ x \} \times_\s{C} \s{C}_{/ \tilde{g}}' \tx{ \; and \; } (v_g)_x : \{ x \} \times_\s{C} \s{C}_{/ \tilde{g}} \hookrightarrow  \{ x \} \times_\s{C} \s{C}_{/ \tilde{g}}  \]
are both weak equivalences (of Kan complexes).  Since $(v_g)_x$ is a monomorphism of Kan complexes, it suffices to show that $\{ x \} \times_\s{C} \s{C}_{/ \tilde{g}}' = (\s{C}_{/ \tilde{g}}')_x $ is a strong deformation retract of $\{ x \} \times_\s{C} \s{C}_{/ \tilde{g}}  = (\s{C}_{/ \tilde{g}})_x$.  We'll construct a deformation retraction 
\[ \Theta = \Theta^x :  (\s{C}_{/ \tilde{g}})_x \to \tx{Fun}(\Delta^1 , (\s{C}_{/ \tilde{g}})_x)  \]
by simplicial induction on $(\s{C}_{/ \tilde{g}})_x$:

\hfill

Let $A$ be the collection of non-degenerate simplicies of $(\s{C}_{/ \tilde{g}})_x$, and for $\sigma \in A$ let $d(\sigma )$ be the dimension of $\sigma$ (considered as a simplex of  $(\s{C}_{/ \tilde{g}})_x$).  Pick a well-ordering on $A$ such that $\sigma <\tau$ if $d(\sigma) < d(\tau)$, and let $\alpha$ be the order type of $A$, so we have an order-preserving bijection $\beta \mapsto \sigma_\beta$, where $\beta$ ranges over the ordinals $\{\beta ; \beta <\alpha \}$.  For each $\beta\leq \alpha$, let $\s{D}_\beta \subseteq (\s{C}_{/ \tilde{g}})_x$ be the largest simplicial subset whose nondegenerate simplicies are of the form $\sigma_{\gamma}$ for some $\gamma < \beta$.  We will show that for each $\beta \leq \alpha$, there exists a family of maps $\{ \Theta_{\beta'} : \s{D}_{\beta'} \to \tx{Fun}(\Delta^1 , (\s{C}_{/ \tilde{g}})_x) \}_{\beta' \leq \beta}$ such that the following two conditions hold: 
\begin{enumerate}
\item If $\beta'' \leq \beta' \leq \beta$ then $\Theta_{\beta'} |_{\s{D}_{\beta''}} = \Theta_{\beta''}$, and \item For each $\beta' \leq \beta$, if we consider $\Theta_{\beta'}$ as a map $\Delta^1 \times \s{D}_{\beta'} \to (\s{C}_{/\tilde{g}})_x$ then $\Theta_{\beta'} |_{\Delta^{\{0\}} \times \s{D}_{\beta'}}$ factors through $v_g$ and the following diagram commutes:
\begin{equation}\label{cohnultrianglepropdiag1} \begin{tikzcd} &  \Delta^{\{1\}} \times \s{D}_{\beta'} \arrow[hook]{d} \arrow[hook]{dr} & \\
\Delta^1 \times ( \s{D}_{\beta'} \cap (\s{C}_{/\tilde{g}}')_x) \arrow[hook]{r} \arrow[]{d}{\tx{pr}_2} & \Delta^1 \times \s{D}_{\beta'} \arrow[]{r}[swap]{\Theta_{\beta'}} & (\s{C}_{/\tilde{g}})_x \\
\s{D}_{\beta'} \cap (\s{C}_{/\tilde{g}}')_x \arrow[hook, bend right=20]{urr} & &
\end{tikzcd}\end{equation}
\end{enumerate}
\begin{itemize}
\item 
At $\beta =0$, we have $\s{D}_0 = \Delta^{-1}$; pick $\Theta_0$ to be the unique map $\Delta^{-1} \to \tx{Fun}(\Delta^1 , (\s{C}_{/ \tilde{g}})_x)$.  

\item 
Next, suppose $\beta = \varpi +1 \leq \alpha$ is a successor and that the induction hypothesis holds for $\beta' < \beta$.  Let $n=d(\sigma_{\varpi})$, and let the family $\{\Theta_{\beta'} \}_{\beta' \leq \varpi}$ satisfy the above conditions. We have a pushout of simplicial sets
\[\begin{tikzcd}
\partial \Delta^n \arrow[hook]{d} \arrow[]{r}{\gamma} & \s{D}_{\varpi} \arrow[hook]{d} \\
\Delta^n \arrow[]{r}{\sigma_{\varpi}} & \s{D}_{\varpi +1},
\end{tikzcd}\]
so to define $\Theta_{\beta}$ we need to define an $n$-simplex $\theta_{\varpi} : \Delta^n \to \tx{Fun}(\Delta^1 , (\s{C}_{/ \tilde{g}})_x)$ such that $\theta_{\varpi} |_{\partial \Delta^n} = \Theta_\varpi \circ \gamma$, $\tx{ev}_1 (\theta_{\varpi}) = \sigma_\varpi$, $\tx{ev}_0 (\theta_{\varpi})$ is an $n$-simplex of $(\s{C}'_{/ \tilde{g}})_x$, and such that $\theta_{\varpi}$ is the constant homotopy at $\sigma_\varpi$ if $\sigma_\varpi$ is already an $n$-simplex of $(\s{C}_{/\tilde{g}}')_x$.

If $\sigma_\varpi$ is an $n$-simplex of $(\s{C}_{/\tilde{g}}')_x$, let $\theta_{\varpi}$ be the currying of the composite
\[ \Delta^1 \times \Delta^n \xrightarrow{\tx{pr}_2} \Delta^n \xrightarrow{\sigma_\varpi} (\s{C}_{/\tilde{g}})_x . \]

If  $\sigma_\varpi$ is not an $n$-simplex of $(\s{C}_{/\tilde{g}}')_x$, suppose $\sigma_{\varpi}$ is given by the map $\chi : \Delta^n \star \Lambda^2_2 \to \s{C}$, and suppose $\Theta_\varpi \circ \gamma$ is given by $\rho : (\Delta^1 \times \partial \Delta^n) \star \Lambda^2_2 \to \s{C}$.  We want to construct a map $\tilde{\chi} : (\Delta^1 \times \Delta^n) \star \Lambda^2_2 \to \s{C}$ such that $\tilde{\chi}|_ {(\Delta^{1} \times \Delta^n )\star \Delta^{-1}}$ is the constant map to $x$, $\tilde{\chi}|_ {(\Delta^{\{1\}} \times \Delta^n )\star \Lambda^2_2} = \chi$, $\tilde{\chi}|_{(\Delta^{1} \times \partial \Delta^n )\star \Lambda^2_2} = \rho$, and such that $\tilde{\chi}|_ {(\Delta^{\{0\}} \times \Delta^n )\star \Lambda^2_2}$ corresponds to an $n$-simplex of $(\s{C}_{/\tilde{g}}')_x$, i.e. the restriction $\tilde{\chi}|_ {(\Delta^{\{0\}} \times \Delta^n )\star \Delta^{\{1\}}}$ is of the form $\f{z}\left(  \big[ \Delta^n \star \Delta^{\{1\}} \hookleftarrow \Delta^n \star \Delta^{-1} \to \s{C} \big] \right)$.

Since $\chi$ takes the vertex $1 \in \Lambda^2_2$ to $\f{z}(z_\s{C})$, which is a terminal object of $\s{C}$, (\cite{thesis}, 1.3.6) implies we can find a map $r : (\Delta^1 \times \Delta^n) \star \Delta^{\{1 \}} \to \s{C}$ such that $r|_{ (\Delta^1 \times \Delta^n) \star \Delta^{-1} }$ is the constant map to $x$, $r|_{(\Delta^{\{0 \}} \times \Delta^n) \star \Delta^{\{1 \}} } = \chi |_{\Delta^n \star \Delta^{\{1\}}} $, $r|_{ (\Delta^1 \times \partial \Delta^n) \star \Delta^{\{1 \}} } = \rho |_{(\Delta^1 \times \partial \Delta^n) \star \Delta^{\{1 \}}}$,
and $r|_{(\Delta^{\{0 \}} \times \Delta^n) \star \Delta^{\{1 \}} }$ is of the form $\f{z}\left(  \big[ \Delta^n \star \Delta^{\{1\}} \hookleftarrow \Delta^n \star \Delta^{-1} \to \s{C} \big] \right)$.  

Thus we are left with an extension problem
\[ \begin{tikzcd}
\big( (\Delta^1 \times \Delta^n) \star \Delta^{\{1 \}} \big) \cup \big( (\Delta^{\{1\}} \times \Delta^n )\star \Lambda^2_2 \big) \cup \big( (\Delta^{1} \times \partial \Delta^n )\star \Lambda^2_2 \big) \arrow[]{r}{(r, \chi , \rho)} \arrow[equal]{d} & \s{C} \\
\big( (\Delta^1 \times \Delta^n) \star \Delta^{\{1 \}} \big) \cup \big( ( \Delta^{\{1\}} \times  \Delta^n  \cup \Delta^{1} \times \partial \Delta^n )\star \Lambda^2_2 \big) \arrow[hook]{d} & \\
(\Delta^1 \times \Delta^n) \star \Lambda^2_2 \arrow[dashed, bend right =30]{uur}[swap]{\tilde{\chi}}, &
\end{tikzcd}\]
which is solvable because the left vertical is inner anodyne by (\cite{htt}, 2.1.2.6 and 2.1.2.3).  The map $\tilde{\chi}$ then corresponds an $n$-simplex $\theta_{\varpi} : \Delta^n \to \tx{Fun}(\Delta^1 , (\s{C}_{/ \tilde{g}})_x)$ with the desired properties. 
\item Suppose $\beta \leq \alpha$ is a limit ordinal, and that the induction hypothesis holds for all $\beta' < \beta$, so by Zorn's lemma we can find a family $\{ \Theta_{\beta'} : \s{D}_{\beta'} \to \tx{Fun}(\Delta^1 , (\s{C}_{/ \tilde{g}})_x) \}_{\beta' < \beta}$ such that if $\beta'' \leq \beta' < \beta$ then $\Theta_{\beta'} |_{\s{D}_{\beta''}} = \Theta_{\beta''}$, and if $\beta' < \beta$ then diagram (\ref{cohnultrianglepropdiag1}) commutes.  Then $\s{D}_\beta \cong \colim_{\beta' < \beta} \s{D}_{\beta'}$, and letting $\Theta_\beta$ be the map $\s{D}_\beta \to \tx{Fun}(\Delta^1 , (\s{C}_{/ \tilde{g}})_x)$ induced by the family $\{ \Theta_{\beta'} \}_{\beta' < \beta}$ gives the desired family $\{ \Theta_{\beta'} \}_{\beta' \leq \beta}$.
\end{itemize}
Thus there exists a family $ \{ \Theta_\beta \}_{\beta \leq \alpha}$ satisfying the above properties, and picking $\Theta^x = \Theta_\alpha$ gives the desired strong deformation retraction of $(\s{C}_{/ \tilde{g}})_x$ to $(\s{C}'_{/ \tilde{g}})_x$.  This completes the proof that $v_g$ is a categorical equivalence.

\hfill

Next we'll construct a strong deformation retraction 
\[ \Phi = \Phi^x :  (\s{C}'_{/ \tilde{g}})_x \to \tx{Fun}(\Delta^1 , (\s{C}'_{/ \tilde{g}})_x)  \]
to the image of $(u_g)_x$ by a similar simplicial induction on $(\s{C}'_{/ \tilde{g}})_x$, with the only real difference being a different extension problem in the successor step; for convenience, we will write $\s{C}_{\perp g}$ for the image of $\s{C}_{\perp g}$ in $\s{C}'_{/ \tilde{g}}$ under the monomorphism $u_g$.

\hfill

Let $B$ be the collection of non-degenerate simplicies of $(\s{C}'_{/ \tilde{g}})_x$, and for $\sigma \in B$ let $d(\sigma )$ be the dimension of $\sigma$ (considered as a simplex of  $(\s{C}'_{/ \tilde{g}})_x$).  Pick a well-ordering on $B$ such that $\sigma <\tau$ if $d(\sigma) < d(\tau)$, and let $\delta$ be the order type of $B$, so we have an order-preserving bijection $\beta \mapsto \sigma_\beta$, where $\beta$ ranges over the ordinals $\{\beta ; \beta <\delta \}$.  For each $\beta\leq \delta$, let $\s{E}_\beta \subseteq (\s{C}'_{/ \tilde{g}})_x$ be the largest simplicial subset whose nondegenerate simplicies are of the form $\sigma_{\gamma}$ for some $\gamma < \beta$.  We will show that for each $\beta \leq \delta$, there exists a family of maps $\{ \Phi_{\beta'} : \s{E}_{\beta'} \to \tx{Fun}(\Delta^1 , (\s{C}'_{/ \tilde{g}})_x) \}_{\beta' \leq \beta}$ such that the following two conditions hold: 
\begin{enumerate}[label=\arabic*'.]
\item\label{cohnultriangleprop1prime} If $\beta'' \leq \beta' \leq \beta$ then $\Phi_{\beta'} |_{\s{E}_{\beta''}} = \Phi_{\beta''}$, and 
\item\label{cohnultriangleprop2prime} For each $\beta' \leq \beta$, if we consider $\Phi_{\beta'}$ as a map $\Delta^1 \times \s{E}_{\beta'} \to (\s{C}'_{/\tilde{g}})_x$ then $\Phi_{\beta'} |_{\Delta^{\{0\}} \times \s{E}_{\beta'}}$ factors through the simplicial subset $(\s{C}_{\perp g})_x \subseteq  (\s{C}'_{/\tilde{g}})_x $ and the following diagram commutes:
\begin{equation}\label{cohnultrianglepropdiag2} \begin{tikzcd} &  \Delta^{\{1\}} \times \s{E}_{\beta'} \arrow[hook]{d} \arrow[hook]{dr} & \\
\Delta^1 \times ( \s{E}_{\beta'} \cap (\s{C}_{\perp g})_x) \arrow[hook]{r} \arrow[]{d}{\tx{pr}_2} & \Delta^1 \times \s{E}_{\beta'} \arrow[]{r}[swap]{\Phi_{\beta'}} & (\s{C}'_{/\tilde{g}})_x \\
\s{E}_{\beta'} \cap (\s{C}_{\perp g})_x \arrow[hook, bend right=20]{urr} & &
\end{tikzcd}\end{equation}
\end{enumerate}
\begin{itemize}
\item At $\beta =0$, we can pick $\Phi_0$ to be the unique map $\Delta^{-1} = \s{E}_0 \to  \tx{Fun}(\Delta^1 , (\s{C}'_{/ \tilde{g}})_x)$.
\item Next, suppose $\beta = \varpi +1 \leq \delta$ is a successor and that the induction hypothesis holds for all $\beta' < \beta$.  Let $n=d(\sigma_{\varpi})$, and let the family $\{\Phi_{\beta'} \}_{\beta' \leq \varpi}$ satisfy the above conditions. We have a pushout of simplicial sets
\[\begin{tikzcd}
\partial \Delta^n \arrow[hook]{d} \arrow[]{r}{\gamma} & \s{E}_{\varpi} \arrow[hook]{d} \\
\Delta^n \arrow[]{r}{\sigma_{\varpi}} & \s{E}_{\varpi +1},
\end{tikzcd}\]
so to define $\Phi_{\beta}$ it suffices to define an $n$-simplex $\phi_{\varpi} : \Delta^n \to \tx{Fun}(\Delta^1 , (\s{C}'_{/ \tilde{g}})_x)$ such that $\phi_{\varpi} |_{\partial \Delta^n} = \Phi_\varpi \circ \gamma$, $\tx{ev}_1 (\phi_{\varpi}) = \sigma_\varpi$, $\tx{ev}_0 (\phi_{\varpi})$ is an $n$-simplex of $(\s{C}_{\perp g})_x$, and such that $\phi_{\varpi}$ is the constant homotopy at $\sigma_\varpi$ if $\sigma_\varpi$ is already an $n$-simplex of $(\s{C}_{\perp g})_x$.

If $\sigma_\varpi$ is an $n$-simplex of $(\s{C}_{\perp g})_x$, let $\phi_{\varpi}$ be the currying of the composite
\[ \Delta^1 \times \Delta^n \xrightarrow{\tx{pr}_2} \Delta^n \xrightarrow{\sigma_\varpi} (\s{C}'_{/\tilde{g}})_x . \]

If  $\sigma_\varpi$ is not an $n$-simplex of $(\s{C}_{\perp g})_x$, suppose $\sigma_{\varpi}$ is given by the map $\chi : \Delta^n \star \Lambda^2_2 \to \s{C}$, and suppose $\Phi_\varpi \circ \gamma$ is given by $\rho : (\Delta^1 \times \partial \Delta^n) \star \Lambda^2_2 \to \s{C}$.  We want to construct a map $\tilde{\chi} : (\Delta^1 \times \Delta^n) \star \Lambda^2_2 \to \s{C}$ such that $\tilde{\chi}|_ {(\Delta^{\{1\}} \times \Delta^n )\star \Lambda^2_2} = \chi$, $\tilde{\chi}|_{(\Delta^{1} \times \partial \Delta^n )\star \Lambda^2_2} = \rho$, and $\tilde{\chi}|_ {(\Delta^{\{0\}} \times \Delta^n )\star \Delta^{\{ 1,2 \}}}$ is of the form $\f{z}  \left( \big[ \Delta^n \star \Delta^{\{1,2 \}} \hookleftarrow \Delta^n \coprod \Delta^{\{2 \}} \to \s{C}   \big]  \right)$. 

Since $\sigma_\varpi$ is a simplex of $(\s{C}'_{/ \tilde{g}})_x$, the restriction $\chi |_{\Delta^n \star \Delta^{-1}}$ is the constant map to $x$, $\chi |_{\Delta^n \star \Delta^{\{1\}}}$ is of the form $\f{z}\left( \big[ \Delta^n \star \Delta^{\{1\}} \hookleftarrow \Delta^n \star \Delta^{-1} \to \s{C} \big] \right)$, and $\chi |_{\Delta^{-1} \star \Delta^{\{1,2 \}}}$ is given by $\f{z}\left( \big[ \Delta^{\{1,2\}} \hookleftarrow \Delta^{\{ 2 \}}  \xrightarrow{c} \s{C} \big] \right)$.  Since $\Phi_\varpi$ is a homotopy in $(\s{C}'_{/ \tilde{g}})_x$, we have that for each nondegenerate simplex $\mu : \Delta^k \hookrightarrow \Delta^1 \times \partial \Delta^n$ of $\Delta^1 \times \partial \Delta^n$, the restriction $\rho |_{\mu \star \Delta^{\{1 \}}}$ is of the form $\f{z} \left( \big[ \Delta^k \star \Delta^{\{1 \}} \hookleftarrow \Delta^k \xrightarrow{\ul{x}} \s{C} \big] \right)$, where $\ul{x}$ denotes the constant map $ \Delta^k \to \Delta^0 \xrightarrow{x} \s{C}$.  By the induction hypothesis, for each nondegenerate simplex $\nu : \Delta^m \hookrightarrow \Delta^{\{0\}} \times \partial \Delta^n$ of $\Delta^{\{0\}} \times \partial \Delta^n$, the restriction  $\rho |_{\nu \star \Delta^{\{1,2 \}}}$ is of the form \[\f{z} \big( \big[ \Delta^m \star \Delta^{\{1,2 \}} \hookleftarrow \Delta^k \coprod \Delta^{\{2\}} \xrightarrow{(\ul{x} , c)} \s{C} \big] \big).\]

Unwinding what this means in terms of convexity, this imples there is a map \[r : (\Delta^1 \times \Delta^n) \star \Delta^{\{1, 2 \}} \to \s{C}\] where $r |_{(\Delta^1 \times \Delta^n) \star \Delta^{-1}}$ is the constant map to $x$, $r |_{(\Delta^1 \times \partial \Delta^n) \star \Delta^{\{1,2\}}}= \rho |_{(\Delta^1 \times \partial \Delta^n) \star \Delta^{\{1,2 \}}}$, $r|_{\Delta^{-1} \star \Delta^{\{1,2\}}} = \f{z}\left( \big[ \Delta^{\{1,2\}} \hookleftarrow \Delta^{\{ 2 \}}  \xrightarrow{c} \s{C} \big] \right)$, $r |_{(\Delta^{\{1 \}} \times \Delta^n) \star \Delta^{\{1, 2 \}} } = \chi |_{\Delta^n \star \Delta^{\{1, 2 \}}}$, and such that $r |_{(\Delta^{\{0 \}} \times \Delta^n) \star \Delta^{\{1, 2 \}} }$ is of the form $\f{z}\left( \big[ \Delta^n \star \Delta^{\{1, 2 \}} \hookleftarrow \Delta^n \coprod \Delta^{\{2 \}} \to \s{C} \big]  \right)$.

Thus we are left with an extension problem
\[ \begin{tikzcd}
\big( (\Delta^1 \times \Delta^n) \star \Delta^{\{1,2 \}} \big) \cup \big( (\Delta^{\{1\}} \times \Delta^n )\star \Lambda^2_2 \big) \cup \big( (\Delta^{1} \times \partial \Delta^n )\star \Lambda^2_2 \big) \arrow[]{r}{(r, \chi , \rho)} \arrow[equal]{d} & \s{C} \\
\big( (\Delta^1 \times \Delta^n) \star \Delta^{\{1,2 \}} \big) \cup \big( ( \Delta^{\{1\}} \times  \Delta^n  \cup \Delta^{1} \times \partial \Delta^n )\star \Lambda^2_2 \big) \arrow[hook]{d} & \\
(\Delta^1 \times \Delta^n) \star \Lambda^2_2 \arrow[dashed, bend right =30]{uur}[swap]{\tilde{\chi}}, &
\end{tikzcd}\]
which is solvable because the left vertical is inner anodyne by (\cite{htt}, 2.1.2.6 and 2.1.2.3).  The map $\tilde{\chi}$ then corresponds to an $n$-simplex $\phi_{\varpi} : \Delta^n \to \tx{Fun}(\Delta^1 , (\s{C}'_{/ \tilde{g}})_x)$ with the desired properties.
\item Suppose $\beta \leq \delta$ is a limit ordinal, and that the induction hypothesis holds for all $\beta' < \beta$.  Then the induction hypothesis holds for $\beta$, by the same argument as in the case of $v_g$.
\end{itemize}
Thus there exists a family $\{\Phi_\beta \}_{\beta \leq \alpha}$ with properties~\ref{cohnultriangleprop1prime} and~\ref{cohnultriangleprop2prime}, and picking $\Phi^x = \Phi_\alpha$ gives the desired strong deformation retraction.  Thus $u_g$ is a categorical equivalence, so $\ell_g$ is a categorical equivalence.
\end{proof}
\begin{rmk}\footnote{Here we are using a version of the argument from (\cite{ha}, 1.1.1.7).} Suppose $\s{C}$ is a pointed $\infty$-category with finite limits, let $\f{z} : \s{C}_\diamond \to \s{C}$ be a coherent choice of nullhomotopies, and consider the sequence of restrictions
\[ \tx{Fun}(\{ x \}\star \Lambda^2_2 ,  \s{C}) \to \tx{Fun}( \Lambda^2_2 ,  \s{C}) \to \tx{Fun}( \Delta^{\{0,2 \}} ,  \s{C}). \]
If $\tx{Fun}'( \Lambda^2_2 ,  \s{C}) \subseteq \tx{Fun}( \Lambda^2_2 ,  \s{C})$ is the full subcategory spanned by the functors $F : \Lambda^2_2 \to \s{C}$ such that $F$ is a left Kan extension of $F|_{\Delta^{\{0,2 \}}}$ (i.e., such that $F(1)$ is a zero object of $\s{C}$), then the restriction map 
\[ \tx{Fun}'( \Lambda^2_2 ,  \s{C}) \to  \tx{Fun}( \Delta^{\{0,2 \}} ,  \s{C}) \]
is a trivial Kan fibration by (\cite{htt}, 4.3.2.15).  Likewise, if $\tx{Fun}''(\{ x \}\star \Lambda^2_2 ,  \s{C}) \subseteq \tx{Fun}(\{ x \}\star \Lambda^2_2 ,  \s{C})$ is the full subcategory spanned by the functors $G : \{ x\} \star \Lambda^2_2 \to \s{C}$ such that $G$ is a right Kan extension of $G|_{\Lambda^2_2}$ (i.e., such that $G$ is a limit diagram), then the restriction 
\[ \tx{Fun}''(\{ x \}\star \Lambda^2_2 ,  \s{C}) \to \tx{Fun}( \Lambda^2_2 ,  \s{C}) \]
is a trivial Kan fibration.  In particular, if $\s{E}\subseteq \tx{Fun}(\{ x \}\star \Lambda^2_2 ,  \s{C})$ is the full subcategory spanned by the functors $\{ x \}\star \Lambda^2_2 \to  \s{C}$ which are fiber sequences, then the restriction maps $\s{E} \to \tx{Fun}'( \Lambda^2_2 ,  \s{C})$ and $\s{E} \to \tx{Fun}( \Delta^{\{0,2 \}} ,  \s{C})$ are both trivial Kan fibrations.

Let $g : b\to c$ be a morphism in $\s{C}$, and let $\tilde{g} : \Lambda^2_2 \to \s{C}$ be as in Remark~\ref{cohperpslicecomparison}, i.e. $\tilde{g}|_{\Delta^{\{0,2 \}}} =g$, and $\tilde{g}|_{\Delta^{\{1,2\}}} = \f{z} ( \big[ \Delta^1 \hookleftarrow \Delta^{\{1 \}} \xrightarrow{c} \s{C} \big] )$.  For an $\infty$-category $\s{D}$, let $\s{D}^\tx{fin} \subseteq \s{D}$ denote the full subcategory spanned by the final objects of $\s{D}$.  By (\cite{kerodon}, 04J2), we have a categorical equivalence of $\infty$-categories
\[ \theta : \s{C}_{/\tilde{g}} \to \tx{Fun}(\{ x \}\star \Lambda^2_2 ,  \s{C})\times_{\tx{Fun}(\Lambda^2_2 ,  \s{C})} \{ \tilde{g} \}, \]
which restricts to a weak equivalence of contractible Kan complexes
\[ \theta |_{(\s{C}_{/\tilde{g}})^\tx{fin}} : (\s{C}_{/\tilde{g}})^\tx{fin} \xrightarrow{\sim} \Big( \tx{Fun}(\{ x \}\star \Lambda^2_2 ,  \s{C})\times_{\tx{Fun}(\Lambda^2_2 ,  \s{C})} \{ \tilde{g} \} \Big)^\tx{fin} \simeq \s{E} \times_{\tx{Fun}(\Lambda^2_2 ,  \s{C})} \{ \tilde{g} \}. \]
By Proposition~\ref{cohnultriangleprop}, we have a commutative diagram of contractible Kan complexes
\[ \begin{tikzcd}
(\s{C}_{\perp g})^\tx{fin} \arrow[]{r}{\sim} & (\s{C}_{/\tilde{g}})^\tx{fin} \arrow[]{r}{\sim} &  \s{E}  \times_{ \tx{Fun}' (\Lambda^2_2 , \s{C}) } \{ \tilde{g} \} \arrow[]{r}{\sim} \arrow[two heads]{d}{\rotatebox[origin=c]{90}{$\sim$}} & \s{E} \times_{\tx{Fun}(\Delta^{\{0,2 \}} , \s{C})}   \{ g \} \arrow[two heads]{d}{\rotatebox[origin=c]{90}{$\sim$}} \\
& & \{ \tilde{g} \} \arrow[]{r}{\sim} & \tx{Fun}' (\Lambda^2_2 , \s{C})  \times_{\tx{Fun}(\Delta^{\{0,2 \}} , \s{C})} \{ g \} ,
\end{tikzcd}\]
and an object $y \in \tx{Fun}(\{x\} \star \Lambda^2_2 , \s{C})   \times_{\tx{Fun}(\Delta^{\{0,2 \}} , \s{C})} \{ g \}$ is a fiber sequence in $\s{C}$ if and only if there is a terminal object $y' \in \s{C}_{\perp g}$ such that $y$ is isomorphic to $\ell_g (y')$ in the $\infty$-category $\tx{Fun}(\{x\} \star \Lambda^2_2 , \s{C})   \times_{\tx{Fun}(\Delta^{\{0,2 \}} , \s{C})} \{ g \}$.
\end{rmk}


\begin{references*} \normalsize 
\bibitem{ben} Antieau, Benjamin. ``Math 521 - derived algebraic geometry.'' 2020, \url{https://antieau.github.io/assets/pdf/9weeks_0.1.pdf}.
\bibitem{beke} Beke, Tibor. ``Sheafifiable Homotopy Model Categories, II.” \emph{Journal of Pure and Applied Algebra}, vol. 164, no. 3, Nov. 2001, pp. 307–324.
\bibitem{bergner} Bergner, Julia. \emph{The Homotopy Theory of $(\infty , 1)$-Categories}.  Cambridge University Press, 2018.
\bibitem{thesis} Brown, Thomas Rush.  ``Coherent Choices of Nullhomotopies in Pointed Infinity-Categories.'' 2026, \url{https://escholarship.org/uc/item/2bx4b3x3}.
\bibitem{stevensoncounter} Campbell, Alexander. ``A Counterexample in Quasi-Category Theory.” \emph{Proceedings of the American Mathematical Society}, vol. 148, no. 1, 9 July 2019, pp. 37–40.
\bibitem{cisinski} Cisinski, Denis-Charles. \emph{Higher Categories and Homotopical Algebra}. Cambridge University Press, 2019.
\bibitem{cordierporterhomotopycoherent} Cordier, Jean-Marc, and Timothy Porter. ``Vogt’s Theorem on Categories of Homotopy Coherent Diagrams.” \emph{Mathematical Proceedings of the Cambridge Philosophical Society}, vol. 100, no. 1, 1 July 1986, pp. 65–90.
\bibitem{duggerprimer} Dugger, Daniel. ``A Primer on Homotopy Colimits.” Preprint, 2008, \url{https://pages.uoregon.edu/ddugger/hocolim.pdf}.
\bibitem{dsrigid} Dugger, Daniel, and David I Spivak. ``Rigidification of Quasi-Categories.” \emph{Algebraic \& Geometric Topology}, vol. 11, no. 1, 7 Jan. 2011, pp. 225–261.
\bibitem{univloopmachine} Gepner, D., M. Groth, and T. Nikolaus. ``Universality of Multiplicative Infinite Loop Space Machines.” \emph{Algebraic \& Geometric Topology}, vol. 15, no. 6, 31 Dec. 2015, pp. 3107–3153.
\bibitem{hirsch} Hirschhorn, Philip S. \emph{Model Categories and Their Localizations}. American Mathematical Society, 24 Aug. 2009.
\bibitem{hovey} Hovey, Mark. \emph{Model Categories}. American Mathematical Society, 2007.
\bibitem{bravenew1} Khan, Adeel.  ``Brave new motivic homotopy theory, I.'' Preprint, 2017, \url{https://www.preschema.com/documents/bravenew1.pdf}.
\bibitem{htt} Lurie, Jacob. \emph{Higher Topos Theory}. Princeton, N.J., Princeton University Press, 2009.
\bibitem{ha} Lurie, Jacob. \emph{Higher Algebra}. Harvard Mathematics Department, Sept. 18, 2017, \url{https://people.math.harvard.edu/~lurie/papers/HA.pdf}.
\bibitem{kerodon} Lurie, Jacob. ``Kerodon: an online resource for homotopy-coherent mathematics.” 2026, \url{www.kerodon.net}. 
\bibitem{rezk} Rezk, Charles. ``Introduction to quasicategories.'' Nov. 29, 2022, \url{https://rezk.web.illinois.edu/quasicats.pdf}.
\bibitem{stevenson} Stevenson, Danny. ``Stability for Inner Fibrations Revisited.” \emph{Theory and Applications of Categories}, vol. 33, 2018, pp. 523–536.
\end{references*}

\end{document}